\documentclass{article}[12]

\usepackage{color}

\usepackage{amsmath,amsfonts,amssymb,latexsym}

\begin{document}

\Large 
\begin{center}

{\bf Why did the Germans arrest and release Emile Borel in 1941?}
 
\vskip 3mm 
\normalsize 

Laurent \textsc{Mazliak}\footnote{Laboratoire de Probabilit\'es et Mod\`eles al\'eatoires, {Universit\'e Pierre et Marie Curie, Paris, 
France}. laurent.mazliak@upmc.fr}
 and Glenn \textsc{Shafer}\footnote{Rutgers Business School, Newark, New Jersey, and Royal Holloway University of London.  gshafer@rutgers.edu}

\today

\end{center}

\vskip 3mm 

\small

\normalsize

\begin{center} \bf Abstract \end{center} \rm \vskip 2mm

The Germans occupying Paris arrested Emile Borel and three other members of the \textit{Acad\'emie des Sciences} in October 1941 and released them about five weeks later.  Why?  We examine some relevant German and French archives and other sources and propose some hypotheses.  In the process, we review how the Occupation was structured and how it dealt with French higher education and some French mathematicians.


{\bf MSC-Code} : 01A60


\section{Introduction}

In late 1941, the German occupiers of Paris arrested and released four members of the \textit{Acad\'emie des Sciences}, the most politically prominent of whom was the mathematician Emile Borel (1871--1956).%
\footnote{The others were the physicist Aim\'e Cotton (1869--1951), the physiologist Louis Lapicque (1866--1952), and the mineralogist Charles Mauguin (1878--1958).}  
Borel's wife, recounting the events in her 1968 autobiography \cite{marbo:1968}, admitted that she had no clue about the thinking behind the Germans' actions.  In this paper, we provide some information on this point, primarily from documents the Germans left behind when they fled Paris.  These documents were explored by a number of historians in the 1990s, especially Thalmann \cite{thalmann:1991}, Michels \cite{michels:1993}, and Burrin \cite{burrin:1995}.  They were thoroughly cataloged in 2002 \cite{beaujouan/etal:2002}.  But they have not been previously examined with an eye to Borel's case.%
\footnote{In general, there has been little research on how the Germans interacted with French mathematicians during the Occupation.  Audin's work on Jacques Feldbau \cite{audin:2008} is one the few recent contributions on this topic.}

The German documents we examine are in the AJ/40 series in the French national archives in Paris.%
\footnote{At CARAN, \textit{Centre d'accueil et de recherche des Archives nationales}.  We draw mainly on material in Boxes 558, 563, 566 and 567.}
This material was originally in the archives of the \textit{Milit\"ar\-befehls\-haber in Frankreich} (Military Command in France), known by the acro\-nym MBF.  We also draw on French documents from various sources, including the national archives, the archives of the \textit{Acad\'emie des Sciences}, and the archives of the French police.  There may be further relevant information in these archives and in the German archive at Freiburg \cite{martens:2002}.

Pierre Laborie has written eloquently about how difficult it has been for the French to find the distance to write dispassionate history of their interaction with the German occupiers (p.~182 of \cite{laborie:2000}):
\begin{quotation}
  The judgements passed on the collective attitudes and behaviors of the period between 1940 and 1944 are characteristic of this mixture of respectable intentions, fearfulness, and anxiety over all that is at stake in the realm of memory.  The extraordinary variety of personal experiences passed on by friends and family, as well as the topic's sensitivity and its popularity---everyone has an opinion on the matter---limit the dispassionate perspective of historians and their efforts to explain what happened.  When these perspectives stray too far from what is touchily guarded as `memorially correct' to a particular group or community, they are poorly received, and sometimes even suspected of insidiously striving to justify the unjustifiable.  The troubling question of behaviors during the Occupation is a recurring central theme in a debate that has been more about pronouncing judgement than about dealing with the issues and understanding their complexity.  Such questions are deeply relevant to our times because of their moral dimension, and yet too often they are reduced to the level of excessive generalizations, simplistic alternatives, or even summary judgments of the `all guilty, all collaborationist' variety.
\end{quotation}
In the case of science, this analysis often applies even to those who are not French, for the emotional ties between scientists and their historians easily cross temporal and geographic boundaries.  

The persecution of the Jews was incontestably the greatest shame of Vichy France.  As a condition for keeping their own jobs, almost all the decision makers in France, French and German, helped implement or at least accepted the dismissal of Jews from most employment, facilitating the murders that followed.  But we aim for the dispassion that Laborie urges, refraining from premature efforts to pass judgement, in the hope of thereby obtaining some degree of access to the ambiguous context of everyday life in this troubling period, without pretending to complete the picture.

The paper is organized as follows. To set the stage, we review Borel's life and career prior to the second world war (Section~\ref{sec:borel}) and the complexity and evolution of the German occupation of Paris up to the beginning of 1942 (Section~\ref{sec:complexity}).  We aim, while remaining brief, to provide enough information to make the story comprehensible to readers who know nothing about Borel or about the Occupation.

The German deliberations that led up to the arrest of Borel and his colleagues in October 1941 began a year earlier, when Paul Langevin was arrested.  So we look at what the archives tell us about relevant aspects of Langevin's arrest and subsequent developments during 1940--1941 (Section~\ref{sec:academy}) before turning to what they tell us about the events of October and November 1941 (Section~\ref{sec:arrest}).  

To conclude, we recount Borel's effort after the war to revisit the 1942 decision by the \textit{Acad\'emie des Sciences} not to select him as permanent secretary because his arrest had shown him to be unacceptable to the Germans (Section~\ref{sec:revisit}), and we discuss three other French mathematicians, Emile Picard, Albert Ch\^atelet, and Ludovic Zoretti, who played roles on the periphery of our story and whose own stories illustrate the complexities of this period (Section~\ref{sec:PicCha}).

As supplementary material, we provide transcriptions and translations of some of the most informative documents and texts on which we have drawn, some in German (Appendix~\ref{sec:german}), and some in French (Appendix~\ref{sec:french}).

\section{About Emile Borel}\label{sec:borel}

Emile Borel was born in a middle-class Protestant family in Saint-Affrique, in Aveyron in the center of southwest France.  He kept close ties with Saint-Affrique throughout his life.  After brilliant secondary studies, he went to Paris to prepare for the competitions leading to the \emph{grandes \'ecoles}, the schools where the French scientific and administrative elites are trained.  There he studied under the famous teacher Boleslas Niewenglowski along with the son of the mathematician Gaston Darboux, and he later recounted that it was at Darboux's home that he discovered his passion for scientific and especially mathematical research.  The \textit{Ecole Normale Sup\'erieure} was the place to pursue this passion. 

Borel immediately specialized in mathematics at the \textit{Ecole Normale}, beginning fundamental studies on divergent series, for which he introduced different modes of summability.  This soon led him to fundamental work on the measure of sets, which cleared the way for Lebesgue to construct his integral and revolutionize analysis \cite{hawkins:1975}.  Measure theory also led Borel, starting in 1905, to focus on probability theory.  He was the leading light in renewing mathematical probability at the beginning of 20th century, opening the way to the axiomatic formalization based on measure theory propounded by Kolmogorov in his {\it Grundbegriffe der Warscheinlichkeitsrechnung} \cite{kolmogorov:1933,shafer/vovk:2006}. 

Borel saw the mathematician as a citizen, and he put this conception into practice with works of popularization and philosophy \cite{bru/bru/chung:1999}.  From early in his career, he engaged in an active social and public life, especially through circles connected to the family of his talented wife Marguerite and his father-in-law, the mathematician Paul Appell.  Marguerite wrote fiction under the pen name Camille Marbo (for MARguerite BOrel).  In 1913, she won the Femina prize for her novel {\it La statue voil\'ee}.  In 1905, Borel and Marbo founded a monthly journal, the {\it Revue du Mois}, which for 10 years was a leading general intellectual outlet for the moderate French left.  Borel was active, along with Paul Langevin, in the \textit{Ligue des Droits de l'Homme} and its fight on behalf of Dreyfus.%
\footnote{See \cite{guieu:1998}; \cite{laberenne:1964}, pp.~253--254.}
In 1911, he and Marbo harbored Marie Curie in their home when Curie, then a widow, was under attack from the right-wing press over her affair with Langevin.%
\footnote{See \cite{guiraldenq:1999}, pp.~74--79.}

During World War I, Borel was a leader in putting the French technical and scientific elites at the service of the military.  In 1915, at the age of 44, he volunteered for the army himself in order to test acoustical devices for locating guns on the battlefield.  The same year, the mathematician Paul Painlev\'e, minister of public education, asked Borel to head a new office devoted to assessing and implementing inventions that could be used in the war.  For more information on his role at this time, see \cite{mazliak/tazzioli:2009}. 

Having been close to centers of power during World War I, and having been personally devastated by the war's slaughter of young graduates of the \textit{Ecole Normale}, where he had been deputy director, Borel turned increasingly to politics after the war.   Determined to work for greater social justice and more understanding between nations, he became prominent in the radical-socialist party, a very moderately leftist party that attracted many scientists and other scholars, so much so that its role in the French governments between the wars led some to call France {\it the republic of professors}.  In 1924, Borel was elected mayor of Saint-Affrique and member of parliament from Aveyron.  When Painlev\'e became the new Prime Minister, he named Borel minister of the navy, a position he held only for a few months.

Having been elected to the French \textit{Acad\'emie des Sciences} in 1921, Borel was also keen to use his political influence to help develop science and its applications.  He played a fundamental role in the creation of several major institutes of higher education, most importantly the Institut Henri Poincar\'e (IHP), inaugurated in Paris in 1928, which became the principal research center in France for mathematical physics and probability.  He saw to it that the IHP hosted the leading mathematicians and physicists of the 1930s, including Soviet scientists for whom travel was difficult, and German refugees fleeing the Nazis after 1933.

In January 1940, the University of Paris celebrated Borel's scientific jubilee, the fiftieth anniversary of his entrance to the \textit{Ecole Normale}.  All the great names of French mathematics and physics of the time were present, joined by foreign scientists who could come to Paris in spite of the war with Germany that had been declared in September 1939.  The local newspaper of Aveyron, {\it Journal de l'Aveyron}, celebrated him on the front pages of its 21 and 28 January 1940 issues. In June 1940, when the French army collapsed and the Germans occupied Paris and northern France, Borel was in Paris.  He and Marbo returned to Saint-Affrique, in the non-occupied zone, that summer, but they were back in Paris in the autumn of 1940, after he was dismissed from the mayoralty of Saint-Affrique by the new Vichy government.

\section{The German presence in Paris}\label{sec:complexity}

To a large extent, the Germans reproduced in Paris the complexity of the Nazi regime in Berlin, where various bureaucracies and militarized agencies competed for power without clear lines of authority among them.  This overview of the picture up to early 1942 draws on work by Burrin \cite{burrin:1995}, Frank \cite{frank:1993}, Nielen \cite{nielen:2002}, and Thalmann \cite{thalmann:1991}.

\paragraph{\textit{Milit\"ar\-befehls\-haber in Frankreich}.}

The MBF, headed from October 1940 to February 1942 by General Otto von St\"ulpnagel, was the most substantial German presence in Paris.  It was headquartered in the Hotel Majestic, near the \textit{Arc de Triomphe}.  It consisted of a security division and an administrative division, the \textit{Verwaltungstab}.  The \textit{Verwaltungstab}, consisting of 22,000 people in German military uniforms, was charged with overseeing the French governmental bureaucracy.  Its senior staff, numbering about 1,500, were mostly professionals detailed to Paris from various German government agencies, companies, and professional organizations.  It was divided into three large sections, a section responsible for coordination and personnel, an immense section responsible for economic matters, and a section responsible for other administrative matters, the \textit{Verwaltungsabteilung}.  

In 1941, the \textit{Verwaltungsabteilung} was headed by Werner Best (1903--1989), a prominent member of the SS.  It was organized into more than ten groups, charged with supervising domains of the French bureaucracy ranging from the police to the veterinary service.  We will be interested mainly in Group 4, responsible for schools and culture (\textit{Schule und Kultur}).  In the discussion of Borel's arrest, the group was represented by war administration adviser (\textit{Kriegsverwaltungsrat}) Dr.\ Dahnke.%
\footnote{Dahnke is identified only by his last name in the documents we have seen, and his first name is also absent from the secondary literature on the Occupation.  Almost certainly, however, he was Heinrich Dahnke, a Nazi bureaucrat who handled other aspects of international cultural affairs before and after his assignment in Paris \cite{lehmann/oexle:2004}.  According to Thalmann, the same Heinrich Dahnke served in the ministry of cultural affairs in Lower Saxony after the war (\cite{thalmann:1991}, p.~102).  There he assessed and supported claims for restitution by victims of the Nazis \cite{szabo:2000}.}

As we will see, Group 4 was not responsible for the arrest of Borel and his colleagues, but it may have played a role in their release, and its archives are our most substantial source of information.  When the Allies arrived in Paris, they found that most German archives had been burned or removed.  Group 4's archives were an exception,%
\footnote{See \cite{nielen:2002}, p.~47.}
and they are among the archives now preserved in the AJ/40 series at CARAN.  The archives of the units of the SS that investigated Borel and his colleagues were probably destroyed.  Some archives of the Abwehr (military intelligence), which actually arrested Borel and his colleagues, also survived the war and are now at Freiburg in Germany, but we have not been able to examine them.

The MBF's primary assignment was to put the French economy to work for the German war effort as effectively as possible, with a minimal expenditure of German manpower.  It achieved this through its control of the French governmental bureaucracy, by directing the allocation of raw materials and requiring French companies to fill orders for the German military.  The armistice agreement Marshall Philippe P\'etain signed on 22 June 1940 gave the Germans the right to require payment for the cost of the Occupation, and as they controlled the amount of this payment, they could pay French companies with money from the French treasury.

The armistice had authorized P\'etain to move his government to Paris, but the Germans never allowed this.  P\'etain remained isolated at Vichy, in the southern zone.  He was represented in Paris by the \textit{D\'el\'egu\'e g\'en\'eral du Gouvernement Fran\c{c}ais dans les Territoires Occup\'es} (DGTO), through which all French government communication with the Germans had to be directed.  During the period we are studying, this office was headed by Fernand de Brinon.  The southern zone was not occupied by the Germans during this period, and even after November 1942 when it was occupied, German permission was required for travel between the two zones.  The MBF required that P\'etain obtain its approval in advance for sensitive legislation and appointments, and it supervised the French bureaucracy, in Paris and in the prefectures, to make sure measures it cared about were implemented to its liking.

\paragraph{Otto Abetz, the German Ambassador.}

A second center of power was the German embassy on Rue de Lille on the left bank.  Otto Abetz held the rank of Ambassador and represented the German ministry of foreign affairs.  In theory, P\'etain's government was still at war with Germany, and the two countries did not have diplomatic relations.  But in practice Abetz was responsible for German relations with Vichy, overseeing P\'etain's supposed authority to legislate for France and to appoint the ministers in charge of the bureaucracy in Paris. 

Abetz was charged by Hitler to manage the politics of France.  He received his instructions directly from Hitler at the beginning of August 1940, when he was summoned to Hitler's summer home in Berghof.  Germany's immediate goals, Hitler explained to him, were to keep France weak and isolated from its neighbors.  It should have an authoritarian government, because this would help isolate it from England and the United States.  But there should be no real support for \textit{v\"olkisch} and nationalist forces.  Abetz should support both the left and the right in French politics, leaning at any time whichever way would maximize division.  The communists should not be wiped out, but they should not be allowed to become too strong, and in the immediate future the socialists should be supported as a counterweight to them.  Abetz shared these instructions with Werner Best, who passed them on to his lieutenants.%
\footnote{Thalmann \cite{thalmann:1991}, pp.\ 42--43, citing a German report in CARAN AJ/40/443:  \textit{Lagebericht} MBF III, December 1940--January 1941, 31 January 1941.}
Abetz and Best, both long-time Nazis, considered the military leadership of the MBF too conservative and too soft on the French.

Abetz was an accomplished student of French history and literature, and he saw the conflict between Germany and France in intellectual terms.  French culture should be purged of its degenerate elements, just as the Nazis had purged German culture, and the French needed to renounce their own claims to universalism and recognize the leadership role of German culture.  Abetz's ideas for improving France were never embraced by Hitler, who valued propaganda but was skeptical about changing the French and thought it better for Germany that France should continue to degenerate.

Abetz's organization in Paris had two main branches, a propaganda section and a cultural section.  The propaganda section organized and funded collaborationist groups, competing with and sometimes coming into conflict with Goebbels's \textit{Propaganda-Abteilung}.  It funded and largely controlled the French press and radio in the occupied zone.  The cultural section, which was involved in Borel's detention, was headed by Karl Epting, whom Abetz had also put in charge of the German Institute at the Sagan Hotel, also on the left bank.

Following Hitler's instructions, Abetz supported collaborationist groups in Paris with roots on the left as well as ones with roots on the right, and he used these groups to put pressure on P\'etain as the need arose.  He was especially supportive of Pierre Laval.  Laval had limited interest in P\'etain's national revolution but sought to convince the Germans that France had a future as a subservient junior partner.  He made unilateral concessions to the Germans while serving as P\'etain's minister of state during the first six months of the Vichy regime, and P\'etain dismissed him in December 1940, replacing him with Admiral Fran\c{c}ois Darlan.  When Borel was under arrest, in October and November 1941, Laval was in Paris under Abetz's protection.  When Laval returned to Vichy as P\'etain's prime minister in April 1942, he further aligned the regime with Germany.

\paragraph{The Nazi secret police.}

At the outset of the Occupation, the German military hoped that the SS, headed by Himmler and Heydrich in Berlin, would not be present in France, and so they included a small secret police unit within the MBF itself, in Section Ic.  But Himmler and Heydrich soon found ways to convince Hitler that their agents were needed in France.

The SS (\textit{Schutzstaffel}) was originally an organ of the Nazi party, not of the German state, and individuals could be deployed anywhere in the state bureaucracy while remaining SS officers; Werner Best is an example.  But by the time of the war, the German state police was fully integrated into the SS.  It was a vast organization.  Its intelligence service (\textit{Sicherdienst}, or SD) and secret police (\textit{Geheimes Staatspolizei}, or Gestapo) were just two of a myriad of units.

In October 1941, the SS presence in Paris was still small, and its agents were supposedly subordinate to the MBF, but it was increasingly coming into conflict with von St\"ulpnagel.%
\footnote{See \cite{burrin:1995}, pp.~96--97.}  
On the night of 2 to 3 October, a right-wing French group led by Eug\`ene Deloncle used explosives supplied by the SS to blow up seven Paris synagogues, wounding two German soldiers in the process.  This enraged von St\"ulpnagel, who was also concerned that the French population was being alienated by the increasing harshness of the Occupation, especially the execution of large numbers of hostages in retaliation for attacks against the occupiers.

Von St\"ulpnagel's reluctance to alienate the French more than necessary did not carry the day with Hitler.  In February 1942, von St\"ulpnagel resigned as military commander in Paris and was replaced by his cousin, Carl Heinrich von St\"ulpnagel.%
\footnote{Carl Heinrich von St\"ulpnagel remained in this position until 20 July 1944, when he participated in the unsuccessful plot against Hitler.  He was convicted of treason and hung on 30 August 1944 in Berlin.  Otto von St\"ulpnagel was put on trial by the allies after the war and committed suicide in 1946.} 
On 9 March 1942, Hitler decreed the appointment of Carl-Albrecht Oberg as top SS commander (\textit{H\"oherer SS- und Polizeif\"uhrer}) in Paris.  According to Thalmann,%
\footnote{See \cite{thalmann:2000}, p.\ 611.}
Oberg's appointment marked a turning point towards repression for both the German forces and the Vichy government.  It was followed by Laval's return to power, the appointment of Ren\'e Bousquet as general secretary of the French police, and intense negotiations between Oberg and Bousquet concerning what the French police would have to do in order to avoid being brought under direct German command.  The \textit{Verwaltungsabteilung} was shrunk, leaving the MBF even more focused on the economy.  Werner Best was sent to oversee the government of occupied Denmark.

As we will see, Borel and his colleagues were arrested neither by the MBF nor by the SS but by the \textit{Abwehr}, the branch of the German military responsible for espionage and counterespionage.  Its office in Paris, the \textit{Abwehrleitstelle} (Central office of military intelligence), was located at the Hotel Lutetia on the left bank.  In February 1944 Hitler abolished the \textit{Abwehr} and turned its functions over to the SS, but in 1941 it was independent of the SS, and its office in Paris reported to the military command in Berlin, not to the MBF.

\section{Events of 1940--1941}\label{sec:academy}

Borel's arrest in October 1941 should be seen in the context of German efforts during the previous year to rid the French professoriat of elements they considered undesirable, and also in a larger context that includes the struggle between Abetz and P\'etain and tensions within the German forces.

\paragraph{Purifying the professoriat.}

Beginning in the summer of 1940, P\'etain launch\-ed his National Revolution by promulgating a series of laws that facilitated the removal of undesirable individuals from education and other branches of the civil service.  The first, promulgated on 17 July 1940, authorized the removal of civil servants without cause.  Others forbade government employment of Freemasons and individuals of foreign origin and restricted government employment of women.  The most notorious, promulgated on 3 October 1940, forbade government employment of Jews.  Later legislation further discriminated against Jews and Freemasons.  

The French laws were more extreme in some respects than the Germans might have demanded, but the Germans used them to make the French bureaucracy take moral and practical responsibility for purges they intended to impose in any case.  The law excluding Jews allowed for exceptions, and the French thought this was important.  Some hoped exceptions would minimize the law's effects, or at least permit the separation of the wheat from the tares, as Darquier de Pellepoix, Laval's zealous {\it Commissaire aux Questions Juives}, put it.  Others were outraged that exceptions were possible.  But the Germans saw to it that were no exceptions.  According to Singer \cite{singer:1992}, the Vichy government authorized exceptions for fourteen leading Jewish professors at the University of Paris, but the Germans vetoed all fourteen.  

As the academic year 1940--1941 drew near, the Germans were obtaining lists of professors from the French bureaucracy and gathering intelligence to decide which ones should be eliminated.  Their deliberations were a three-way affair, involving Group 4, the embassy, and the SS.  The first discussions seem to have centered on the \emph{Institut libre des sciences politiques}, a private school, and the \emph{\textit{Coll\`ege de France}}, a prestigious government institution (\cite{thalmann:1991}, p.~103).

On 30 October 1940, before any recommendations had been made to von St\"ulpnagel, the secret police, at Epting's request, arrested Paul Langevin (1872--1946), a member of the \textit{Coll\`ege de France}.  Langevin was well known both for his left-wing political activity and for his accomplishments in physics.  According to Group 4's report for 18--24 November,%
\footnote{Transcribed and translated in Appendix~\ref{subsec:18_24_nov1940}.}
Langevin was not accused of continuing his political activity under the Occupation, but Epting wanted to make an example of him, to intimidate those in the university who might want to resist the occupiers, and to encourage those interested in collaboration.

Langevin's arrest came as the German forces, coordinated by Group 4, were involved in two tests of will with the \textit{Coll\`ege de France}.%
\footnote{See the archives from the French ministry of education in CARAN F/17/13385.}  
On 16 September, the MBF had informed the school that results from research in nuclear physics being conducted with the school's cyclotron by Fr\'ed\'eric Joliot-Curie together with German researchers would be exclusively for German eyes and could not be communicated to French authorities.  On 21 October, the MBF had informed the school that Langevin's presence on their faculty, along with that of two of his left-wing colleagues, Ernest Tonnelat and Henri Wallon, was incompatible with German interests and prestige.

For all we know, Werner Best and Group 4 at the MBF were complicit in Langevin's arrest, but von St\"ulpnagel was not pleased that such a step had been taken without his knowledge.  According to Group 4's report for 18--24 November, von St\"upnagel had been consulted neither before nor after the arrest.  Best was obliged to call together Epting, SS-commando Biederbick, and representatives of other secret police units to make it clear that no arrests of political significance were to be undertaken without consulting the MBF.  He also made it clear that political activity before the war was not grounds for measures against French scientists; von St\"ulpnagel wanted to be very cautious about interfering in personnel matters in French higher education.  It was quite another matter if Abetz, through diplomacy with Vichy, could get the French government to impose measures on which the embassy and the MBF were agreed.

On the morning of 11 November, the anniversary of the armistice of 1918, students mounted an impressive anti-German demonstration on the \textit{Champs Elys\'ees}.  This demonstration, often considered a turning point in French opinion, was violently repressed.  The Germans arrested 150 demonstrators, mostly lyc\'ee students.  Gustave Roussy, the rector of the Academy of Paris,%
\footnote{The Academy of Paris included the university faculties (science, letters, law, etc.) and the lyc\'ees in Paris; its rector was responsible for all faculty appointments in these institutions.} 
was dismissed, and the university was closed.  Roussy was replaced by J\'er\^ome Carcopino, a distinguished classicist and personal acquaintance of P\'etain's, who had been director of the \textit{Ecole Normale} since August.

Group 4 was instructed to continue their investigation, with Epting and Biederbick, of the politics of French academics, in case the reopening of Paris higher education could be made conditional on personnel changes.  The outcome of this investigation was explained by Dahnke in a Group 4 report dated 13 December 1940.%
\footnote{Transcribed and translated in Appendix~\ref{subsec:13_dec1940}.}
Even though Abetz's hand was strengthened by a new decree from Hitler, on 20 November 1940, giving the embassy exclusive authority over political matters and telling the MBF to attend to its military duties,%
\footnote{Group 4's report for the week of 25 November--1 December 1940, dated 2 December, CARAN AJ/40/563.}
Abetz and Epting did not manage to mount a purge of non-Jewish professors.  As Dahnke explains, Epting prepared a list of academics to be excluded, and then Biederbick and Group 4 pared it down by omitting some of the non-Jews.   The result, published in 1991 by Thalmann,%
\footnote{\cite{thalmann:1991}, pp.~354--361.} 
includes (with some repetitions) names of 109 individuals identified as Jews and 16 identified as non-Jews.  But on 7 December, the embassy contacted Group 4 orally to withdraw the suggestion that the French be required to exclude the 16 non-Jews before the university was opened, because not enough consideration had been given to excluding right-wing as well as left-wing enemies of Germany.  So on 8 December it was decided to demand only the expulsion of Jewish faculty members.  Harald Turner, representing the MBF, informed the education ministry%
\footnote{On 13 December, the date of Dahnke's report, Georges Ripert was replaced as minister of education by Jacques Chevalier.}
and Carcopino, the Paris  rector, that the MBF expected the French laws excluding the Jews to be applied without exception%
\footnote{``\dots der Miltaerbefehlshaber die restlose Durchfuehrung des franzoesischen Judengesetzes an den franzoesischen Hochschulen und ausserdem die Entfernung der juedischen Lehrkraefte in Bereich des Enseignement libre erwartet.''  The word \textit{restlos} is usually translated as ``complete'', but its literal meaning is ``without anything remaining''.}
to higher education and that Jewish teachers also be excluded from private schools.  Consideration of measures against academics because of their political past, especially those on the right, would have to await further study and an appropriate occasion, Dahnke concluded.

Although Borel had not yet retired from the university when these lists were prepared, his name does not appear on any of them.  Epting's list included the name of Fr\'ed\'eric Joliot-Curie, but the list approved by Group 4 did not.  The lists contain many errors, suggesting that the Germans and their informers were not very knowledgeable about the French academic world.  The mathematicians Maurice Fr\'echet and Georges Darmois for instance are identified as Jews.  A French Germanist active with Epting's institute is listed as anti-German.  In some cases, there are remarks about the political activity of which the individual is accused.  Some are labelled, probably erroneously, as Freemasons.  Others are labelled {\it Kolonialaufruf} (colonial call); perhaps this indicated that they supported de Gaulle's call to continue the fight against Germany from the French colonies.  

As for Langevin, the education ministry dismissed him from the \textit{Coll\`ege de France} on 19 November, under the authority of the law of 17 July, and the Germans released him from prison and put him under surveillance at Troyes, 180 kilometers from Paris, on 7 December.  Tonnelat and Wallon were eventually allowed to resume teaching at the \textit{Coll\`ege de France}.  On 18 January 1942, the MBF wrote to the French government to insist that it had exclusive rights under the Hague convention to the research in Joliot's laboratory because it was a military installation.  On request, however, it would inform the chief of state (P\'etain) about work currently in progress.%
\footnote{Much has been written about Joliot-Curie, who managed to participate in the communist Resistance while conducting his research in collaboration with the Germans.  Metzler \cite{metzler:2000} discusses the Occupation's impact on Joliot-Curie's scientific standing and gives additional references.}

Also relevant to our story is the launching of the clandestine newspaper \textit{l'Universit\'e libre} by three young communists, Jacques Decour, George Politzer, and Langevin's son-in-law Jacques Solomon.%
\footnote{See Favre \cite{favre:2002}, pp.~175--178; Racine \cite{racine:1987}; Raphael \cite{raphael:2000}, p.~722.}
The first four-page issue, in late November 1940, protested Langevin's arrest and dismissal, denounced the inadequacy of the university's protests against the arrest, and reported nearly unanimous reprobation for the government's antisemitic measures.  Joliot, it reported, had suspended his collaboration with the Germans in protest of Langevin's arrest, and the Faculty of Sciences, with only the physicist Eug\`ene Darmois dissenting, had voted that all its members should continue to teach regardless of their ``race''.%
\footnote{The quotation marks are in the original.  This issue other issues of \textit{l'Universit\'e libre} have been posted by the \textit{Conservatoire des m\'emoires \'etudiantes} at www.cme-u.fr.}
The three founders of \textit{l'Universit\'e libre} were arrested by the French police along with others in late February and early March 1941, turned over to the Germans, and executed in late May 1941, but \textit{l'Universit\`e libre} continued to appear irregularly until the Liberation, publishing altogether 104 issues.

\paragraph{Accusations against academicians.}

In the spring of 1941, the SS continued to keep an eye on the French professoriat, extending their attention to the academies in the \textit{Institut de France}.%
\footnote{The \textit{Institut de France} consists of five academies.  The most prestigious is the \textit{Acad\'emie fran\c{c}aise}.  Another is the \textit{Acad\'emie des Sciences}.}

On 15 April 1941, the SS addressed a memorandum to the MBF%
\footnote{Transcribed and translated in Appendix~\ref{subsec:SS_15_april_1941}.}
making two unrelated accusations concerning these academies:  
\begin{enumerate}
\item
A council of the \textit{Acad\'emie fran\c{c}aise} had held a vote on collaboration, in which the anti-German forces had a narrow majority.
\item
Two members of the \textit{Acad\'emie des Sciences}, Aim\'e Cotton and Charles Mauguin, had been involved in distributing {\it l'Universit\'e libre}.  
\end{enumerate}
Cotton and Mauguin were arrested along with Borel the following October.  

According to documents preserved by the MBF, the SS continued its investigation.  In a brief note for the week of 1--13 May 1941,%
\footnote{Transcribed and translated in Appendix~\ref{subsec:SS_13_may_1941}.}
the SD reported that the \textit{Acad\'emie des Sciences} did not have any distinct political orientation as a whole.  In an undated memorandum stamped with the date 24 May 1941,%
\footnote{Transcribed and translated in Appendix~\ref{subsec:SS_24_may_1941}.}
SS Major Biederbick reported that \textit{l'Universit\'e libre} was being distributed in the Latin Quarter by professors who were Freemasons, including Cotton, Mauguin, and Joliot.

\paragraph{Emile Borel in 1940--1941.}

Borel turned 70 on 7 January 1941, and so 1940--1941 was his last year at the University of Paris.  He officially retired at the end of the academic year, 30 September 1941.%
\footnote{In theory, French professors must retire by their 70th birthday, but in practice they are always allowed to work until the end of the academic year.  The decree extending Borel's service to 30 September 1941 was signed by the minister of education on 17 May 1941.}  In retirement, he remained active in the \textit{Acad\'emie des Sciences}, and he remained Borel, to whom others in the French higher education would turn for advice on appointments and prizes in mathematics.%
\footnote{The letters in Appendix~\ref{subsec:lacroix} give some glimpses of his activity in 1941--1944.}

Although we know little about Borel's activities at the university during 1940--1941, his last year there, we catch a glimpse of him in a story told by Carcopino in his 1953 memoirs.  Carcopino's account, if it is to be credited, is an interesting example of the solidarity that existed, at least at that date, within the French establishment.  In Borel's eyes, it seems, Vichy's ministry of education was not the enemy.  The story takes place in a monthly meeting of the university's council in which Carcopino, as rector and president of the council, sought to obtain the council's support for his opposition to a change in governance decreed by Jacques Chevalier, then Minister of Education, on 20 January 1941.  The rector served as president of the council, but its vice presidency was rotated among the deans of the five faculties.  Chevalier had proposed a change:  a permanent vice president would be appointed by the minister of education.  Carcopino saw this as a dangerous encroachment on the independence of the university, and all his colleagues on the council agreed, except Borel.  According to Carcopino, Borel thought that the appointment of a vice president by the government would have more advantages than disadvantages, because it would give the university administration greater unity and continuity.%
\footnote{``M.\ Emile Borel opina que la nomination, par le pouvoir central, d'un vice-pr\'esident du Conseil de l'Universit\'e qui serait permanent, pr\'esentait moins d'inconv\'enients que d'avantages et munirait l'administration universitaire de l'unit\'e et de la continuit\'e qui lui manquaient dans le syst\`eme actuel.  Ce petit discours de l'ancien ministre de la Marine d\'emontrait uniquement, par un exemple personnel, que chez un savant, longtemps m\^el\'e \`a la politique, l'homme du gouvernement primait l'universitaire.  Son intervention m'avait surpris et son argumentation choqu\'e.  Mais dans les dispositions o\`u je me trouvais, il me d\'eplaisait d'ouvrir une discussion avec le grand math\'ematicien qu'est M.\ Emile Borel\dots'' \cite{carcopino:1953}, pp.~273--274.}

It is plausible that Borel would have seen advantages in having a permanent vice president, but it is also understandable that Carcopino should have remembered their exchange as he did, for it underlines his own insistence on the independence of the university, which he needed to document in every possible way to refute accusations that he had conceded too much to the Germans.  And information in Borel's personnel file%
\footnote{In CARAN F/17/24854.} 
gives us some reason to doubt that he was present at the Sorbonne for the meeting.  Carcopino placed the meeting on the morning of Monday, 24 February 1941, more than a month after Chevalier's decree and on the very day Vichy announced Carcopino's appointment to succeed Chevalier as minister of education.  On 22 January 1941,
Borel had requested a leave of six weeks, from 20 February until the Easter vacation, saying that he hoped to obtain by 20 February a pass to cross the demarcation line so he could rest at Saint-Affrique.  An attached note from his physician states that he has a bad cold, complicated by acute pulmonary congestion, and recommends convalescence in the south of France for three months.  He was granted the leave he requested on 5 February 1941.  There is no indication in his personnel file that he returned to the Sorbonne to teach after this official leave expired.

On the other hand, the Germans did not hand out the needed passes readily, and we have no other evidence that Borel was in Saint-Affrique in 1941.  In the little she says about 1940--1941 in her autobiography,%
\footnote{The relevant passages are in the excerpt translated into English in Appendix~\ref{subsec:marbo}.}
Borel's wife does not mention any trip to Saint-Affrique.  She does report that the couple moved within Paris, from their elegant apartment on Boulevard Hausmann, where her mother and husband were shivering for lack of fuel, to a small apartment in Montparnasse, where Borel could be close to an intellectual milieu.  She recalls that after a brief vacation south of the Loire at the beginning of the summer, they had returned to Paris, where Borel had made contacts [\textit{o\`u il a pris des contacts}].  The wording suggests that the Germans were right to suspect that this 70-year-old man was plotting resistance.

\paragraph{J\'er\^ome Carcopino (1881--1970).}

As P\'etain repeatedly shuffled his cabinet during the first year of the Occupation, four ministers of education%
\footnote{Even the title and scope of the ministry varied.  Emile Mireaux was Ministre de l'Instruction publique et des Beaux-Arts, while Carcopino was Secr\'etaire d'Etat \`a l'Instruction Publique et \`a la Jeunesse.}
came and went:  
\begin{enumerate}
\item
Albert Rivaud (1876--1956), who left the post in July 1940, when Laval became Vice President of the Council,
\item
Emile Mireaux (1885--1969), who left the post in September 1940, 
\item
Georges Ripert (1880--1958), who left the post in December 1940, when P\'etain dismissed Laval, and 
\item
Jacques Chevalier (1882-1962), who left the post in February 1941, when Darlan became P\'etain's vice president.
\end{enumerate}
Chevalier, P\'etain's godson and a clerically-minded philosopher, had run into resistance with his plans for religious instruction in the primary schools, and Darlan replaced him with J\'er\^ome Carcopino, who had demonstrated his ability to negotiation with the Germans and keep the peace as director of the \textit{Ecole Normale} from August to November 1940 and then as rector after the 11 November 1940 demonstrations.  Carcopino remained minister of education until Laval returned to power in April 1942 and replaced him with the well known author and journalist Abel Bonnard (1883--1968), an avowed collaborator but less effective administrator, who remained minister until the Liberation.%
\footnote{Bonnard had been elected to the \textit{Acad\'emie fran\c{c}aise} in 1932.  Abetz had pushed for Bonnard's appointment as minister of education in February 1941, when Carcopino received the appointment.  Bonnard lived in exile in Spain after the war.}

Stubbornly loyal to P\'etain, Carcopino implemented Vichy's measures against Jewish teachers and students and appeased the Germans as he found necessary.  In 2008, after reviewing Carcopino's communications with the Germans in the MBF archives, the historian Alan Mitchell opined that collaborationism had no more fervent advocate \cite{mitchell:2008}.  According to \textit{l'Universit\'e libre}, Carcopino had been a \textit{Gauleiter} as rector; as minister he tried to be more flexible, more hypocritical, and more demagogic than his Vichy predecessors.%
\footnote{See Issue 14, 1 April 1941.}
The collaborationist press in Paris praised Carcopino at the outset of his tenure but soon found plenty to criticize.  When he softened Chevalier's proposals to strengthen religious instruction, the anti-clerical papers that had decried the measures fell silent, while the pro-clerical papers derided him.

In his memoirs, Carcopino complained bitterly about \textit{l'Universit\'e libre}'s attacks, but he acknowledged that the Vichy government in 1941 was a military dictatorship, with each minister exercising immense power within his own domain so long as he retained P\'etain's and Darlan's confidence.%
\footnote{See \cite{carcopino:1953}, pp.~298--299.}  
One power Carcopino exercised was to review dismissals made under the 17 July 1940 law.  He used the review to purge teachers considered incompetent or corrupt.  He also restored some individuals to their positions; even \textit{l'Universit\'e libre} conceded that he rehabilitated a few Freemasons.%
\footnote{See Issue 18, 15 May 1941.}
According to his own account after the war, he recognized merit regardless of political opinion, and he found reasons to dismiss individuals whose collaboration with the Germans threatened the university's independence.  This is documented not only in his memoirs but also in the criticism of the Paris press.  In the archives of the MBF, for example, we find a clipping of an article from the 28 August 1941 issue of \textit{Le cri du peuple}%
\footnote{Transcribed and translated in Appendix~\ref{subsec:godefroi}.}
that deplores Carcopino's appointment or retention of a whole list of individuals:  Charmoillaux, Maurain, Piobetta, Luc, Masbou, Chattelun, Santelli, and Hatinguais, and his dismissal of another, Jeanneret.  In its weekly reports, Group 4 took note of such articles matter-of-factly, as if they expressed public opinion.

After being replaced by Bonnard as minister, Carcopino returned to the \textit{Ecole Normale}, where he was director until the Liberation.  After the war, he was tried for collaboration and cleared on the grounds that he had made up for what he had done as a minister by his assistance to the Resistance afterwards.  He was elected to the \textit{Acad\'emie fran\c{c}aise} in 1955.  The ambiguities of his action continue to fascinate French historians \cite{corcy:2001}.  A recent assessment of his directorship at the \textit{Ecole Normale} credits him for finding ways to allow the Jewish students already there to complete their studies but faults him for barring the admission of additional Jews \cite{israel:1994}.

\paragraph{Raymond Voize and Albert Peyron.}

Of particular interest to our story is Carcopino's dismissal of Raymond Voize, a 51-year-old professor of German at the lyc\'ee Louis le Grand.  Voize's personnel file%
\footnote{In CARAN AJ/16/6176.}
indicates that he was an individual of ability and industry.  Born in 1889, he worked in commerce until he was 16, but then he managed to study at the lyc\'ee Voltaire, the Sorbonne, and the University of Halle.  He passed the \textit{agr\'egation} in German on his first attempt in 1913.  He also studied law and was interested in political science.  He was seriously wounded during the war, obtaining a 75\% disability pension.  His \textit{agr\'egation} entitled him to a position teaching in a lyc\'ee, but he refused assignments outside of Paris, teaching in a private school and working in a center for ``social and political documentation'' until finally being assigned to a Paris lyc\'ee in 1936.  In 1939, he ran into trouble with his superiors for using the name of his lyc\'ee to promote courses he was teaching during the holidays.  He remained in his position at that time in spite of exchanging very nasty letters with the Paris rector, Gustave Roussy, but in late 1940 his superiors became nervous about his relations with the German authorities.  One note in his file indicates that he had unloaded such a mass of denunciations on the Germans that they had not hidden their disgust from the French authorities.

According to a vicious article in the Paris newspaper \emph{l'Appel} on 7 August 1941,%
\footnote{Transcribed and translated in Appendix~\ref{subsec:voize1}.} 
Carcopino forced Voize to retire in July 1941 because of his connections with the Germans.  Carcopino confirmed this in his memoirs.%
\footnote{See \cite{carcopino:1953}, p.~351.}
Also on 7 August 1941, Voize published a long article in the weekly \emph{La Gerbe}%
\footnote{Transcribed and translated in Appendix~\ref{subsec:voize2}.}
proposing a high commission for French-German intellectual relations, which would award books on Germany to prize-winning French students.  This article began with a diatribe against the prospect that Carcopino might bring Roussy back to the rectorship in Paris.  When Roussy had been ousted from the rectorship in the aftermath of the 11 November 1940 demonstrations, he no longer held the position in the Faculty of Medicine he had held before being rector.  Carcopino, by his own account,%
\footnote{See \cite{carcopino:1953}, pp.~351--352.} 
had no intention of trying to reverse Roussy's removal from the rectorship, which had been signed by P\'etain himself, and he did not have the means to return Roussy to the Faculty of Medicine.  But he had wanted to appoint him director of the Institut Pasteur, and he had dropped the idea in April 1941 only because Roussy had reacted angrily at being offered so little.  In Carcopino's view,%
\footnote{See \cite{carcopino:1953}, pp.~550--551.}
Voize's article was the crudest and most unfair attack the Paris press ever made on him; he saw Abetz's hand behind it.
  
In his attacks on Carcopino and Roussy, Voize had an ally in Albert Peyron (1884--1947), a distinguished scientist at the Institut Pasteur who also fell victim to Carcopino's willingness to dismiss those who threatened the university's independence by their dealings with the Germans.  A dossier on Roussy in the German archives contains both letters of accusation and letters defending him from ``Peyron and Voize's calumnies''. 

Voize's article was noticed by the MBF's press service, which considered it outstanding.  It was also praised by Edmond Pistre-Carguel, the new Aryan Commissioner-Administrator for the publishing house Fernand Nathan under the German ordinance of 18 October 1940.%
\footnote{Pistre published under the name Caraguel.  His anti-English book {\it Angleterre contre la paix} was republished with additional chapters in 1940.  He died in 1942, and his book was banned in 1944.}
On 11 August 1941, Pistre-Carguel wrote to the Propaganda Abteilung Frankreich to ask them to grant Voize an opportunity to present his idea for a commission.  Such a commission, Pistre-Carguel argued, would give his own mission better support than he was getting from the education ministry.%
\footnote{See Appendix~\ref{subsec:pistre}.}

Voize was soon in contact with Group 4.  On 6 September 1941, he wrote to Dahnke at the MBF%
\footnote{The letter, transcribed and translated in Appendix~\ref{subsec:voizedahnke}, is addressed only to ``Monsieur'', but Dahnke annotated it and responded to it.}
that he wanted to go into private work instead of heading his proposed High Commission and wanted to see Dahnke again about creating an ``Institute for Languages and Culture''.  Dahnke responded On 20 September 1941%
\footnote{In a letter transcribed and translated in Appendix~\ref{subsec:dahnke_voize}.}
that he had not yet found the means for creating the institute and that he would be out of the office for five weeks, but that Voize should contact him again at the end of October.

\section{Arrest and release of the four academicians}\label{sec:arrest}

The four members of the \textit{Acad\'emie des Sciences} were arrested by the \textit{Abwehr} on 10 or 11 October 1941, imprisoned at Fresnes, and released on 13 November.  No explanation was ever given by the German or French authorities for their arrest or their release.

We will recount first how these events looked to the French at the time and then what we know from the German archives.

\subsection{From the French viewpoint}

In an autobiography published a year before her death in 1969 \cite{marbo:1968}, Borel's wife Camille Marbo gave a five-page account of her husband's arrest.%
\footnote{Transcribed and translated in Appendix~\ref{subsec:marbo}.}
A German officer, accompanied by four soldiers and a sergeant, came to their apartment at 2:00 in the afternoon, searched it, and then took Borel away at 5:00 with no explanation.  Marbo's brother managed to learn that Borel was at Fresnes only by taking a package there and getting it accepted.  Borel was never allowed visitors, but Marbo's packages, including clean clothes, were sometimes accepted.  

Marbo does not tell us the day in October 1941 when Borel was arrested, but it was almost certainly either 10 or 11 October.  According to a document in the archives of the French education ministry%
\footnote{In CARAN F/17/13385.}
it was Saturday, October 11.  
Mauguin was arrested on 10 October in a university building.  The next day the ministry contacted DGTO, asking them to find out from the Germans where Mauguin was being held and whether the arrest had any bearing on the university as a whole.  Two days later, they contacted the DGTO again to add that Borel and Cotton had been arrested on 11 October.  

News of the arrests spread quickly.  In the 12 October 1941 entry in his diary, Jean Gu\'ehenno writes that Borel has been arrested and that Langevin has been arrested anew.  The new Paris police commissioner, Gu\'ehenno says, is boasting of having arrested 1100 communists and anglophiles, and the Gestapo has declared the whole university suspect.%
\footnote{``Le nouveau pr\'efet de police, un amiral, bien entendu, se vante d'avoir d\`es maintenant fait arr\^eter onze cents communistes ou anglophiles.  Langevin, qui \'etait en r\'esidence surveill\'ee, est de nouveau emprisonn\'e.  Borel (soixante-seize ans) [sic] est aussi arr\^et\'e.  La Gestapo d\'eclare toute l'universit\'e suspecte.''  \cite{guehenno:1947}, p.~154.}
Langevin's biographers confirm that he was arrested at second time while at Troyes, interrogated, and then released after a few days.  Biquard places the arrest on a Wednesday at the end of September, and notes that the local German forces who arrested and interrogated him were not aware of his scientific stature.%
\footnote{See \cite{biquard:1969}, p.~95.  Lab\'erenne gives January 1942 as the date of the second arrest (\cite{laberenne:1964}, p.~302), but this is surely an error.}

As Marbo's account indicates, the French police were not involved in Borel's arrest.  All they could do was check after the fact on whether it had happened.  A report in archives of the Paris police,%
\footnote{Transcribed and translated in Appendix~\ref{subsec:frenchpolice}.}
dated 16 October 1941, states that they investigated the reported arrests of Langevin, Lapicque, Mauguin, Borel, and Cotton.  They confirmed the arrests of Lapicque, Mauguin, and Borel, giving 10 October as the date in each case.  They visited the address in Paris where Langevin had lived, but learned only that he was retired and now lived in Troyes.  A subsequent report in the same police archives, dated 7 November 1941, inveighs against communist militants whom it accuses of using the arrests to stir up fear and anti-German feeling to the detriment of France.

The arrests soon came to international attention.  On 18 October 1941, an article on occupied France in the London Times concluded with the comment, ``No reason has so far been given for the arrest in Paris by the German authorities of the well-known mathematician Emile Borel, a former Minister of the Navy.''  Two days later, on 20 October 1920, the Times devoted an entire article to the arrests.  As this article tells us, the reasons for the arrests remained a matter of speculation.  

\begin{quote}
From Our Special Correspondent, French frontier,\footnote{References to the Haute Savoie suggest that the correspondent was stationed in Switzerland.} Oct.\ 19

\smallskip
The Vichy Government today confirmed the arrest in Paris by the German authorities of five prominent professors of the University -- namely, MM.\ Borel, Langevin, Lapicque, Mauguin and Cotton. 

According to some sources they are charged with spreading de Gaullist propaganda, according to others with pro-British sentiments, while some newspapers lay emphasis on the fact that the political activity of MM.\ Langevin and Borel has been well known since the time of the Front Populaire. The brother of the former Prefect of the Seine Department, M.\ Villey, has also been arrested, together with his son and daughter, on a charge of alleged de Gaullist activity. 

Judging by opinion in Haute Savoie, these arrests are causing bewilderment, as even the former political opponents of these scientists cannot believe that they have been arrested on account of their personal views. Some light may be thrown on the affair by a recent article published by Laval in his newspaper the {\it Moniteur du Puy-de-D\^ome}. In this he says that now that Germany has conquered her enemies, who are those of France, the latter must conquer her disorder and errors and hold out her hand to Germany -- the Queen of Europe. Laval then declares that all French persons who are still imbued with anti-German prejudice should be at once dismissed from public offices. He adds that this prejudice now exists mainly among the intellectuals, where it may be regarded as a remnant of anti-Fascism. 

In Haute Savoie the view is expressed that the above `ultimatum' by Laval inspired the Vichy Government to act accordingly, as the French authorities certainly lent a hand in the arrest of the Paris professors.
\end{quote}
Although the \textit{Times} correspondent was mistaken to believe that the French government had a hand in the arrests, the collaborationist press was loudly supporting Laval's demand.  On 23 October 1941, for example, in a violent article opposing Roussy's return the university, {\it l'Appel} asked how Roussy could be allowed to hold his position when Professors Lapicque, Cotton, Mauguin, Borel, Villey, and Saintelag\"ue had been locked up for Bolshevist Gaullism.  

Marbo hoped that the \textit{Acad\'emie des Sciences} would petition for their members' release, and she marshaled support from three members, Maurice de Broglie, Elie Cartan, and Paul Montel.  But the leadership of the Academy feared that speaking out would risk the Academy's abolition.  Marbo also went to talk with Carcopino; he told her his hands were tied.  Carcopino confirms this is his 1953 memoir, where he writes of feeling sad and helpless when she told him about Borel's not getting the blankets she had brought to the prison for him.%
\footnote{See \cite{carcopino:1953}, p.~472.}
He also recalls that in November 1940, when he was rector and Georges Ripert was minister of education, Ripert had asked him to approach the Germans about Langevin.  Neither the \textit{Coll\`ege de France} nor the other institutions where Langevin worked fell under the rector's jurisdiction, but Ripert appealed to Carcopino for help because of Carcopino's reputation for dealing with the Germans.  Dissuaded from approaching the MBF directly by officials at the DGTO, Carcopino asked the scientist and industrialist Georges Claude, the most outspoken proponent of Collaboration in the \textit{Acad\'emie des Sciences}, to talk with the Germans about Langevin.  Claude told Carcopino he would but changed his mind after finding that the permanent secretaries of the \textit{Acad\'emie des Sciences} would not lend their names to the effort.%
\footnote{For Carcopino's account, see \cite{carcopino:1953}, pp.~347--349.  According to Corcy-Debray (\cite{corcy:2001}, p.~231), Claude's letter to Carcopino explaining the refusal of the permanent secretaries to be involved is in the French National Archives, AN 3W122.}

During the period when Borel and his colleagues were imprisoned, the representative of the ministry of education at the DGTO was Maurice Roy.  Born in 1888, Roy had been a professor of German at the Lyc\'ee Saint-Louis until being promoted to the rank of inspector in 1940.  He was delegated to the DGTO in March 1941.  According to Carcopino's biographer St\'ephanie Corcy-Debray%
\footnote{See \cite{corcy:2001}, p.~76.}
Roy intervened with the Germans almost systematically on behalf of students and professors who had been arrested, at the same time as he worked on many other conflicts involving education and youth movements.  It was his job, for example, to submit legislation and regulations to the German censors before they were published by the French government.  Roy was doubtlessly in contact with the Germans about Borel and his colleagues as soon as their arrests were known to the education ministry, and, as we shall see, he was eventually involved in their release.

In her biography of her husband Emile Cotton, Eug\'enie Cotton states that the four prisoners were released on 13 November,%
\footnote{See Appendix~\ref{subsec:cotton}.}
and this date is confirmed by a handwritten note on an MBF document prepared in Feburary 1942.%
\footnote{Transcribed and translated in Appendix~\ref{subsec:mauguin_11_2_1942}.} 
Other information is consistent with this date.  Marbo tells us that Borel fell ill with double pneumonia the day after he was released, and by 19 November he was writing a note to his colleague Albert Lacroix.%
\footnote{See the letter in Appendix~\ref{subsec:lacroix}.}

We have not found any account of Borel's recollections about his interrogation.  The biologist Maurice Caullery, a fellow member of the \textit{Acad\'emie des Sciences} and friend of Borel, Cotton, Mauguin, and Lapicque, stated in his own memoirs that the interrogation of the four was a sham.%
\footnote{When he himself was imprisoned at Fresnes on 2 April 1941, Caullery recalled that ``... mes amis Borel, Cotton, Mauguin et Lapicque, arr\^et\'es ainsi en octobre 1941, \'etaient rest\'es \`a Fresnes plus d'un mois, avant d'\^etre rel\^ach\'es apr\`es un simulacre d'interrogatoire.''  About a dozen academicians, including Caullery, were held for two days, 2 April to 4 April.  See \cite{telkes:1993}, pp.~239--243.}

\subsection{From German documents}

The German documents we have examined do not tell us definitively why Borel and his colleagues were arrested or why they were released, but they provide many hints.

\paragraph{The accusations by Voize and Peyron.}

We have already seen something of the activities of Raymond Voize and Albert Peyron from August to October 1941.  In his 6 October note to Voize, Dahnke had told Voize that he would be gone for five weeks and asked him to contact him again at the end of October.  On 23 October, \textit{l'Appel} published its denunciation of the possibility of Roussy's return as rector.

We may surmise that Voize and Peyron were both in contact with Dahnke at the end of October, because the MBF archives include a memorandum,%
\footnote{Transcribed and translated in Appendix~\ref{subsec:dahnke_1_nov_1941}.}
signed by Dahnke and dated 1 November 1941, that lists 15 individuals they had accused of being ``representatives of the Freemason and Bolshevist view in higher education and the administration of public education.''  Gustave Roussy heads the list, and Dahnke reports that Peyron claims told Carcopino him personally that he wanted to make Roussy rector again.  Next are eight individuals that Carcopino had appointed or retained or granted a pension even though their collaborationist bona fides were suspect:  Luc, Guyot, Hatinguais, Zoretti, Masbou, Santelli, and Chattelun.  These are people Voize and the Paris press had been denouncing since August.  Then Fr\'ed\'eric and Ir\`ene Joliot-Curie.  Then our four academicians:  Mauguin, Cotton, Lapicque, and Borel.  And then Gustave Monod, who had been in the education ministry before the war, had been demoted to teaching in a lyc\'ee because he refused to enforce the expulsion of Jews, and had then retired.

Voize and Peyron surely knew in late October that Borel, Cotton, Lapicque, and Mauguin had already been arrested.  Was Dahnke, catching up on his work after returning to Paris, compiling information that Voize and Peyron had provided to him or others earlier?  Or did Voize and Peyron add names of individuals they knew had already arrested in order to make their accusations against Carcopino more persuasive?

It is noteworthy that Voize and Peyron are singling people out only for their opinions, not for actions against the occupiers.  For most of those named, ``Bolshevist'' was a wild exaggeration, but they were on left, and the Joliot-Curies really were communists.  

Were any of those named Freemasons?  From the beginning, the Vichy regime took measures against Freemasonry as well as against the Jews.  A law of 13 August 1940 prescribed that Freemason lodges be closed, their properties impounded and sold.  Civil servants and public officials were ordered to break any links with the dissolved lodges and not to affiliate anew if they were reconstituted.  This repression increased dramatically in the summer of 1941.  On 11 August 1941, a new law was promulgated stipulating the publication of lists of Freemasons in the {\it Journal Officiel} and directing that civil servants who had been Freemason dignitaries would automatically lose their jobs.   Huge lists of names were published in the subsequent weeks, and many did lose their jobs.  Borel was not on any of the lists, but Voize himself, who was a Freemason before returning to Catholicism, was.

As Jean Gu\'ehenno noted in his diary in 1941, the lists of Freemasons refuted the myth of the Freemasons' power, for hardly anyone prominent in the Third Republic was on the lists.%
\footnote{``Vichy, pour orienter la haine des Fran\c{c}ais, a fait publier les noms des franc-ma\c{c}ons.  Mais la publication n'a pas eu l'effet esp\'er\'e.  On ne pouvait mieux faire pour d\'etruire la l\'egende de la puissance de la franc-ma\c{c}onnerie.  Cette liste montre avec \'evidence qu'\^etre franc-ma\c{c}on pouvait assez bien conduire \'a \^etre instituteur, voire percepteur, mais presque aucun des grands noms de la troisi\`eme R\'epublique ne s'y retrouve...''  \cite{guehenno:1947}, p.~155, entry for 11 October 1941.}  
Carcopino noted that there was only one Freemason among the 80 members of the faculty of letters at the University of Paris.%
\footnote{``La Facult\'e des Lettres de l'Universit\'e de Paris, \`a laquelle j'ai eu l'honneur d'appartenir, n'aurait compt\'e, en 1941, sur plus de 80 professeurs, ma\^itres de conf\'erences, et charg\'es de cours, qu'un seul franc-ma\c{c}on.''  \cite{carcopino:1953}, p.~378.}

\paragraph{The selection of hostages.}

Attacks on German soldiers by young French communists began in France after Hitler broke with Stalin and attacked the Soviet Union in June 1941.  A German officer was killed in Paris in late August, and in early September the Germans killed more than twenty people in retaliation.  Hitler considered this response too moderate; he demanded that at least a hundred French be killed for every German.  Another German officer was killed in Nantes on 20 October, and the Germans killed nearly a hundred hostages in retaliation on 22 and 24 October.  Revulsion over these murders was important in shifting French publich opinion towards de Gaulle and the Resistance.

Von St\"ulpnagel, keenly aware of negative impact the killing of hostages would have on French public opinion, wanted the French public to perceive it as being directed towards communists and Jews, not towards the French population in general.  So in September, the hunt was on for potential victims who could be identified as communists.  This may have played a role in the arrest of Borel and his colleagues in October.

When Dahnke returned to his office in Group 4 of the Verwaltungsabteilung in late October, after an absence of five weeks, he found a circular issued on 18 September by his colleagues in Group 8, responsible for German oversight of the French courts.  The circular asked for names to be added to the hostage lists, and it specifically asked for names of students and professors.  In his report for the week of 20--26 October,%
\footnote{Transcribed and translated in Appendix~\ref{subsec:hostage_27_10_1941}.}
Dahnke conceded that many French natural scientists were radical leftists and could be held responsible for inspiring communist subversion.  But, he pointed out, if the Germans shot one of them as a hostage, the act would be held against them in the university community for a very long time.

\paragraph{The intervention of the German Navy.}

The earliest mention of Borel's name that we have found in the MBF archives comes in a letter dated 25 October 1941, from the office of the Commanding Admiral of the German Navy in France to the secret police section of the MBF.%
\footnote{See Appendix~\ref{subsec:navy} for a transcription of the German original.}  
Here is the body of the letter in English translation:
\begin{quote}
The research section of the naval weapons office at the Headquarters of the War Navy is currently working in Paris on important problems of nuclear physics together with the Parisian ``Curie'' Institute.  The German scientists depend on
the French scientists in this work.  To be named, among others, are the mathematician Prof.\ Borell [sic], the physicists Prof.\ Langevin and Cotton, the crystallographer Prof.\ Mauguin and the mineralogist Prof.\ La Picque [sic], the last two of the Sorbonne.

According to information from the 
representatives of OKM at the Commanding Admiral in France, the forenamed French scientists have been under arrest for some time.  Because the collaboration between the German and French scientists will be very difficult under these circumstances, and the continuation of the military scientific research may become impossible, we would like to be advised about whether the misdeeds committed by the arrested French scientists are so serious that their arrests must be upheld.
\end{quote}
It seems unlikely that Fr\'ed\'eric Joliot-Curie and his French assistants would have been getting help from their aging colleagues in Paris, let alone from Langevin in Troyes, or that Joliot-Curie could have convinced the German scientists in the laboratory that this was the case.  But it is plausible that the German scientists would have solicited this letter as a gesture to Joliot-Curie.  Joliot-Curie's role in obtaining Langevin's transfer from prison to Troyes in December 1940 is documented by Burrin (\cite{burrin:1995}, pp.~315--322).

\paragraph{Dahnke's reports on the arrest and release.}

By 10 November, at least, Dahnke knew that the academicians had been arrested.  On that day he writes:%
\footnote{We transcribe the German original in Appendix~\ref{subsec:dahnke_10_nov_1941}.}
\begin{quote}
The Alst arrested professors Mauguin, Cotton, Lapicque and Borel (see the attached note of 1 November 41), and lyc\'ee teachers Aubert and Cazalas (see the attached memorandum). I have gotten in contact with the Alst (Major Dr.\ Reile), spoken with the expert in charge of the file, Captain Kr\"ull, and transmitted to him the note of 1 November 41 in order to bring connections of which he had been unaware to his attention.  He intends to extend his investigations to this circle.
\end{quote}
``Alst'' was the Nazis' abbreviation for the Paris \textit{Abwehrleitstelle}, the agency of the military responsible for counterespionage in France.  On 9 October 1941, the Alst had launched a wave of arrests (code-named \textit{Fall Porto}, or Operation Porto) against a loosely organized network of the Resistance, the Hector network, that was funded by Vichy's air force and led by Alfred Heurteaux.  Launched because German infiltrators reported that the Hector network was forming cells that might engage in assassinations of Germans as the communists had done, the operation also extended to other groups suspected of espionage against the Germans.  Many of those arrested were executed or deported; others were released.%
\footnote{See the memoirs of Oskar Reile \cite{reile:1962}.  There may be more information about the arrests of Borel and his colleagues in the archives of the Paris Alst in Freiburg (\cite{martens:2002}, pp.~243--248.  German dossiers on many of those arrested on 9 October 1941 are in CARAN AJ/40/1500--1521 (\cite{beaujouan/etal:2002}, pp.~510--512).}

In a handwritten note added on 18 November to the memorandum of 10 November, Dahnke reports that the Borel, Cotton, and Mauguin had been released in the meantime after a note from Roy -- presumably Maurice Roy, Carcopino's liaison with the Germans.

Dahnke seems to have been meeting with the Alst weekly.  In his next report, he writes:%
\footnote{The memorandum transcribed in Appendix~\ref{subsec:dahnke_25nov41}.}
\begin{quote}
Consultation with Captain Kr\"ull on 25 November 41.  He released Professors Mauguin, Cotton, Lapicque and Borel.  Their interrogation showed that all of them, especially Cotton, still candidly stand by the political ideas they advocated before and during the war.  The openly declared that they expect England and America's political system to rescue France.  But they emphatically denied that they had in any way acted on their views, especially with students.  The intelligence service is not in position to prove such activity, though our informers claim it has taken place.  In particular, it is impossible to arrange a confrontation with students from the circles these professors were supposed to have influenced, because the informers did not identify any such students by name.

During a consultation between this expert, Dr.\ Epting, and Dr.\ Biederbick, we considered putting the four professors under police surveillance outside Paris, as was done with Langevin in Troyes.  In view of such a political attitude on the part of the four professors, who are described by rightist circles as the center of extreme leftist and anti-German tendencies in the Sorbonne --- see the references in the newspapers --- it cannot be expected that they will refrain from expressing their opinions in French scientific circles.  This would justify a measure of this kind, even if it is impossible to produce witnesses to prove their political activity.  No active participation by Langevin was proven either.  But we must deliberate carefully before implementing such a measure; I have therefore asked Dr.\ Epting to personally discuss the matter at the Embassy as soon as possible. 
\end{quote}
We are left to conjecture that Abetz decided that it was not worthwhile to exile the four academicians from Paris.

We learn a little more from documents prepared by Group 4 in early 1942, when Charles Mauguin again came to the Germans' attention.  Mauguin was elected to the council of the university, and the clandestine press, noting that a collaborationist had also stood for the post, celebrated the election as a defeat for the Collaboration.  Group 4 duly looked into the matter, and in a report to the Verwaltungsabteilung's police section, dated 11 February 1942,%
\footnote{Transcribed and translated in Appendix~\ref{subsec:mauguin_11_2_1942}.  This report was not prepared by Dahnke; it is signed by ``P.''}
they noted that Mauguin had been accused by his right-wing opponents of being close to the communists but that no proof had yet been found; they would investigate further.  A handwritten note on this report notes that Mauguin had been already been arrested in October 1941 and released on 13 November 1941.

\subsection{Why did the Germans release the four academicians?}

The German documents answer this question reasonably clearly.  They tell us repeatedly, in a variety of contexts, that von St\"ulpnagel was not willing to punish French scientists who could not be proven to have taken action against the Occupation, because to do so would interfere with keeping order in the university.  

Perhaps the intervention of the German Navy had some influence, but on the whole, we must take at face value Dahnke's report that the Alst released the academicians for lack of evidence that they were distributing \textit{l'Universit\'e libre} or other pamphlets or otherwise actively resisting the Occupation.  On 15 December 1941, Gilbert Gidel, then the Paris rector, met with von St\"ulpnagel, and according to Gidel's account of the conversation,%
\footnote{Transcribed and translated in Appendix~\ref{subsec:gidel}.}
von St\"upnagel made explicit the bargain the Germans were offering university authorities:  the Germans would refrain from harsh measures so long as students and professors stuck to their academic business.

Dahnke's reports do not tell us whether his communication with the Alst, reported on November 10, played a role in the prisoners' release on November 13, but he ties it to a note from Maurice Roy, Carcopino's liaison with the MBF.  Roy could communicate with the Germans only through the MBF, and Dahnke or someone else at the MBF would have been responsible for conveying the note and its context to the Alst.

\subsection{Why did the Germans arrest the four academicians?}

The Alst was surely working with the SS, and the accusations the Alst made against the four academicians, concerning the distribution of pamphlets, were those the SS had been investigating, at least in the case of Cotton and Mauguin, since April 1941.  The Alst's failure to interrogate the academicians seriously suggests that they still did not have much evidence.  So the arrests must have been precipitated by other events.  

Perhaps the Alst swept the four academicians up as part of Operation Porto, without paying attention to their prominence, simply because they were on some long list of possible subversives.  But the Alst and the SS were surely aware of von St\"ulpnagel's policy of requiring clear evidence of specific acts for politically significant arrests.  So something more must have been in play.  Several hypotheses, not mutually exclusive, are suggested by what we have learned:
\begin{enumerate}
\item
  The SS may have been looking for ways to challenge von St\"ulpnagel.
\item
  The Alst and the SS may have wanted to add prominent scientists to the reserve of prisoners ready to be shot in retaliation for attacks by the communists.
\item
  Abetz and Epting may have pushed for the arrests as part of a campaign to destabilize Carcopino and Darlan.  \end{enumerate}

\subsection{The arrest of academicians in April 1942}

On 2 April 1942, the Germans arrested and imprisoned for two days a dozen members of the \textit{Institut de France}, including two members of the \textit{Acad\'emie des Sciences}, Aim\'ee Cotton, who had been arrested with Borel in October, and the biologist Maurice Caullery.  Carcopino's liaison at the \textit{Milit\"arbefehlshaber} was told by the Germans that the academics were arrested because they had been receiving the clandestine periodical \textit{La France Continue}.  Caullery, who compared notes with the others who were arrested, confirms that the Germans accused them all of receiving clandestine material.  Caullery told them that this was true; he had received the material in a sealed envelope, had not known who sent it, and had destroyed it.  Caullery conjectured that the Germans' only evidence against those arrested was the presence of their names on a distribution list for the clandestine material.  Aside from the questioning when they were arrested, the academicians were not interrogated further before being released on 4 April.%
\footnote{See \cite{telkes:1993}, pp.~239--243.}

According to Carcopino the arrests were discussed by P\'etain and his ministers.  P\'etain tried at first to minimize the significance of the arrests, observing that none of those arrested were in the \textit{Acad\'emie fran\c{c}aise}.  But Carcopino countered that some of them were P\'etain's own colleagues in \textit{l'Acad\'emie des Sciences morales}.  P\'etain was persuaded to act on behalf of the prisoners in a way he had not in previous cases:  Darlan was told to have de Brinon intervene with the Germans.  When the prisoners were released on April 4, Carcopino rejoiced in the effectiveness of P\'etain's action, but he later concluded that the Germans had released the prisoners because the arrests had achieved their goal of bringing P\'etain around to changing his government.  If the Germans' purpose had been to attack the Resistance, Carcopino reasoned, they would have arrested more members of the \textit{Acad\'emie des Sciences}, the citadel of the Resistance within the \textit{Institut de France}, rather than arresting eight or nine members of Carcopino's own academy, the \textit{Acad\'emie des Inscriptions}.  After the war he found support for this opinion in telegrams from the German Embassy in Paris to von Ribbentrop in Berlin, one on 21 March 1942 proposing that the crisis created by the fiasco of the Riom trial be used to oust Carcopino and other undesirable ministers, a second on 3 April 1942 indicating that P\'etain had now conceded to Laval the departure of Carcopino and the minister of agriculture.%
\footnote{See \cite{carcopino:1953}, pp.~560--565; \cite{corcy:2001}, p.~232.}
By the middle of April, Carcopino had left the government, returning to the directorship of the \textit{Ecole Normale} in Paris, and Laval was in charge of the Vichy government. 

Although Carcopino's analysis of the events of April 1942 may have exaggerated his own importance, we should not underestimate the importance of the arrests of French elites in the struggle between the Germans and Vichy.  Vichy was a government of technocrats more than any other in French history, and P\'etain expanded its ranks more than any other French leader.  He may have wanted the French population as a whole to support his National Revolution, but this revolution was primarily concerned with the elites.  They were the ones to be brought back from their decadence, purged of the influence of Jews, Freemasons, and communists, restored to bring France back to its glory.  What P\'etain surely sought to forestall above all was that Hitler should destroy the talented tenth in France as he had done in Poland.  By rushing to expel Jews and Freemasons from government service, P\'etain had established a corps of functionaries who could plausibly be protected from the Germans by their collaboration.  But the arrests of the academicians demonstrated that the Germans could still strike at France's best.  Emile Borel and his other aged colleagues in the \textit{Institut de France} were surely of little practical importance to Vichy or to the Germans.  But if the Germans could eliminate them, how much of the elite that mattered to P\'etain could he have saved?

\subsection{The sequel for Borel and Marbo}

According to Marbo's memoir, Borel was arrested again in 1942, this time by two French policemen, who returned him home after three hours, reporting that the German officer to whom they had taken him had rejected him as a prisoner because of his age.  Marbo does not give a date for this second arrest.  Was it on 2 April when the other academicians were arrested?  Perhaps, but Marbo's testimony that the French police executed the arrest counts against this, as the Germans seem to have carried out the 2 April arrests without French help.

The second arrest persuaded Borel and Marbo to return to Saint-Affrique.  They did so, after obtaining the necessary permission for crossing the line of demarcation, in October 1942.  According to Marbo, they worked with the Resistance in Saint-Affrique.  Borel's contribution was to allow the Resistance to use a forest he owned, but Marbo, who was thirteen years younger and apparently in much better health at that point, did what she could to feed and otherwise help Jews and other fugitives.  By the spring of 1944, Borel and Marbo were back in Paris, where they stayed with Marbo's brother Pierre and Marbo helped him deliver messages for the Resistance.  On D-Day, Borel and Marbo were staying clandestinely in a clinic in the Passy district of Paris, where Borel underwent surgery.%
\footnote{See \cite{marbo:1968}, pp.~306--309.}

After the war, Borel resumed his activity in Paris and served again as mayor of Saint-Affrique,%
\footnote{According to Marbo, it was always she rather than Borel who actually did the work of mayor.} 
from 1945 to 1947.  He died in Paris in 1956.  Camille Marbo died in 1969.

\section{Who should have succeeded Picard?}\label{sec:revisit}

Emile Picard died on 11 December 1941 at the age of 86.  Borel, a relative by marriage,%
\footnote{Borel's father-in-law Paul Appell and Picard had both married nieces of the mathematician Joseph Bertrand.}
would normally have attended the funeral.  But he seems not to have done so, perhaps because he was still recovering from his ordeal or perhaps out of prudence.  The announcement of the funeral lists his wife as attending as part of the Appell family.

Picard had been one of two permanent secretaries for the \textit{Acad\'emie des Sciences}, the other being the mineralogist Fran\c{c}ois Lacroix (1863--1948).  As permanent secretary, Picard had represented the mathematical sciences, which had subsections for geometry, mechanics, astronomy, geography, and physics.  On 2 February 1942, the academy elected the physicist and Nobel laureate Louis de Broglie to succeed him.  

Had he not been arrested, Borel would have been a natural successor to Picard, as the post had been held previously by the academy's most senior pure mathematician.  Jacques Hadamard was more senior than Borel, but he had fled to the United States because he was Jewish.  After the war, Borel contended that he would have been elected but for his arrest and sought to have the election reversed.

\paragraph{De Broglie's election}

The file at the academy's archives for its meeting of 2 February 1942 confirms that the election of de Broglie was a delicate matter.  A letter from the ministry, dated 30 December 1941, reveals that in response to the ministry's official condolences to Picard's widow, Lacroix had mentioned the vacancy of the post of permanent secretary and the academy's intention of proceeding to an election to fill it, and asking for an appointment with Carcopino.  He met with Carcopino on the morning of 7 January 1942.  

The next day, Lacroix wrote to the members of the academy laying out how the new permanent secretary would be elected, with reference to documents going back to 1803, when the post had been introduced.  There would be nothing out of the ordinary about such a notice, but the circumstances seem to have required supplementary documentation.  Perhaps someone had suggested leaving the post vacant, because Lacroix cites an article of the regulations that had been amended in 1816 to state that the academy would \emph{not deliberate about whether or not to elect someone to the post} (Lacroix's emphasis) but would elect a commission of six members from the section (mathematical sciences in this case), which would produce a list of candidates in consultation with the academy's president.  Lacroix then proposed a calendar for the process:  the commission would be named on 19 January, it would name the candidates on 26 January, and the election would take place on 2 February if a quorum of 40 could be assembled; failing the quorum the election would take place on 9 February by simple majority of those present.  The file also contains a tally showing that there had not been a quorum at any of the meetings for December 1941.

Another document, in Lacroix's hand, indicates the results of a secret committee meeting on 12 January that chose the commission members to be elected by the assembly the following week.  The commission consisted, naturally, of the most senior members of the five subsections:  Borel for geometry, Villat for mechanics, Deslandres for astronomy, Bourgeois for geography, and Cotton for physics, along with Maurain as the most senior of the other members of the section.  At the bottom of the document is a discreet acknowledgement of Hadamard's existence:  \emph{It is a matter of the most senior members present in Paris.}

The two candidates proposed by the commission on 19 January were Louis de Broglie and Elie Cartan.  In his January 7 notice, Lacroix had prescribed that the commission would not make a report, on the surprising grounds that its work involved competition among colleagues.  

On 2 February, 41 members of the academy were present, and 39 voted:  22 for de Broglie, 15 for Cartan, and two with blank ballots.  The same day, Lacroix sent Carcopino an excerpt of the minutes proclaiming de Broglie's election and asking him to confirm the choice.  The published minutes of the following week's meeting%
\footnote{See CRAS 214, 16 February 1942, p.\ 294} 
records the confirmation and reproduces de Broglie's very proper acceptance speech.  He merely affirmed that the \textit{Acad\'emie des Sciences} had to play its role in the difficult times being endured for the sake of the country's recovery.  Such rhetoric would have been completely satisfactory from the viewpoint of the Vichy government, which always insisted that it was working for the restoration of the country and for its future triumphs.

\paragraph{Borel's remonstrance}

The Liberation of Paris disturbed the activities of the \textit{Acad\'emie des Sciences} very little.  The assembly did not meet during the week when the Liberation took place, but they did meet the following Monday, August 28, starting late because of the difficulties experienced by the trains, and unanimously voting to join the other academies of the \textit{Institut} in congratulating the provisional government and thanking the allied troops and French forces%
\footnote{See CRAS, 219, p.\ 225.}
The following week, on Monday, 4 September, they went into secret session and decided to expel Georges Claude, the scientist and industrialist who had been so vocal a proponent of collaboration.%
\footnote{See CRAS, 219, p.\ 264.} 

Soon thereafter, on 23 September 1944, Borel wrote to Lacroix proposing a reconsideration of the election of de Broglie as permanent secretary.  In this letter,%
\footnote{Transcribed and translated in Appendix~\ref{subsec:lacroix}.}
Borel recounts what Jean Vincent, president of the academy, has told him about the Academy's actions in 1941.  Vincent had opposed appealing to the Germans for their release, because de Brinon had persuaded him that such an appeal would be dangerous both for the prisoners and the Academy itself.  For the same reason, Vincent had opposed Borel's selection as permanent secretary, and he now agreed that Borel would have been elected had he not been arrested.  At the beginning of 1942, Borel had thought he had enough votes to be elected in spite of Vincent's opposition, but as Lacroix very well knew, Carcopino had forced him to withdraw his candidacy.  In light of this history, Borel felt that he was owed ``reparation''; Louis de Broglie should resign, and Borel should replace him.  Borel had already talked to friends of Louis and his brother Maurice, and he thought they could be persuaded that this was appropriate. 

Borel's plan seems to have gotten off the ground.  He wrote to Lacroix again on 6 October 1944, saying that he had a good conversation with Maurice and was convinced that Louis would indeed resign.  But Louis did not resign, and Borel never became permanent secretary.  We have no further evidence concerning the attitudes of the de Broglies and other members of the Academy, but it is easy to imagine them hesitating to revisit decisions made during the Occupation.  Did reparation mean an admission of fault?  Who else might deserve reparation?  If Borel had a claim, did Hadamard also have one?

\section{Three French mathematicians}\label{sec:PicCha}

Emile Picard was a prominent mathematician of the generation that preceded Borel.  Albert Ch\^atelet and Ludovic Zoretti were younger than Borel.  The three appear in the German documents we have encountered in ways that illustrate the complexity of the interaction between the French and the occupiers.

\paragraph{Emile Picard (1856--1941).}

Picard was known for anti-German outbursts during World War I and for advocating ostracism of Germany afterwards \cite{lamande,mazliak/tazzioli:2009}.  He was also very well known for his right-wing views, and from Vichy's point of view, he was a natural ally.  In August 1940, the {\it Journal de l'Aveyron}, falling in line with Vichy's national revolution, wrote at length about the educational reforms needed to undo the damage from the radical and anti-clerical conceptions that had prevailed during the previous regime and led to the disaster.   Never mentioning the name of their native son Emile Borel, they quoted Picard, who had no connection with Aveyron, at length.  The {\it Journal} made no mention of Borel's arrest in October 1941 or his release in November 1941.

The Nazis, in contrast, had no use for Picard.  In a report on the prospects for collaboration with French mathematicians, dated 20 December 1940, the Nazi mathematician Harold Geppert, who edited the \textit{Zentralblatt f\"ur Mathematik und ihre Grenzgebiete} during the war, mentioned Picard in connection with his leadership of the International Mathematical Union, from whose quadrennial meetings the Germans had been excluded in the years following World War I, and which had not met since 1932.%
\footnote{A copy of this letter is in CARAN AJ/40/56-.  For more on Geppert's role, see Siegmund-Schultze \cite{siegmund-schultze:1993}.}
According to Geppert, the Germans had decided to create a new international organization for mathematics rather than reviving the Union.  
\begin{quote}
So it is useless to look for what remains of the earlier minutes of the Union, which are presumably in the hands of the permanent secretary of the \textit{Acad\'emie des Sciences}, Prof.\ Emile Picard, who is the intellectual leader of the Union and a thoroughly anti-German scientific polemicist.  So the idea of undertaking a search of P.'s house has been dropped.
\end{quote}

We have already mentioned another note in the German archives, dated 15 April 1941, that makes negative mention of Picard.%
\footnote{Appendix~\ref{subsec:SS_15_april_1941}.}
This note claims that a council of the \textit{Acad\'emie fran\c{c}aise} had decided against supporting collaboration, by a secret vote of 4 to 3.  Picard was listed at one of those voting against collaboration, along with Paul Val\'ery, Georges Duhamel, and Maurice de Broglie.  The note erroneously identifies Picard as a member of the French Popular Front of the 1930s.  Perhaps the informer or the German agent confused Emile Picard with Emile Borel, whose party had joined the Popular Front in the election of 1936, even though Borel himself had not stood for election that year.  Borel was never a member of the \textit{Acad\'emie fran\c{c}aise}, let alone its council.  Picard had been elected to the \textit{Acad\'emie fran\c{c}aise} in 1924 as successor to the physicist Charles de Freycinet in Chair 1, then traditionally held by scientists.

Picard's attitude towards collaboration was nuanced.  According to Audin \cite{audin:2009}, Picard argued, in a discussion in the \textit{Acad\'emie des Sciences} in November 1940, that efforts to distribute its \textit{Comptes rendus} should not involve any direct relationship with the Germans.  Yet in correspondence with Alfred Lacroix, his fellow permanent secretary of the \textit{Acad\'emie des Sciences}, he maintained that P\'etain had rightly agreed to a ``very general'' collaboration, necessary for an indefinite time in order to avoid France's being completely crushed.  Audin concludes, on the basis of the correspondence with Lacroix, that Picard was an antisemite of the ordinary French variety (as opposed to the Nazi variety) and that he was almost a collaborationist.  Because he died in December 1941, he did not see where P\'etain's policies led.  Members of his family took different paths.  One of his sons-in-law, Jean Villey, was caught by the French police on 13 October 1941 in the act of distributing Gaullist propaganda, delivered to the Germans, and condemned to two years in prison.%
\footnote{See CARAN AJ/16/7117.}
Another, Louis Dunoyer de S\'egonzac, who was close to the extreme-right \textit{Action Fran\c caise}. was faulted after the war for accepting appointment to Jean Perrin's position at the Sorbonne in 1941, after Perrin had fled to the United States.

\paragraph{Albert Ch\^atelet (1883--1960).}

In February 1941, when Carcopino had left the Paris rectorship to become minister, he postponed the problem of finding a new rector acceptable to the Germans by putting Charles Maurain in the position temporarily.  As dean of the Faculty of Sciences, Maurain had been taking his turn as vice president of the university council, and Carcopino reasoned that he was therefore in line to step in as president of the council, or rector.  But, Maurain was due to retire at the end of September.  Carcopino himself had been appointed to the Paris rectorship after the student demonstrations of 11 November 1940, and his success in calming the students while satisfying the Germans had led to his appointment as minister.  Who could follow this act?

According to a report in Group 4's archives,%
\footnote{Report for the weeks 25--31 August and 1--6 September, dated 8 September 1941, in CARAN AJ/40/563.}
Carcopino went personally to the German embassy and proposed three names:  Paul Hazard, Olivier Martin, and Pierre Renouvin.  Abetz rejected all three.  The comparative literature professor Hazard and the law professor Martin were unacceptable because they knew nothing of Germany.  Renouvin, a historian, was unacceptable because he had written about the causes of the first world war from the French point of view.  In Abetz's view, German cultural-political goals in France required a Paris rector with personal and professional connections with Germany.  So von St\"ulpnagel asked Abetz to make his own suggestions.  This seems not to have been so easy.  Abetz's connections with French journalists and literary and political figures, dating from the 1930s, were not matched by similar connections with academics, and it is not clear that anyone with the stature and skill to be rector satisfied his desiderata.  The only name Abetz had suggested by September 8 was that of Ch\^atelet.  Ch\^atelet had been rector at Lille from 1924 to 1937, and then worked in the Ministry of Education until 1940.  Abetz's argument for Ch\^atelet was that he had promoted French-German exchanges while rector in Lille.  

The appointment needed to be made before Maurain left the rectorship, so von St\"ulpnagel agreed to meet with Carcopino to settle the matter.  In the meeting, Carcopino explained that Ch\^atelet was on the extreme left and therefore unacceptable to the French government.  Carcopino had dismissed Ch\^atelet from the rectorship in Caen, and to reinstate him now would undermine his authority with the students by making it clear that he was acting on the order of the Germans.  Carcopino made a new suggestion:  Gilbert Gidel, who was highly regarded by Friedrich Grimm, a German professor who visited France regularly to lecture, assess the situation, and advise Abetz.  Von St\"ulpnagel agreed to Carcopino's appointing Gidel.  

Group 4's report for December 1941--January 1942 states that Gidel had been received by von St\"upnagel and pledged his loyalty.  Gidel's report on the meeting%
\footnote{Transcribed and translated in Appendix~\ref{subsec:gidel}.}
puts the matter differently:  Gidel was committed to keeping the university calm and creating an atmosphere of work.

Gidel remained rector until the Liberation, then returned to teaching law; he died in 1958.  Ch\^atelet finished his academic career as dean of the Faculty of Sciences at Paris and then went into politics; he was a candidate for president of the republic in 1958.

\paragraph{Ludovic Zoretti (1880--1948).}

Zoretti had brilliantly launched a mathematical career at the \textit{Ecole Normale}, coming to Borel's attention and contributing regularly to the \textit{Revue du Mois}.%
\footnote{For example, he translated into French Volterra's 1901 \textit{Prolusione} at the University of Rome.  The translation was the very first article in the \textit{Revue du Mois}.}
A specialist in the theory of functions in Borel's style, Zoretti was proposed for the Peccot lecture in 1908--1909 and became a professor at Caen, but his mathematical career was damaged by criticisms by Brouwer.  He joined the SFIO in 1914 and became a labor organizer, very active in the CGT.  In the 1930s, his militant pacifism led to his expulsion from the SFIO and suspension from teaching.  The Vichy regime dismissed him completely from teaching, replacing him at Caen with Robert Fortet.  

During the Munich crisis in September 1938, Zoretti created a stir by accusing Blum of risking a war that would destroy a civilization in order to make life easier for a hundred thousand Jews.  In December 1938, he went back to denouncing Nazi atrocities and supporting the international league against antisemitism.  But during the Occupation, he aligned himself with Marcel D\'eat, a collaborationist leader who had also participated in the pacifist movement.  In 1941, Zoretti published a nationalist and antisemite pamphlet \emph{France, forge ton destin}.  His late conversion was mocked by the collaborationist \textit{Je suis partout}.%
\footnote{See \cite{epstein:2008}, pp. 215--217, 264.}
In the spring of 1944 after D\'eat became Minister of Labor, Bonnard and D\'eat gave Zoretti the task of creating a workers' university.  He went into hiding when Paris was liberated.  Condemned for collaboration, then arrested in June 1946, tried again, and sentenced to eight years in prison, he died in the 1948 at Camp Carr\`ere in the Lot.

\bibliographystyle{plain}

\pagebreak

\appendix

\section{Some German documents}\label{sec:german}

These documents are from archives of the MBF (\textit{Milit\"ar\-befehlshabers in Frank\-reich}) now preserved in series AJ/40 at CARAN (\textit{Centre d'accueil et de recherche des Archives nationales}) in Paris.  In most cases, we provide both a transcription of the original German document and an English translation.  Except where otherwise noted, the documents were prepared with typewriters.  We have silently correced obvious typographical errors, but we have preserved misspellings of names.  We omit some headings and other material outside the main body of the documents.  For the most part, we have not tried to transcribe the handwritten notes on the documents, which we found difficult to read.  The documents are in chronological order, from November 1940 to April 1942.

The documents were apparently all in the files of Group 4 (Schule und Kultur) of the \textit{Verwaltungsabteil\-ung}.  Sometimes this group is referred to simply as ``V kult''.  Similarly, Group 2, responsible for the police, is referred to as ``V pol''.

We gratefully acknowledge Anke Piepenbrink's help with the German transcriptions and translations.

\subsection{Group 4 re Langevin's arrest, 18--24 November 1940}\label{subsec:18_24_nov1940}

\textit{Excerpt from the report by Group 4 for the week of 18--24 November 1940.  Otherwise undated, it is now in box AJ/40/563.}

\subsubsection*{German original}

Die Vorg\"ange, die zur Schliessung der franz\"osischen Hoch\-schulen in Paris gef\"uhrt haben, sind Herrn Abteilungsleiter bekannt.  Von weiteren Vorf\"allen gleicher Art in Paris ist bei Gruppe 4 bislang nichts bekannt geworden.  Nach einer Mitteilung des Vertreters des frans\"osischen Unterreichtsministeriums ist nunmehr auch die Universit\"at in Dijon geschlossen worden, weil sich dort anl\"assich der R\"uckkehr der in Dijon beheimateten Pariser Studenten gleichfalls Studentunruhen zugetragen haben.  N\"aheres ist bislang nicht bekannt.  Ein Bericht von Bezirkschef C ist angefordert.

An 30.10.1940 hat das SS-Sonderkommando auf Betreiben Dr.\ Eptings ohne vorherige oder nachherige Verst\"andigung des Milit\"ar\-befehlshabers in Frank\-reich in Zusammen\-wirken mit der Geheimen Feldpolizei den politisch extrem links gerichteten, wissenschaftlich bedeutenden Professor Langevin verhaftet.  Anlass zu dieser Verhaftung haben, wie Ermittlungen in der Berichtswoche gezeigt haben, nicht konkrete Vorw\"urfe gegen Langevin \"uber%
\footnote{Someone has crossed this word out and written in by hand another word, which we cannot decipher.} 
Fortsetzungen seiner politischen T\"atigkeit nach der deutschen Besetzung, sondern die Absicht gegeben, ein Exempel gegen die in der franz\"osischen Intelligens in Paris vorhandenen deutschfeindlichen Kr\"afte zu statuieren und denjenigen Bestrebungen in der franz\"osischen Hochschullehrerschaft Auftrieb zu geben, welche auf eine politische Neueorientierung Frankreichs gerichtet sind.

Herr Ministerialdirektor Dr.\ Best hat in einer Besprechung in der Berichts\-woche Dr.\ Epting als Vertreter der Deutschen Botschaft und Dr.\ Biederbick als Vertreter des SS-Sonderkomamndos nachdr\"ucklich darauf hingewiesen, dass der Oberbefehlshaber f\"ur die Zukunft in jedem einzelnen Fall vorherige Anfrage verlangen m\"usse, wenn eine Verhaftsmassnahme von politischer Bedeutung beabsichtigt sei.  Zu der Frage, inwieweit die Zukunft [sic] Massnahmen gegen durch eine deutschfeindliche Bet\"atigung in der Vergangenheit bekannte franz\"osischen Wissenschaftler angebracht seien, hat Herr Min.\ Direktor Dr.\ Best in dieser Sitzung in Beisein des Vertreters der Botschaft, des SS-Sonderkommandos, der Abwehr und des Ic festgestellt, das eine Einflussnahme deutscherseits auf die Personalverh\"altnisse an den franz\"osischen Hochschulen entsprechend der von Herrn Oberbefehlshaber erteilten Weisung nur mit gr\"osster Zur\"uckhaltung ge\-hand\-habt werden k\"onne.  Etwas anderes sei es, wenn Botschafter Abetz auf diplomatische Wege gegen\"uber der franz\"osischen Regierung in Vichy allgemein solche Personalmassnahmen durchsetzen k\"onne.  Nach Ansicht der Gruppe 4 muss erwartet werden, dass sich die Botschaft ihrerseits bei derartigen Massnahmen mit Mitlit\"arbefehlshaber in Einvernehmen h\"alt.

Die bei der Gruppe 4 laufenden internen Arbeiten zur Aufkl\"arung der politischen Haltung der franz\"osischen Hochschullehrer werden zusammen mit dem SS Sonderkommando und Dr.\ Epting fortgesetzt, um gegebenenfalls die Wieder\-er\"offnung der Pariser Hochschulen, an die einstweilen u.U.\ etwa zu Weihnachten gedacht ist, gegebenenfalls an personelle Bedingungen kn\"upfen zu k\"onnen.

Der neue Minister-Staatssekret\"ar f\"ur \"offentliche Unterrichtswesen, M.\ Ripert, hat seinen Dienstsitz nach Paris verlegt.  Es ist anzunehmen, dass er mit den Milit\"arsbefehlshaber in Frankreich pers\"onliche Verbindungen aufnehmen und m\"oglicherweise die vorstehend geschilderten Zussamenh\"ange zur Spreche bringen wird.  Der bisherige Vertreter des frans\"osicher Unterrichtsministeriums in Paris hat um Ueberlassung der Liste der am 11.11.\ Verhafteten gebeten.  Ueber diesen Wunsch ist bis jetzt nicht entschieden, da der Vorgang beim Ic und der Gruppe Polizei noch l\"auft.  Hierau [sic] ist erw\"ahnenswert, dass sich unter den \"uber hundert Verhafteten nur 19 Studenten befinden.

\subsubsection*{English translation}

The events that led to the closing of French higher education%
\footnote{In general, we have translated the German word \textit{Hochschulen} as ``higher education'' rather than ``universities'', because many of the institutions that were closed, such as the \textit{Coll\`ege de France} and the \textit{grandes \'ecoles}, are not usually called universities.  Academic laboratories and most of the lyc\'ees remained open.}
in Paris are known to the department head.%
\footnote{The weekly report is evidently addressed to the Werner Best, head of the \textit{Verwaltungsabsteilung}.}  
So far Group 4 is not aware of any further incidents of this kind.  According to a report from the representative of the French Education Ministry, the University in Dijon is now also closed, because the student unrest also occurred there when students living in Paris returned home there.  No details are yet known.  We have asked for a report from regional chief C.

On 30.10.1940, on Dr.\ Epting's initiative and without consulting the MBF before or after the fact, SS special units, working together with the secret military police, arrested Professor Langevin, who is oriented to the extreme left politically and is scientifically important.  Inquiries during the past week revealed that the reason for the arrest was not any concrete accusation against Langevin concerning continuation of his political activity after the German occupation.  Instead, the intention was to make an example against the anti-German forces that exist in the Paris intelligentsia and to support those in the French academy who are working towards a new political orientation for France.

Mr.\ Ministerial Director Dr.\ Best, in a meeting this week with Dr.\ Epting as representative of the German Embassy and Dr.\ Biederbick as representative of the SS special units, indicated emphatically that the military commander must be consulted in advance in the future in every single case where an arrest of political significance is being considered.  As for the question of the extent to which measures are to be taken in the future against well known French scientists because of past anti-German activity, Mr.\ Ministerial Director Dr.\ Best explained, in this same session, in the presence of representatives of the Embassy, the SS special units, the military security service, and the Ic,%
\footnote{The secret military police was located in section Ic of the MBF.}
that according to instructions from the military commander, German interference in personnel matters in French higher education should be taken only with the greatest reluctance.  It would be different if Ambassador Abetz could push through such general personnel measures through diplomacy with the French government in Vichy.  In Group 4's view, the Embassy must be expected come to agreement with the MBF on such measures.

Group 4's continuing internal work on clarifying the political mindset of the professoriat in French higher education will be carried out together with the SS special units and Dr.\ Epting, in case the reopening of French higher education, which is now contemplated for around Christmas if circumstances permit, can be tied to conditions on personnel.

The new minister for education, Mr.\ Ripert, has moved his office to Paris.  It is to be expected that he will establish personal relations with the MBF, and perhaps the issues just mentioned can be raised.  The current representative of the French education ministry in Paris has asked for a list of those arrested on 11.11.  We have not yet decided on the response to this request, because the operation by Ic and the police group is still in progress.  Here it is noteworthy that there were only 19 university students among the over one hundred arrested.

\subsection{Group 4 re proposed purge, 13 December 1940}\label{subsec:13_dec1940}


\textit{Group 4's report, signed by Dahnke, on excluding anti-German teachers from higher education in the occupied zone.  Now in CARAN AJ/\-40/\-567.}

\subsubsection*{German original}

\noindent
\underline{Betreff:}  Entfernung deutschfeindlicher Hochschullehrer von den franzoesischen Hochschulen im besetzten Gebiet.

\smallskip
\noindent
\underline{Sachbearbeiter:}  KVR Dr.\ Dahnke

\smallskip

Durch Schreiben vom 6.12.1940 an den Herrn Oberfehlshaber hat Herr Bot\-schaf\-ter Abetz angeregt, anlasslich der Wiedereroeffnung der Pariser Hochschu\-len zum 1.1.41 von der franzoesischen Regierung die Entfernung juedischer und deutsch\-feind\-licher Hochschullehrer von den Pariser Hochschulen zu fordern.  Der Sachbearbeiter fuer kulturpolitische Fragen bei der deutschen Botschaft, Dr.\ Epting, hat sich entsprechend einer Weisung des Botschafters Abetz mit dem zustaendigen Sachbearbeiter bei der Gruppe 4 Schule und Kultur wegen Durchfuehrung dieser Anregung in Verbindung gesetzt und die mit 1 bezeichnete einliegende Liste ueberreicht.  Bei der Bearbeitung dieser Liste, zu der auch das SS-Sonderkommando (Dr.\ Biederbick) herangezogen wurde (vergl.\ Anlage 2), ist die Entfernung der in der Anlage 1 blau angemerkten franzoesischen Hochschullehrer von den Parisier Hochschulen ins Auge gefasst worden.

Die deutsche Botshaft hat auf muendlichen Wege am 7.12.1940 ihre Anregung insoweit rueckgaengig gemacht, als sie sich auf franzoesische Hoch\-schul\-lehrer bezog, deren Entfernung wegen ihrer deutschfeindlichen Betaetigung gefor\-dert werden sollte.  Massgebend hierfuer war, dass die getroffene Auswahl den Bot\-schafter Abetz insofern nicht befriedigte, als die auf der ehemaligen franzoesischen Rechten stehenden deutschfeindlichen Kraefte bei dieser noch nicht genuegend beruecksichtight werden seien.

Die Botschaft hat ihre Anregung insoweit aufrechterhalten, als sie sich auf die Entfernung der juedischen Hochschullehrer bezog.  Nach einer Besprechung bei Herrn Miniserialdirektor Dr.\ Best am 8.12, bei der die Botschaft durch Dr.\ Epting vertreten war, hat Staatsrat Turner dem franzoesischen Unterrichterministerium und dem Rektor der Akademie und der Universitaet Paris demgemaess eroeffnet, dass der Miltaerbefehlshaber die restlose Durchfuehrung des franzoesischen Judengesetzes an den franzoesischen Hochschulen und ausserdem die Entfernung der juedischen Lehrkraefte in Bereich des Enseignement libre erwartet.  Die geplanten Massnahmen bezueglich der franzoesischen Hochschul\-lehrer mit einer ausgesprochenen deutschfeindlichen Vergangenheit soll durch weitere Ermittlungen, insbesondere bezueglich der franzoesischen politischen Rechten, vorbereitet werden und einer geeigneten Situation vorbehalten bleiben.

\subsubsection*{English translation}

\noindent
Re: Weeding out anti-German professors from French higher education in the occupied zone.

\smallskip
\noindent
Expert:  KVR Dr.\ Dahnke

\smallskip

In a letter of 6 December 1940 to the military commander, Ambassador Abetz proposed that the French government be required to dismiss Jewish and anti-German university professors on the occasion of the reopening of the Paris institutions of higher education on 1 January 1941.  On Ambassador Abetz's instructions, the cultural politics expert of the German Embassy, Dr.\ Epting, has contacted the expert responsible at group 4 ``Schools and Culture'' and provided the list designated as attachment 1.  On the basis of an examination of this list, which also involved the special command of the SS (Dr.\ Biederbick; see attachment 2), the dismissal of the French professors marked in blue in attachment 1 was being considered. 

The German Embassy, in an oral communication on 7 December 1940, has asked to withdraw its proposal with respect to the French professors whose dismissal was to be demanded because of anti-German activities.  The reason given was that Ambassador Abetz was not satisfied with the selection, because it had not taken sufficient account of the anti-German forces that had formerly stood on the French right.  

The Embassy continues, however, to support its proposal with respect to the Jewish university professors.  After a meeting in the office of section director Dr.\ Best on 8 December, where the Embassy was represented by Dr.\ Epting, Staatsrat Turner\footnote{The SS officer Harald Turner had earned the title ``Staatsrat'' in his previous career as a jurist.} informed the French Minister of Public Instruction and the rector of the Academy and University of Paris that the Military Command expects the complete application of the French laws on Jews in French institutions of higher education, and also the exclusion of Jewish teachers from non-public schools.  

The planned measures with respect to French university professors with definite past anti-German activity will be reconsidered at an appropriate moment after further investigation, especially with respect the French political right.

\subsection{SS re \textit{l'Institut de France}, 15 April 1941}\label{subsec:SS_15_april_1941}

\textit{This memorandum from the SS to the MBF is in CARAN AJ/\-40/\-566.}

\subsubsection*{German original}

\noindent
Paris, den 15 4 1941

\smallskip
\noindent
An den Milit\"arbefehlshaber in Frankreich Kommandostab

\smallskip
\noindent
Paris Hotel Majestic

\smallskip
\noindent
Betr: Institut de France, Paris und angeschlossene Akademien

\smallskip

Eine vor einiger Zeit im Vorstand der Acad\'emie Fran\c{c}aise stattgefundene geheime Abstimmung \"uber die Frage der Politik des Marschalls P\'etain und die Collaboration ergab, dass von den sieben anwesenden Vorstandsmitgliedern 3 f\"ur und 4 gegen Petain bezw.\ gegen die Collaboration gestimmt haben. Negativ haben sich ge\"aussert:
\begin{itemize}
\item
Paul Val\'ery, Schriftsteller
\item
Georges Duhamel, Schriftsteller
\item
Duc Maurice de Broglie, Historiker [sic]
\item
Emile Picard, Angeh\"origer der ehem. franz.\ Volksfront [sic]
\end{itemize}
Ihre Zustimmung gaben:  Andr\'e Bellessort, Schriftsteller und st\"andiger Sekre\-t\"ar der Acad\'emie Fran\c{c}aise
Kardinal Baudrillart, sowie der bekannte Publizist Abel Bonnard.

Im \"ubrigen sind sowohl in der Acad\'emie Fran\c{c}aise, wie auch in der Acad\'emie des Inscriptions et belles lettres, der acad\'emie des sciences, der acad\'emie des beaux arts und der acad\'emie des sciences morales et politiques Juden, Pro\-bolsche\-wisten, Freimaurer, aber auch der Action Fran\c{c}aise angeh\"orende und nahestende Kreise Mitglieder gewesen und sind es noch heute. In einer Reihe von F\"allen waren diese Wissenschaftler zudem auch Mitglieder der Union des Intellectuels Fran\c{c}ais pour la Justice, la Libert\'e et la Paix, die als diejenige Organisation angesprochen werden muss, die weite Kreise der franz\"osischen Intellektuellen auf die alleinigen Basis der absoluten Deutschfeindlichkeit vereinigte. Es ist erwiesen, dass Mitglieder der Acad\'emie des sciences : Aim\'e Cotton und Charles Mauguin, beide der offiziell nicht mehr t\"atigen Union des Intellectuels Fran\c{c}ais angeh\"orig, in illegalen Flugbl\"attern - z.B. ``l'Universit\'e libre'' - gegen die im besetzten Frankreich bestehenden Verh\"altnisse und gegen Deutschland gehetzt haben. Die Untersuchungen hierzu schweben noch.

Aufgrund dieses Eindruckes dr\"angt sich die Notwendigkeit einer Durchsuch\-ung der B\"uros des Institut de France und der angeschlossenen Akademien auf diese deutschfeindliche T\"atigkeit - teils vor dem Kriege, teils jetzt geradezu auf.

Eine erg\"anzende Haussuchung bei den jeweils am meisten belasteten Personen wird sich dar\"uber hinaus noch als notwendig erweisen. 

Es ist beabsichtigt, auch in dieser Frage eng mit dem zust\"andigen Kulturreferat der Deutschen Botschaft, Paris zusammen zu arbeiten.

Um Mitteilung dort evtl.\ bekannter Einzelheiten oder besonderer Anhalts\-punkte bei der Durchf\"uhrung dieser geplanten Massnahmen darf gebeten werden.

\smallskip
\noindent
SS-Obersturmbannf\"uhrer

\subsubsection*{English translation}

\noindent
To the Military Command in France

\smallskip
\noindent
Headquarters

\smallskip
\noindent
Paris, Hotel Majestic

\smallskip
\noindent
Re: Institut de France, Paris, and attached academies

\smallskip

A secret vote on Marshall P\'etain's policies and collaboration, which took place some time ago in the council of the \textit{Acad\'emie fran\c{c}aise}, resulted in 3 of the council members present voting for and 4 voting against P\'etain and against collaboration.  Voting negatively were:
\begin{itemize}
\item
Paul Val\'ery, writer
\item
G.\ Duhamel, writer
\item
Duke Maurice de Broglie, historian [sic]
\item
Emile Picard, adherent of the former French Popular Front [sic]
\end{itemize}
Voting in favor:  Andr\'e Bellessort, writer and permanent secretary of the \textit{Acad\'e\-mie fran\c{c}aise}, Cardinal Baudrillart, as well as the well known journalist Abel Bonnard.%
\footnote{Andr\'e Bellessort (born 1866) and Cardinal Baudrillart (born 1859) both died in 1942.  The official website of the French Academy asserts that Bellessort was permanent secretary too briefly for his openly collaborationist views to damage the Academy.  Abel Bonnard (1883--1968) lived in exile in Spain after the war.}

Not only in the \textit{Acad\'emie fran\c{c}aise}, but also in the \textit{Acad\'emie des Inscriptions et belles letters}, the \textit{Acad\'emie des Sciences}, the \textit{Acad\'emie des beaux arts}, and the \textit{Acad\'emie des sciences morales et politiques}, there have been and still are Jews, pro-Bolshevists, Freemasons, but also adherents of the \textit{Action Fran{\c c}aise} and members of associated groups.  In a number of cases, these scientists were also members of the \textit{Union des Intellectuels Fran{\c c}ais pour la Justice, la Libert\'e et la Paix}, which must be mentioned as the organization that united wide circles of French intellectuals on the sole basis of their absolute hostility to anything German.  It has been proven that Aim\'e Cotton and Charles Mauguin, members of the \textit{Acad\'emie des Sciences} who both belong to the \textit{Union des Intellectuels Fran{\c c}ais}, though officially it is no longer active, have agitated against the existing situation in occupied France and against Germany with illegal pamphlets such as \textit{l'Universit\'e libre}.  The matter is still under investigation.

This makes it absolutely necessary to search the offices of the \textit{Institut de France} and its associated Academies for information on this anti-German activity, partly before the war, but partly ongoing.

A complementary search of the homes of those most incriminated is also necessary.

On this issue also, we plan to work in close concert with the appropriate cultural officials of the German Embassy in Paris.

We request any details that may be known or any particular leads that should be followed in carrying out these planned measures.

\smallskip
\noindent
[signature]

\smallskip
\noindent
SS-Major

\subsection{SD re \textit{l'Acad\'emie des Sciences}, 13 May 1941}\label{subsec:SS_13_may_1941}

\textit{CARAN AJ/\-40/\-567.}

\subsubsection*{German original}

\noindent
SD, Dienststelle Paris

\noindent
Lagebericht Nr.\ 22/41

\noindent
f.d.Z. vom 1.-13.5.41

\smallskip

Die \underline{Acad\'emie des Sciences} ist in der Zusammensetzung ihrer Mitglieder etwas weniger politisch ausgepr\"agt.  Als die hervorragendsten und wissenschaft\-lich bedeutendsten Vertreter werden der Doktor Gosset, der Doktor Roussy und der Prince de Broglie bezeichnet. Einige andere Mitglieder, wie der Jude Hadamard, Verwandter von Dreyfuss [sic], sind zu den deutschfeindlichen Elementen zu z\"ahlen.

\subsubsection*{English translation}

\noindent
SD, Paris Office

\noindent
Situation report Nr.\ 22/41

\noindent
For the period of 1.-13.5.41

\smallskip

In the composition of its membership, the \textit{Acad\'emie des Sciences} has somewhat less political orientation.   Dr.\ Gosset, Dr.\ Roussy, and the Prince de Broglie%
\footnote{Louis de Broglie, who became Duke rather than Prince only after his brother's death.} 
are among its most outstanding and scientifically prominent members.  Some of the other members, such as the Jew Hadamard, a relative of Dreyfuss [sic], are to be counted among its anti-German elements.

\subsection{SS re \textit{l'Universit\'e libre}, 24 May 1941}\label{subsec:SS_24_may_1941}


\textit{This document, signed by Biederbick of the SS, is in CARAN AJ/\-40/\-567.  It is not dated, but it is stamped with the date 24 May 1941, presumably the date it was received by Group 4.}

\subsubsection*{German original}

\noindent
Re:  Die Agitation unter der Studenten von Paris und die T\"atigkeit einer pro-britannischen Organisation, genannt ``L'Universit\'e Libre'' (Die freie Universit\"at).

\smallskip

Man meldet die T\"atigkeit im besetzten Gebiet der gaullistischen und britannischen Propagandaorganisation, genannt ``L'Universit\'e Libre''.  Die Ver\-breit\-ung der von dieser Organisation herausgegeben im Abzugsverfahren hergestellten Wochenschrift geschieht in einer bedeutenden Anzahl im geschlossenen Umschlage und adressiert an Pers\"onlichkeiten, die sich f\"ur diese Sache interessieren k\"onnten und zur Mitarbeit bereit w\"aren.  Die haupts\"achlisten Pers\"onlichkeiten der Universit\"at, welche an der vogenannten Aktion mitarbeiten, sind alles Professoren, die der Freimaurerei angeh\"ort haben und ihre Beziehungen zu dieser nicht gel\"ost haben.  Sie werden eifrig von einer Anzahl israelitischer Professoren unterst\"utzt, die ihres Amtes bisher noch nicht enthoben worden sind.  Unter den T\"atigsten findet man an der Spitze der ``Universit\'e Libre'' folgende Herren:
\begin{itemize}
\item
Cotton, professeur \`a la Facult\'e des Sciences
\item
Tiffeneau, prof.\ \`a la Facult\'e de M\'edecine, ancien Doyen wohnh.: Paris, 85, Bv.\ St.\ Germain
\item
Mauguin, prof.\ \`a la Facult\'e des Sciences
\item
Joliot, prof.\ au Coll\`ege de France
\item
Sainte Lague [sic], Prof.\ au Conservatoire des Arts et M\'etiers
\item
Daniel Auger, Charg\'e de Recherches
\item
Albert Bayet, k\"urzlich abberufen.
\end{itemize}

Alle diese Professoren \"uben gegenw\"artig ihr Amt in Paris aus.  Eine schnell durchgef\"uhrte Nachforschung hat ergeben, dass die meisten gleichfalls an der Spitze von Gruppen stehen, welche f\"ur den Widerstand gegen die Zusammen\-arbeit und die Besetzung sind und von densen es eine Menge im lateinischen Viertel (Quartier Latin) gibt.

Die Organisation ``L'Universit\'e Libre'' disponiert, wie es bereits gesagt wurde, \"uber zahlreiche Korrespondenten in der nicht besetzten Zone, die ihrerseits sich bem\"uhen, dieselbe T\"atigkeit unter der Universit\"atsjugend der freien Zone auszu\"uben.  Man sagt, dass es die Gegenden von Lyon, Grenoble, Marseille, Montpellier und Toulouse sind, welche den Gegenstand der gr\"ossten T\"atigkeit von Seiten der Vertrauensm\"anner dieser Organisation bilden.  Man f\"uhrt inbesondere die nachstehenden Herren auf:
\begin{itemize}
\item
Jean Perrin in Lyon
\item
Paul L\'evy, prof.\ \`a l'Ecole Polytechnique de Lyon
\item
Hadamard in Toulouse
\item
Sampaix, Pierre, in Grenoble
\item
Masson, Bernard, in Marseille
\item
Luc, Directeur de l'Enseignement technique \`a Vichy
\end{itemize}

Die neuerliche Verbreitung von anglophilen Flugschriften unter der Universit\"atsjugend, die sich sowohl in der besetzten als auch in der freien Zone entwickelt hat, r\"uhrt von dieser Organisation her, wie man versichert.  Diese ist weit davon entfernt, sich mit der Arbeit zu begn\"ugen, und sie organisiert ausserdem im Quartier Latin und in den vorgenannten St\"adten Konferenzen, in denen die Universit\"atsjugend insbesondere ermutigt wird, gegen jede Orientierung der Regierung zu einer Politik der Zusammenarbeit mit dem Reiche Widerstand zu leisten und zu manifestieren.

\smallskip
Dr.\ Biederbick [signature by hand]

\smallskip
SS-Sturmbannf\"uhrer.

\subsubsection*{English translation}

\noindent
Re:  Agitation among Parisian students and the activity of a pro-British organization called ``l'Universit\'e Libre'' (Free University)

\smallskip

There are reports of activity in the occupied zone by a Gaullist and British propaganda organization called \textit{l'Universit\'e libre}.  A significant number of copies of a weekly mimeographed newsletter published by this organization are distributed in sealed envelopes, addressed to people who could be interested and ready to help. The most prominent university figures participating in this activity are all professors who adhered to Freemasonry and have not cut their links with it.  They are eagerly supported by a number of Jewish professors who have not yet been dismissed from their positions.  Among the most active leaders of the \textit{Universit\'e libre} are the following people: 
\begin{itemize}
\item
Cotton, professor at the Faculty of Sciences
\item
Tiffeneau, professor at the Faculty of Medicine, former dean.  Living at: Paris, 85 Bd.\ St.\ Germain
\item
Mauguin, professor at the Faculty of Sciences
\item
Joliot, professor at the College of France
\item
Sainte Lague [sic],%
\footnote{The French mathematician Andr\'e Sainte-Lagu\"e (188--1950) was a pioneer in graph theory.} 
professor at the Conservatoire des Arts et M\'etiers
\item
Daniel Auger, researcher
\item
Albert Bayet, recently recalled
\end{itemize}

All these professors presently work in Paris.  A quick investigation has revealed that most of them are also leaders of groups that are for resistance against collaboration and the occupation, and that many of these are in the Latin Quarter.

The organization \textit{l'Universit\'e libre} has, as was already said, numerous correspondents in the unoccupied zone, who for their part are trying to carry out the same activity among the university youth in the free zone.  It is said that trusted members of this organization are most active in the regions of Lyon, Grenoble, Marseille, Montpellier, and Toulouse.  The following men are especially mentioned:
\begin{itemize}
\item
Jean Perrin, in Lyon
\item
Paul L\'evy, professor at the Ecole Polytechnique of Lyon
\item
Hadamard in Toulouse
\item
Sampaix, Pierre, in Grenoble
\item
Masson, Bernard, in Marseille
\item
Luc, Director of Technical Education at Vichy
\end{itemize}

The recent dissemination of anglophilic pamphlets among the university youth, which has developed in the occupied zone just as in the free zone, comes from this organization, as has been verified.  The organization is far from being satisfied its work, and it also organizes conferences in the Latin Quarter and in the above-named cities, in which the university youth are especially encouraged to offer and demonstrate resistance against any orientation of the government to a policy of collaboration with the Reich.  

\smallskip
Dr.\ Biederbick [signature by hand]

\smallskip
SS-Major.

\subsection{Dahnke to Voize, 20 September 1941}\label{subsec:dahnke_voize}

\textit{This is the body of a letter on MBF letterhead, signed by Dahnke and dated 20 September 1941.  Now in CARAN Box AJ/\-40/\-563.}

\subsubsection*{German original}

\noindent
Sehr geehrter Herr Voize!

\smallskip
\noindent
Es ist mir nicht gelungen, Ihnen bis jetzt die Mittel f\"ur den Aufbau Ihrer Unterrichtsanstalt zu beschaffen.

\smallskip
\noindent
Ich werde die n\"achsten f\"unf Wochen abwesend sein und bitte Sie deshalb, sich Ende Oktober noch einmal mit mir in Verbindung zu setzen.

\subsubsection*{English translation}

\noindent
Very dear Mr.\ Voize!

\smallskip
\noindent
I have not yet been able to obtain for you the means for setting up your educational institute.

\smallskip
\noindent
I will be absent for the next five weeks.  So please get in touch with me once again at the end of October.

\subsection{German Navy re Borel et al., 25 October 1941}\label{subsec:navy}

\textit{From the office of the Commanding Admiral of the German Navy in France to the secret police section of the MBF.  Now in CARAN AJ/\-40/\-558.  We provide an English translation of the body of the letter in Section~\ref{sec:arrest}.}

\subsubsection*{German original}

\noindent
Geheim

\smallskip
\noindent
Abschrift

\smallskip
\noindent
Der Kommandierende Admiral in Frankreich

\noindent
B.Nr.38791 geh. A III c

\smallskip
\noindent
Paris, den 25. Okt. 1941

\smallskip
\noindent
An Milit\"arbefehlshaber in Frankreich Gruppe Ic

\smallskip
\noindent
Betr. : Verhaftung franz\"osischer Wissenschaftler

\smallskip
\noindent
Ohne Vorgang

\smallskip

Der Forschungsstab des Marinewaffenamtes im Oberkommando der Kriegsmarine arbeitet zurzeit in Paris an wichtingen kernphysikelischen Problemen zusammen mit dem Pariser Institut ``Curie''. Bei diesen Arbeiten sind die deutschen Wissenschaftler auf Zusammenarbeit mit franz\"osischen Wissenschaft\-lern angewissen.  U.a. sind hierbei zu nennen der Mathematiker Prof. Borell [sic], die Physiker Prof. Langevin und Cotton, der Christallograf Prof.\ Mauguin und der Mineraloge Prof.\ La Pique [sic], letztere beide von der Sorbonne.

Nach Mitteilung der Vertreter des OKM beim Kommandierenden Admiral Frankreich sind die vorgenannten franz\"osischen Wissenschaftler seit einiger Zeit verhaftet worden. Da sich unter diesen Umst\"anden die Zusammenarbeit zwi\-schen den deutschen und franz\"osischen Wissenschaftlern sehr schwierig gestaltet und die Durchf\"uhrung  milit\"arwissenschaftlich wichtiger Forschungsaufgaben vielleicht unm\"oglich gemacht wird, wird um Mitteilung geboten, ob die Verfehlung der verhafteten franz\"osischen Wissenschaftler so schwerwiegender Natur sind, dess die Verhaftungen aufrecht erhalten werden m\"ussen.

\smallskip
\noindent
F\"ur den Kommandieren Admiral in Frankreich

\smallskip
\noindent
Der Chef des Stabes

\smallskip
\noindent
I.V.

\smallskip
\noindent
gez. Unterschrift

\smallskip
\noindent
F.d.R.d.A.

\smallskip
\noindent
Rittmeister

\subsection{Group 4 re hostage taking, 27 October 1941}\label{subsec:hostage_27_10_1941}

\textit{In CARAN AJ/\-40/\-563.}

\subsubsection*{German original}

Bericht \"uber die Woche vom 20--26.10.1941

\smallskip

In dem Runderlass vom 18.9.1941 (V ju 821.1009.41 g),%
\footnote{``V ju'' refers to the section of the Verwaltungsabsteilung that supervised the French courts.}
\"uber die Geiselnahme haben die Studenten und Hochschullehrere als einzige Berufsgruppe ausdr\"uckliche Erw\"ahnung gefunden.  Es w\"are erw\"unscht, wenn V kult bei der Aufnahme von Hochschullehrern und Studenten in die Geisellisten beteiligt w\"urde.  Ein Teil der franz\"osischen Hochschullehrer hat sich in der Vergangenheit ausgesprochen deutschfeindlich bet\"atigt, insbesondere ist unter den franz\"osischen Naturwissenschaftlern eine linksradikale Haltung weit verbreitet gewesen und besteht zweifellos auch heute noch fort; insofern k\"onnen gewissen franz\"o\-sische Wissenschaftler mit Recht als geistige Urheber kommunistischer Umtriebe mit\-verantwortlichen gemacht werden.  Dasselbe gilt entsprechend f\"ur die chauvenis\-tische Haltung franz\"osischer Geisteswissenschaftler.  Es ist aber erforderlich, bei der Geiselnahme die wissenschaftliche Bedeutung der in Betracht kommenden franz\"osischen Gelehrten in Rechung zu stellen, weil die Erschiessung eines bedeutenden Wissenschaftlers als Geisel (nicht als T\"ater) im akademischen Bewusstsein zweifellos sehr lange nachwirkt.

Das frans\"osischen Unterrichtsministerium hat darauf aufmerksam gemacht, dass die Bekanntmachungen des Milit\"arsbefehlshabers in auffallender Weise an Schul-und Universit\"atseing\"angen angeschlagen w\"urden; es hat darum gebeten, die Bem\"uhungen des Rektors um die Aufrechterhaltung der Ruhe unter den Sch\"ulern und Studenten dadurch zu unterst\"utzen, dass von einer Anbringung der Plakate unmittelbar an den Eing\"angen der Universit\"atsinstitute und der Schulen abgesehen w\"urde.  V pol wird diesem Wunsche Rechnung tragen und dem Kommandanten von Gross-Paris, auf dessen Anordnung der Pr\"afekt die Plakate durch eine Plakatierungsgesellschaft anbringen l\"asst, entsprechende Anweisungen erteilen.  Uber die Stimmung in der Studentenschaft berichtet der Vertreter des franz\"osischen Unterrichtsministeriums, dass die deutschen Erfolge im Osten ihre Wirkung auf die franz\"osischen Studenten anscheinend nicht verfehlt h\"atten.

Die Ermittlungen deutschen Buchgutes, die Verzeichnung deutscher Autographe und die Photokopierung wertvoller Miniaturhandschriften in den Pariser Bibliotheken durch den Bibliotheksschutz wurden planm\"assig fortgef\"uhrt.

Oberregierungsrat Dr.\ Dahnke hat in der Berichtswoche seinen Dienst wie\-der angetreten.

Die B\"uror\"aume im Hause Av.\ Kleber 43--45 wurden in der Berichtswoche nicht geheizt.  Der Fortgang der Arbeiten wurde dadurch erbeblich beeintr\"achtigt.  Dem Vernehmen nach ist auch in Zukunft die regelm\"assige Heizung dieser R\"aume nicht gew\"ahrleistet.

\subsubsection*{English translation}

Report for the week of 20--26 October 1941

\smallskip

In the circular of 18 September 1941 (V ju 821.1009.41 g), concerning the taking of hostages, students and university professors were the only professional group explicitly mentioned. 	It would be desirable for Group 4 to be involved in the inclusion of names of students and professors in the hostage lists.  Part of the French university professoriat participated outspokenly in anti-German activity in the past, and especially among French natural scientists a radical leftist attitude was very widespread and doubtlessly continues today; to this extent certain French scientists can be rightly held responsible as a groupd as instigators of communist subversion.  The same is true of the chauvinistic attitude of French professors in the humanities.  But for the taking of hostages, it is necessary to take into account the scientific significance of the French professors being considered, because the shooting of an important scientist as a hostage (rather than for his actions) would doubtlessly have an affect on academic consciousness for a very long time.

The French Education Ministry has pointed out that the announcements of the MBF are conspicuously displayed at the entrances to schools and universities; it has asked that in order to support the efforts of the rector to maintain peace among the students, the placement of placards directly in the entrances of university institutes and schools be avoided.  The V pol section will take this request into account and will accordingly advise the Commander of Greater Paris, under whose order the prefects had the placards posted by a placarding company.  As for the mood of the student body, the representative of the French Education Ministry reports that the German successes in the East seem not to be failing to have an effect.

The investigation of German book properties, the registration of German manuscripts, and the photocopying of minature manuscripts in the Paris libraries by the Library Protection Service are continuing as planned.

Superior war administration adviser Dahnke returned to his post during the week of this report.

The office room in the house at 43--45 Kleber Avenue was not heated during the week of this report.  This seriously affected the continuation of the work.  As we understand it, regular heating of this room is also not assured in the future.

\subsection{Group 4 re Voize and Peyron's accusations, 1 November 1941}\label{subsec:dahnke_1_nov_1941}


\textit{In CARAN AJ/\-40/\-567.}

\subsubsection*{German original}

\noindent
Paris, den 1.11.1941

\smallskip
\noindent
Betrifft: Die freimaurerischen und bolschevistischen Kreise an der Sorbonne und im Akademiebezirk Paris.

\smallskip
\noindent
Sachbearbeiter: OKVR. Dr Dahnke

\smallskip
\noindent
1.) Vermerk:  Professor Peyron vom Institut Pasteur und der vom Minister Carcopino seines Amtes enthobene Professor vom Lyc\'ee Louis le Grand Voize, haben als Vertreter der freimaurerischen und bolschevistischen Einstellung unter den franz\"osischen Hochschullehren und Verwaltungsbeamten der Erziehungs\-wesens folgende Personen namhaft gemacht:
\begin{enumerate}\renewcommand{\theenumi}{\arabic{enumi})}
\item
G.\ Roussy, von Minister Jean Zay zum Rektor der Sorbonne ernannt; von Minister Ripert oder Minister Chevalier Ende 1940 seines Amtes als Rektor enthoben.  Nach Angabe des Professors Peyron (s.o.) auf dessen Betreiben.  Mit Halbj\"udin verheiratet.  Noch jetzt Direktor des Krebs-Instituts Paris.  Ende 1940 wurde Roussy die franz\"osische Staatsb\"urger\-schaft aberkannt, da er von Haus aus Schwiezer ist; Minister Carcopino hat seine Wieder einb\"urgerung veranlasst und wollte ihn wieder zum Rektor der Sorbonne machen.  Beweis:  Nach Angabe Peyrons hat Carcopino dies selbst dem Peyron mitgeteilt.

\smallskip
\noindent
Roussy ist nach der Schilderung Peyrons und Voizes in den freimaurersischen und bolschevisierenden Kreisen ausserordentlich einflussreich (vergl.\ de Zeitungsabschnitte bei diesem Aktenst\"uck).
\item
C.\ Luc [sic], Direktor des Enseignement Technique im franz\"osischen Unterrichtsministerium (Vorgang bei V kult 402) ebensfalls selbst nicht Frei\-mauer, aber ausgesprochener F\"orderer der Freimaurerei und der bolschevis\-ierenden Tendenzen in franz\"osischen Unterrichtswesen.
\item
Guyot, seit 20 jahren Generalsekret\"ar der Universit\"at Paris, von dieser Stellung zugleich mit Roussy seit 1940 entfernt, von Carcopino zum B\"uro\-chef des Enseignement sup\'erieur ernannt.  Freimaurer.
\item
Mme.\ Hattinghais [sic], von Carcopino zur Direktorin de jeunes Filles in S\`evres ernannt.  Bekannt als F\"orderin der Juden und Freimaurer.
\item
Professor Soretti [sic],%
\footnote{Ludovic Zoretti; see Section~\ref{sec:PicCha}.}
von Minister Chevalier seines Amtes als Professor an der Universit\"at Caen wegen seiner bolschevistischen Einstellung seines Amtes [sic] enthoben; Carcopino hat der Amstsenthobung in eine Versetzung in den Ruhestand umwandelt.
\item
Masbou, von dem directeur de l'enseignement, Rosset, als directeur de l'enseignement primaire de la Seine, eingesetzt; selbst nicht Freimaurer, aber F\"orderer derselben und bolschevistischer Tendenzen in der Lehrer\-schaft; von Chevalier seines Amtes enthoben von Carcopino als directeur der Ecole Normale sup\'erieure de l'enseignement technique wieder eingesetzt.
\item
Santelli, inspecteur d'Academie de la Seine et Marne, von Chevlaier wegen kommunistischen Einstellung abgesetzt; von Carcopino als inspecteur general l'enseignement technique wieder eingesetzt.
\item
Chatelun [sic], fr\"uher proviseur de lyc\'ee Louis-le-Grand, durch Carcopino zum directeur de l'enseignement primaire de la Seine en Stelle von Masbou ernannt.  Ist selbst nicht Freimaurer, sympathisiert aber mit diesem und hat deutschfeindliche Anschriften im lyc\'ee Louis-le-Grand geduldet.
\item
Prof.\ Fr\'ed\'eric Joliot, Prof.\ der Chemie am College de France.  Arbeitet auf Anordnung des Milit\"arbefehlshabers mit deutschen Wissenschaftlern in Paris zusammen und steht deswegen under deutschem Schutz.  Massnahmen gegen ihn sind nur nach F\"uhlungnahme mit dem Milit\"arbefehlshaber in Frankreich, Gruppe V kult m\"oglich.  Seine Frau: 
\item
Mme.\ Ir\`ene Joliot-Curie, Prof.\ der Chemie an der Sorbonne.  Beide links\-radikal und gegen die Zusammenarbeit mit Deutschland eingestellt.
\item
Professor D.\ H.\ Maguin [sic], Prof.\ der Mineralogie an der Sorbonne und an der Ecole pratique des Hautes Etudes.
\item
Professor Emil [sic] Cotton, Prof.\ der Physik an der Sorbonne und an der Ecole pratique des Hautes Etudes.
\item
Professor Louis Lapicque, Prof.\ der Physik an der Akademie der Wissenschaften.  Ebenso.
\item
Professor Emil [sic] Borel, Prof.\ der Mathematik an der Sorbonne.  Mitglied der Academie des Sciences.  Ebenso.
\item
Professor Gustav [sic] Monod, professeur de philosophie in Versailles.  Directeur de cabinet under dem j\"udischen Minister Jean Zay.  Ebenso.
\end{enumerate}

\subsubsection*{English translation}

\noindent
Re: Freemason and Bolshevist circles in the Sorbonne and the Paris academic district

\smallskip
\noindent
Expert: OKVR\footnote{Oberkriegsverwaltungsrat (Superior War Administration Advisor)} Dr.\ Dahnke

\smallskip
\noindent
Remarks:  Professor Peyron, of the Institut Pasteur, and Professor Voize, who was removed from his position at the Louis-le-Grand lyc\'ee by Minister Carcopino, have named the following persons as representatives of the Freemason and Bolshevist view in higher education and the administration of public education: 
\begin{enumerate}\renewcommand{\theenumi}{\arabic{enumi})}
\item
G.\ Roussy, who was named rector of the Sorbonne by Minister Jean Zay and removed from the position at the end of 1940 by Minister Ripert 
or Minister Chevalier.  Professor Peyron said the removal was on his initiative.  Married to a half-Jew.  Still today director of the cancer institute in Paris.%
\footnote{This institute, not to be confused with the scientifically much more important \textit{Institut Pasteur}, was one of many regional institutes for cancer in France.  It was founded by Roussy in 1926 and rechristened the \textit{Institut Gustave Roussy} after the war.}
At the end of 1940, Roussy, was stripped of his French citizenship, as he is originally from Switzerland.  Minister Carcopino has arranged for him to regain his citizenship and wants to make him rector of the Sorbonne again.  Proof:  According to Peyron, Carcopino himself told him this.

\smallskip
\noindent
According to Peyron and Voize's accusation, Roussy is exceptionally influential in Freemason and Bolshevist circles (see the attached newspaper clipping).
\item
C.\ Luc [sic],%
\footnote{Hippolyte Luc.} 
director of technical education in the French Ministry of Public Instruction (file in V kult 402), though not himself a Freemason, is 
a decided supporter of Freemasonry and Bolshevist tendencies in French education.
\item
Guyot, for 20 years general secretary of Paris University, removed from this position along with Roussy at the end of 1940 
and named head of the office of higher education by Carcopino. Freemason.
\item
Mme.\ Hattinghais [sic],%
\footnote{Edm\'ee Hatinguais.} 
named by Carcopino director of the \textit{Ecole Normale} for young women in S\`evres. Known as supporter of Jews and Freemasons.
\item
Professor Soretti [sic],%
\footnote{Ludovic Zoretti; see Section~\ref{sec:PicCha}.}
removed by Minister Chevalier from his position as professor at the University of Caen because of his 
Bolshevist opinions; Carcopino has changed his dismissal to retirement.  
\item
Masbou, named director of primary education in the Seine by education director Rosset.  Not himself a Freemason but a supporter of Freemason and Bolshevist tendencies in the teaching corps; removed from his functions by Chevalier, rehired by Carcopino as director of the \textit{Ecole Normale sup\'erieure} for technical education.
\item
Santelli, inspector of the Academy of the Seine and Marne, removed by Chevalier because of his communist opinions, rehired as general inspector of technical schools by Carcopino.
\item
Chatelun [sic],%
\footnote{Lucien Chattelun.} 
former director of the Louis-le-Grand lyc\'ee, named by Carcopino director of primary education of the Seine, replacing Masbou.  Not himself a Freemason but sympathizes with them and has permitted anti-German placards at the Louis-le-Grand lyc\'ee.
\item
Prof.\ Fr\'ed\'eric Joliot, professor of chemistry at the College of France.  Works with German scientists in Paris by arrangement with the military command and is therefore under German protection.  Measures against him are only possible after checking with the military command in France, group V culture.  His wife: 
\item
Mme.\ Ir\`ene Joliot-Curie, professor of chemistry at the Sorbonne.  Both are radical leftists and are engaged against collaboration with Germany.
\item
Professor D.\ H.\ Maguin [sic], professor of mineralogy at the Sorbonne and at the Ecole pratique des hautes \'etudes.
\item
Professor Emil [sic] Cotton,%
\footnote{There is a confusion here between Aim\'e Cotton, the physicist who was a member of the \textit{Acad\'emie des Sciences} and was arrested along with Borel, and his younger brother, the mathematician Emile Cotton, who was not yet a member of the \textit{Acad\'emie des Sciences} in 1941.}
professor of physics at the Sorbonne and at the Ecole pratique des hautes \'etudes.
\item
Professor Louis Lapicque, professor of physiology at the \textit{Acad\'emie des Sciences}. The same.
\item
Professor Emil [sic] Borel, professor of mathematics at the Sorbonne. Member of the \textit{Acad\'emie des Sciences}.  The same.
\item
Professor Gustav [sic] Monod, professor of philosophy at Versailles. Cabinet director for the Jewish minister Jean Zay. The same.
\end{enumerate}

\subsection{Group 4 re arrest of Borel et al., 10 November 1941}\label{subsec:dahnke_10_nov_1941}

\textit{In CARAN AJ/\-40/\-567.}

\subsubsection*{German original}

\noindent
\underline{Betreff:}  Die freimaurerischen und bolschevistischen Kreise an der Sorbonne und im Akademiebezirk Paris.

\smallskip
\noindent
\underline{Sachbearbeiter:}  OKVR.\ Dr.\ Dahnke

\smallskip
\noindent
1.  Vermerk:  Die Alst hat die Professoren:  Mauguin, Cotton, Lapicque und Borel (vergl.\ anliegenden Vermerk vom 1.11.41), sowie die Studienr\"ate Aubert und Cazalas (vergl.\ anliegende Eingabe) verhaftet.  Ich habe mich mit der Alst (Major Dr.\ Reile) in Verbindung gesetzt und mit dem Sachbearbeiter Hauptmann Kr\"ull%
\footnote{Here, and in the next paragraph, the name ``Sch\"urr'' was typed and crossed out, with ``Kr\"ull'' or ``Kr\"ulle'' written in by hand.}
gesprochen und ihm die Aufzeichnung vom 1.11.41 \"ubergeben, um ihn auf die Zusammenh\"ange aufmerksam zu machen, die ihm unbekannt waren.  Er will sein Ermittlungen auf diesen Kreis erstrecken.

\smallskip
\noindent
2.) \underline{Wieder vor zum 24.11.1941} (erneute F\"uhlungsnahme mit Hauptmann Kr\"ull, um sich \"uber seine Ermittlungsergebnise zu unterrichten).

\medskip
\noindent
\textit{Handwritten note, difficult to decipher, initialed by Dahnke and dated 18 November:}  Nach Mitteilung de Roy sind Borel, Cotton \& Mauguin inzwischen wieder entlassen.

\subsubsection*{English translation}

\noindent
Re: Freemason and Bolshevist circles in the Sorbonne and the Paris academic district

\smallskip
\noindent
Expert: OKVR Dr.\ Dahnke

\smallskip
\noindent
1- Note: The Alst arrested professors Mauguin, Cotton, Lapicque and Borel (see the attached note of 1 November 41), and lyc\'ee teachers Aubert and Cazalas (see the attached memorandum). I have gotten in contact with the Alst (Major Dr.\ Reile), spoken with the expert in charge of the file, Captain Kr\"ull, and transmitted to him the note of 1 November 41 in order to bring connections of which he had been unaware to his attention.  He intends to extend his investigations to this circle.

\smallskip
\noindent
[added later] 2- Once again, on 24 November 41 (new contact with Captain Kr\"ull, informing me about the results of his investigations)

\medskip
\noindent
\textit{Handwritten note, difficult to decipher, initialed by Dahnke and dated 18 November:}  In the meantime, Borel, Cotton \& Mauguin released after a note from Roy.

\subsection{Group 4 re release of Borel et al., 25 November 1941}\label{subsec:dahnke_25nov41}


\textit{This document, in CARAN AJ/\-40/\-567, is dated November 1941, with no day of the month indicated, but it was evidently written on 25 November or later.}

\subsubsection*{German original}

\noindent
Der Milit\"arbefehlshaber in Frankreich 

Verwaltungsstab Abteilung Verwaltung

\smallskip
\noindent
Paris, am ... Nov 1941

\smallskip
\noindent
\underline{Betreff:}  Die freimaurerischen und bolschevistischen Kreise an der Sorbonne und im Akademiebezirk Paris.

\smallskip
\noindent
\underline{Sachbearbeiter:} OKVR. Dr Dahnke

\smallskip

1- Vermerk : R\"ucksprache mit Hauptmann Kr\"ull am 25 11 41. Er hat die Professoren: Mauguin, Cotton, Lapicque und Borel s\"amtlich wieder aus der Haft entlassen; ihre Vernahmung hat ergeben, dass sie sich alle, insbesondere aber Cotton aufrichtig zu den von ihnen vor und w\"ahrend des Krieges vertretenen politischen Ideen noch heute bekennen. Sie haben offen zum Ausdruck gebracht, dass sie das politischen- System Deutschlands nicht billigen und dass sie von dem Siege Englands und Amerikas die Rettung Frankreichs erwarte. Sie haben jedoch nachdr\"ucklich in Abrede gestellt, diese Gesinnung in irgendeine Weise, insbesondere unter den Studenten, bet\"atigt zu haben. Die Abwehr ist nich in der Lager gewesen, ihnen eine solche Bet\"atigung nachzuweisen, obwohl die V-M\"anner eine derartige Bet\"atigung behauptet hatten. Insbesondere war eine Gegen\"uberstellung mit Studenten, die aus dem Kreise dieser Professoren politisch beeinflusst sein sollten, nicht m\"oglich weil von den V-M\"annern solche nicht namhaft gemacht wurden.

In einer Besprechung mit dem Sachbearbeiter, Dr Epting und Dr Biederbick ist erwogen worden, die vier Professoren in derselben Weise ausserhalb Paris unter Polizeiaufsicht zu stellen, wie dies mit Langevin in Troyes geshehen ist. Angesichts einer derartigen politischen Haltung der vier Professoren, die von den rechtsgerichteten Kreisen - vergl. die Vorg\"ange in der Presse - als Zentrum der extremen links - und gegen Deutschland gerichteten Tendenzen an der Sorbonne bezeichnet wird, ist nicht anzunehmen, dass sie sich jeder politischen Meinungs\"ausserung im Kreise der franz\"osischen Wissenschaft enthalten werden. Dies rechtfertigt eine derartige Massnahme, auch wenn eine solche politische Bet\"atigung nich durch Zeugenaussagen nachgewiesen werden kann; auch Langevin ist seinerzeit eine aktive Bet\"atigung nicht nachgewiesen worden. Es muss jedoch zuvor die Auswirkung einer derartigen Massnahme wohl erwogen werden; ich habe deshalb Dr. Epting gebeten, die Sache zun\"achst in der Botschaft selbst zur Sprache zu bringen.

2.) Wieder vor zum 8.12.1941 (Besprechung mit Dr.Epting und Dr. Biederbick am 9.12.1941, 12 Uhr).

\subsubsection*{English translation}

\noindent
Re: Freemason and Bolshevist circles in the Sorbonne and the Paris academic district

\smallskip
\noindent
Expert: OKVR Dr.\ Dahnke

\smallskip
\noindent
1- Note:  Consultation with Captain Kr\"ull on 25 November 41.  He released Professors Mauguin, Cotton, Lapicque and Borel.  Their interrogation showed that all of them, especially Cotton, still candidly stand by the political ideas they advocated before and during the war.  The openly declared that do not approve of the German political system and that they expect France to be rescued by England and America's victory.  But they emphatically denied that they had in any way acted on their views, especially with students.  The intelligence service is not in position to prove such activity, though our informers claim it has taken place.  In particular, it is impossible to arrange a confrontation with students from the circles these professors were supposed to have influenced, because the informers did not identify any such students by name.

During a consultation between this expert, Dr.\ Epting, and Dr.\ Biederbick, we considered putting the four professors under police surveillance outside Paris, as was done with Langevin in Troyes.  In view of such a political attitude on the part of the four professors, who are described by rightist circles as the center of extreme leftist and anti-German tendencies in the Sorbonne --- see the references in the newspapers --- it cannot be expected that they will refrain from expressing their opinions in French scientific circles.  This would justify a measure of this kind, even if it is impossible to produce witnesses to prove their political activity.  No active participation by Langevin was proven either.  But we must deliberate carefully before implementing such a measure; I have therefore asked Dr.\ Epting to personally discuss the matter at the Embassy as soon as possible. 

\smallskip
\noindent
2- Back on 8 December 41 (conversation with Dr.\ Epting and Dr.\ Biederbick on 9 December 41 at 12 o'clock).

\subsection{Group 4 re Mauguin, 11 February 1942}\label{subsec:mauguin_11_2_1942}

\textit{In CARAN AJ/\-40/\-567.}

\subsubsection*{German original}

\noindent
1) An Gruppe V pol \underline{im Hause}.%
\footnote{Group 4 of the Verwaltungsabteilung was sometimes designated as ``Gruppe V kult''.  This memorandum, prepared by Group 4, was addressed to Group 2, or ``Gruppe V pol'', responsible for the police.}

\smallskip
\noindent
\underline{Betrifft:}  Kommunistische Propaganda bei den Angeh\"origen der Universit\"at.

\smallskip
Von Professor Mauguin ist von seinen politisch rechtsstehenden Gegnern schon fr\"uher behauptet worden, dass er kommunistischen Kreisen nahest\"unde.  Tats\"acheliches Material, durch das diese Behauptung gest\"utzt w\"urde, leigt jedoch bis jetzt nicht vor.  Im Einverhehmen mit den Sicherheitsdienst wird die Angelegenheit weiter \"uberwatcht.  Zur Zeit ist kenie Handhabe gegeben, um von Seiten des Milit\"arbefehlshabers etwas zu unternehmen.  Von etwaigen weiteren Ermittliungsergebnisse wird V pol verst\"andigt werden.

\smallskip
\noindent
2) Wiedervorlage nach Abgang (Besprechung mit Dr.\ Biederbick.)

\medskip
\noindent
\textit{Handwritten note, difficult to decipher:}  1941 Okt ergriffen, 13 Nov entlassen

\subsubsection*{English translation}

\noindent
1) To Group V pol \underline{in the MBF}.

\smallskip
\noindent
\underline{Subject:}  Communist propaganda directed to members of the University.

\smallskip
Professor Mauguin's right-wing political opponents had already earlier declared that he stood close to communist circles.  Actual evidence to support this declaration has until now not yet been produced.  The matter will be further monitored in cooperation with the security service (SD).  So far the MBF has no grounds for doing anything. Group V pol will be informed of any results from further investigation.

\smallskip
\noindent
2) Review after departure (discussion with Dr.\ Biederbick.)

\medskip
\noindent
\textit{Handwritten note, difficult to decipher:}  1941 October apprehended, 13 November released

\section{Some French texts}\label{sec:french}

In most cases, we provide both the original and a translation into English.

\subsection{On the dismissal of Raymond Voize, 7 August 1941}\label{subsec:voize1}

\textit{This article appeared unsigned in \textit{l'Appel}.}

\subsubsection*{French original}

\noindent
Justice universitaire!

\smallskip

Nous sommes les premiers \`a annoncer que M.\ R.\ Voize, professeur d'Allemand au lyc\'ee Louis le Grand, mutil\'e de la grande guerre \`a 75\%, officier de la L\'egion d'honneur \`a titre militaire, a \'et\'e mis \`a la retraite d'office, sur l'initiative de M.\ Carcopino.

Malgr\'e plusieurs d\'emarches du doyen des professeurs du Lyc\'ee Louis-le-Grand et de diverses personnalit\'es, dont M.\ Lafont, ancien directeur de la Famille fran\c caise au minist\`ere de la Sant\'e, malgr\'e des notes exceptionnelles, la d\'ecision minist\'erielle est devenue irr\'evocable, gr\^ace \`a l'intol\'erance du repr\'esen\-tant de M.\ Carcopino, \`a Paris, un certain M.\ Verrier, qui s'est refus\'e \`a ouvrir le dossier de l'honorable professeur. 

En fait, le repr\'esentant de M.\ Carcopino est fix\'e sur toute cette affaire.  Il sait que l'on ne peut rien reprocher au professeur Voize.  Il s'en tire donc par des arguit\'es aussi habile que d\'eshonnorantes.  Les voici!

D'\'education catholique, M.\ Voize, apr\`es avoir connu une crise de conscience et d'agnosticisme, a retrouv\'e la voie de son enfance et s'est distingu\'e depuis nombre d'ann\'ees comme un catholique militant (nous en demandons pardon \`a M.\ de la Fouchardi\`ere).  C'est ainsi que, dans les premi\`
eres mois de l'ann\'ee 1939, M.\ Voize avait publi\'e une brochure \'educative o\`u il ne cachait pas ses opinions fond\'ees et respect\'ees.

Des mouchards intervinrent.

Ceux-ci alert\`erent de tr\`es veilles connaisances:  l'ignoble Jean Zay et son chef de ``cabinet'' Abraham, le Suisse Roussy -- mari\'e \`a la juive Thompson -- recteur naturalis\'e de l'Universit\'e et propri\'etaire de la Farine Nestl\'e, agent provocateur du Front populaire Quartier Latin, le F.\ M$\therefore$ Guyot,%
\footnote{A pyrimid of three dots was sometimes used by Freemasons, especially in France, in place of a period in initials.}
ci-devant  secr\'etaire g\'en\'eral de la Sorbonne, et Gustave Monod, Inspecteur g\'en\'eral de l'Acad\'emie de Paris.

L'affaire \'etait belle et bonne.  Et le repr\'esantant de M.\ Carcopino, le d\'ej\^a nomm\'e Verrier, qui se dit ``FILS DE QUARANTE HUITARD ET ATTACHE AUX IMMORTELS PRINCIPES DE 1789'', en a profit\'e pour vider la petite querelle, ajoutant ``QU'IL ENTENDAIT EPURER L'UNIVERSTE DE TOUS LES ELEMENTS REACTIONNAIRES''.

On retrouva m\^eme une note \'etablie en 1930 par le feu Charl\'ety, o\`u M.\ Voize \'etait accus\'e d'avoir mal administr\'e une fondation ext\'erieure \`a l'Universit\'e\dots  . Rien de plus et rien de moins!

Mieux encore M.\ Verrier a su rappeler \`a M.\ Carcopino qu'il avait soumis, \'etant encore recteur, \`a la signature de M.\ Chevalier, alors ministre, en date du 22 f\'evier 1941, un arr\^et\'e de mise \`a la retraite d'office de M.\ Voize\dots

Car M.\  Voize inqu\'etait l'Universit\'e! Et pourquoi? Parce que M.\ Voize \'etait "SUSPECT DE SYMPATHIE POUR L'ALLEMAGNE !". Ce sont les termes du rapport.

Cette d\'ecision prouve \'evidemment que M.\ Carcopino est un humaniste fort prompt et fort passionn\'e. Cette d\'ecision prouve encore - nous voulons le croire - que M.\ Carcopino a \'et\'e mal renseign\'e. Car cette mesure, \`a la fois faite pour r\'ejouir les Juifs et la s\'equelle gaulliste, n'est certainement pas dans la vocation spirituelle du romanisant illustre qu'est M.\ Carcopino.

Il est donc probable que M.\ Carcopino voudra bien mettre fin \`a ces d\'eloyaux proc\'ed\'es. On r\'ehabilite - ou presque - un r\'epugnant Zoretti, chef de ``nervis'' communistes, ex-doyen de la facult\'e des Sciences de Caen, personnage aussi suspect que... riche ! Et l'on tente de frapper un patriote intelligent, qui sait que la France a \'et\'e vaincue par des Zoretti et leurs laudateurs masqu\'es.

\subsubsection*{English translation}

We are the first to announce that Mr.\ R.\ Voize, professor of German at the lyc\'ee Louis le Grand, 75\% disabled in the great war, officer of the Legion d'honneur for his military service, has been retired on Mr.\ Carcopino's initiative.  

Despite many remonstrances by the dean of the professors at the lyc\'ee Louis le Grand and many others, including Mr.\ Lafont, former director of the French Family at the health ministry, and despite exceptional evaluations, the ministerial decision has been made final, thanks to the intolerance of Mr.\ Carcopino's representative in Paris, a certain Mr.\ Verrier, who has refused to reopen the file on this honorable professor.

In fact, Mr.\ Carcopino's representative is all over this case.  He knows that Professor Voize cannot be reproached for anything.  He gets around this by stunts that are as clever as they are dishonorable.  Here they are!

Catholic by education, Mr.\ Voize found the path of his childhood again after a crisis of conscience and agnosticism, and he has distinguished himself for many years as a militant Catholic (we ask pardon of Mr.\ de la Fouchardi\`ere).  Thus it was that in the summer of 1939, Mr.\ Voize had published an educational pamphlet in which he did not hide his well founded and respected opinions.

Then the informers intervened.

The alerted some very old acquaintances:  the vile Jean Zay and his chief of ``staff'' Abraham, the Swiss Roussy -- married to the Jewess Thompson -- naturalized rector of the University and proprietor of Nestl\'e Flour, \textit{agent provocateur} for the Popular Front in the Latin Quarter, the Freemason Guyot, previously general secretary of the Sorbonne, and Gustave Monod, Inspector General of the Paris Academy.

It was a good scandal.  And Mr.\ Carcopino's representative, the Verrier already mentioned, who calls himself ``SON OF A FORTY-EIGHTER AND ATTACHED TO THE IMMORTAL PRINCIPLES OF 1789'', took advantage of it to settle the little quarrel, adding ``THAT HE INTENDED TO PURIFY THE UNIVERSITY OF ALL THE REACTIONARY ELEMENTS''.

They even found a note set down in 1930 by the late Charl\'ety, where Mr.\ Voize had been accused of have poorly administered a foundation outside the university\dots.  Nothing more and nothing less!

Better yet, Mr.\ Verrier was able to remind Mr.\ Carcopino that while he was still rector, he had submitted for the signature of Mr.\ Chevalier, then minister, a decree dated 22 February 1941 that would retire Mr.\ Voize from office\dots

Because Mr.\ Voize disturbed the university!  And why?  Because Mr.\ Voize was ``SUSPECTED OF SYMPATHY FOR GERMANY!''.  This is the phrasing of the report.  

This decision obviously proves that Mr.\ Carcopino is a very diligent and very passionate humanist.  It also proves -- so we want to believe -- that Mr.\ Carcopino was misinformed.  Because this decision, made to please both the Jews and the Gaullist aftereffect, is certainly not in the spirit of a distinguished Latinizer like Mr.\ Carcopino.  So it is quite likely that Mr.\ Carcopino will indeed want to put an end to these disloyal processes.  They are rehabilitating -- or almost -- a repugnant Zoretti, chief of communist assassins, former dean of the Faculty of Sciences at Caen, a personality as suspicious as he is\dots rich!  And they try to punish an intelligent patriot, who knows that France has been conquered by the Zorettis and their masked extollers.%
\footnote{We provide some information about the mathematician Ludovic Zoretti in Section~\ref{sec:PicCha}.}

\subsection{Voize's proposal for commission, 7 August 1941}\label{subsec:voize2}

\textit{This article appeared in the weekly \textit{La Gerbe}.}

\subsubsection*{French original}

\noindent
Les rapports intellectuels franco-allemands:  Pour la cr\'eation d'un Haut Commissariat

\smallskip

\noindent
par R.\ VOIZE, agr\'eg\'e de l'Universit\'e, professeur au lyc\'ee Louis-le-Grand

\smallskip

Donc, Carcopino, ministre de l'Education Nationale est en train de s'effondrer.

Il s'\'ecroule sans retentissement, honteusement, dans l'\'etat de d\'ecomposition o\`u l'ont amen\'e ses manquements es ses f\'elonies: une chute mate.  Manquements et f\'elonies qui l'avaient d\'enounc\'e d'eux-m\^emes au Mar\'echal, \`a l'amiral, \`a l'ambassadeur repr\'esentant le gouvernement franc\c{c}ais en zone occup\'ees, et aux secr\'etaires g\'en\'eraux de la pr\'esidence du Conseil.

Il a \'et\'e achev\'e gr\^ace \`a l'action d\'etermin\'ee, concert\'ee, implacable et justici\`ere d'universitaires clairvoyants et r\'esolus, anim\'es de l'esprit du Mar\'echal, adeptes de la R\'evolution nationale, ayant constitu\'e ``l'\'equipe de reconstruction universitaire'', qui lui ont arrach\'e le masque aux rubans coup\'es, qu'il plaquait encore sur sa figure, de ses mains crisp\'ees. 

Janus bifrons, il protestait de son esprit de R\'evolution nationale et de sa bonne volont\'e de collaboration:  et il se courvrait du c\^ot\'e du Front populaire, dont il escomptait sournoisement le retour, et les r\'ecompenses.

Il donnait des gages.

Luc \'etait maintenu par lui dans sa grasse pr\'ebende de directeur g\'en\'eral de l'Engsignement technique, et install\'e en outre comme contr\^oleur et briseur d'\'elan de Lamirand.  (Ici m\^eme on a cri\'e alerte:  voir \textit{La Gerbe} du 3 juillet).

Pis encore:  Il s'appr\^etait \`a r\'einstaller, apr\`es une renaturalisation effectu\'ee sous ses auspices, dans le haut poste de recteur de l'universit\'e de Paris, Roussy le m\'et\`eque que nos courageux universitaires nationaux avaient fait jeter \`a la porte de la Sorbonne cet hiver. Roussy, Suisse d'origine, qui apportait de Suisse sa camelote de lait condens\'e, en m\^eme temps que des tracts bolcheviques, Roussy ce ploutocrate d\'emagogue (Nestl\'e lui assure un million par an) qui, recteur, se faisait photographier part \textit{l'Humanit\'e}, le poing lev\'e, \`a c\^ot\'e de Vaillant-Couturier:  scientifique ind\'esirable de pacotille, qui a peupl\'e l'Universit\'e de m\'ediocres \`a son image.  On va le renvoyer en Suisse, \`a ses vaches \`a lait, ou le maintenir en France, dans un camp qui ne sera pas un camp de Jeunesse. 

Roussy avait un f\'eal servant:  Guyot.  Carcopino l'a rescap\'e.  Guyot, ex-secr\'etaire g\'en\'eral de l'universit\'e de Paris, avorton perfide, larve mal\'efique, le d\'etenteur de tous les sales secrets des recteurs qui, depuis vingt ans, ont d\'esorient\'e moralement la jeunesse fran\c{c}aise.  Ce Guyot, chass\'e de la Sorbonne en m\^eme temps que Roussy, a \'et\'e rep\^ech\'e par Carcopino, et il se livre de nouveau \`a sa dissolvante activit\'e\dots comme chef de bureau de l'Engseignement sup\'erieur.

Reculant des limites de l'inf\^amie, Carcopino allait profiter de la torpeur universitaire du mois d'ao\^ut pour effectuer, dans tous les ordres d'enseignement, des nominations et des promotions massives de professeurs et d'administrateurs anim\'es d'une hostilit\'e syst\'ematique et irr\'eductible \`a l'\'egard du gouvernement de R\'evolution nationale et de collaboration du Mar\'echal.  C'est alors que l'\'equipe de reconstruction universitaire, l'\'equipe des hommes du Mar\'echal, se dressa d'un bloc:  le malfaiteur de l'Universit\'e fut stopp\'e net; et on le cassa aux gages.

Son activit\'e n\'efaste peut se r\'esumer ainsi:  il a cauteleusement soutenu, encourag\'e, dans l'Universit\'e, tous ceux qui excitaient stupidement ou criminellement la jeuenesse \`a fronder l'autorit\'e allemande et \`a se proclamer gaulliste; quant \`a ceux qui tentaient d'orienter cette jeunesse vers la voie de l'avenir indiqu\'ee par le Mar\'echal et l'amiral, Carcopino les a d\'esavou\'es, brim\'es, pers\'ecut\'es, \'elimin\'es.  Bref, il s'est appliqu\'e \`a parquer nos jeunes gens dans une esp\`ece de camp de concentration intellectuel ferm\'e \`a l'id\'eal comme au r\'eel, et il leur a interdit toute \'echapp\'ee, tout \'elan, toute tentative de compr\'ehension et de collaboration.

Laissons de c\^ot\'e, dans les bas c\^ot\'es, ce Bonaparte en carton sans p\^ate , cet histrion avantageux qui jouait les imperators et les traitres de tragi-com\'edie, ce macaque chamarr\'e au grima\c{c}ant sourire.  Qu'il aille rejoindre hors du Temple, dans les catacombes, la poussi\`ere de ses dieux morts!  Ne nous acharnons pas sur ce cadavre, emp\^echons-le simplement de se faire embaumer dans le sarcophage de l'Ecole normale sup\'erieure; il est \`a jamais indigne de guider la jeunesse, ce tenant de la basse maffia sectairement antireligieuse, d\'emagogue, bolchevisante et profiteuse, dont le chor\'ephore fut le derviche tourneur Durckheim et dont l'aboyeur tr\`es amplement patent\'e fut le butor intellectuel Bougl\'e.

\begin{center}
* * *
\end{center}

Et parlons constructivement:  pour que votre r\`egne arrive, collaboration, salut de la France et de l'Europe.

La collaboration a pour base la compr\'ehension, et pour but, la paix.

Des efforts patients, tenaces et sereins, sont d\'eploy\'es depuis un an, depuis l'exode, par les Allemands dans tous les domaines:  assistance \'emouvantes aux hordes mis\'erables des populations fran\c{c}aises (regardez sur les sentiers de la d\'eroute); travail procur\'e en France et en Allemagne aux Fran\c{c}ais et Fran\c{c}aises qui ont \`a assurer leur subistance et celle de leur famille; mansu\'etude des torts \`a l'\'egard des vaincus fr\'emissants\dots

Limitons, dans cet article, notre examen aux efforts de collaboration sur le plan intellectuel.  C'est sur le plan intellectuel que la collaboration pourra pr\'esenter les moins grandes difficult\'es:  car, selon le mot de Leibniz, ``les esprits sont les corps qui s'empeschent le moins.''  Et puis, nous croyons \`a la valeur pr\'edominante de l'\'elite des \'elites:  les \'elites sont faites pour diriger, et les masses pour \^etre digir\'ees par les \'elites.  Et les \'elites, chez qui le r\'eflexe est surmont\'e par la r\'eflexion, sont plus particuli\`erement aptes \`a comprendre que vivre, c'est surmonter le pass\'e, un pass\'e hypoth\'equ\'e de toutes nos carences, nos veuleries et nos abdications.

Le rapprochement des \'elites, les Allemands y travaillent.  L'Institut allemand, magistralement dirig\'e par le docteur Epting et son adjoint le docteur Bremer, met les Fran\c{c}ais de bonne volont\'e en contact avec la langue et la civilisation allemandes:  entreprise couronn\'ee de succ\`es, aux r\'epercussions incalculables.  Par ailleurs, des Allemands eminents, savants et techniciens, viennent prendre contact avec les milieux fran\c{c}ais correspondant \`a leurs sp\'ecialit\'es:  recemment encore, le ministre de l'Agriculture du Reich.

Que font les Fran\c{c}ais?  Les Fran\c{c}ais agissent.  Le groupe Collaboration,%
\footnote{A footnote gives the address:  26, rue Bassano, Paris.  T\'el.: Kl\'eber 71-14.}
avec ses diverses sections si fournies et si actives; l'Institut d'\'etudes germaniques de la Sorbonne; les diff\'erents groupements et journaux qui se r\'eclament du Mar\'echal et sont anim\'es de l'esprit de la R\'evolution nationale, les cahiers franco-allemands.

Donc, les Allemands agissent, des Fran\c{c}ais agissent.

\begin{center}
* * *
\end{center}

Mais ce qui manque, c'est un \textit{organisme de coordination} de tous ces efforts \'emouvants. 

Seul, un organisme officiel aurait les possibilit\'es de r\'ealiser cette indispensable, cette urgente, coordination:  parce que, seul, il appara\^itrait, aux yeux de tous, comme parfaitement qualifi\'e.

Cet organisme de coordination pourrait rev\^etir la forme d'un \textit{haut commissariat aux rapports intellectuels franco-allemands}. 

Il devrait \^etre rattach\'e directement \`a la pr\'esidence du Conseil, constituant ainsi un organisme super-minist\'eriel, de fa\c{c}on \`a imposer son autorit\'e aux divers minist\`eres, notamment au minist\`ere de l'Education nationale, et \`a \'edicter des mesures qui vaudraient \`a la fois pour la zone non occup\'ee et la zone occup\'ee.

S'il avait exist\'e, il e\^ut ordonn\'e au ministre de l'Education nationale de profiter de l'occasion des distributions des prix pour faire tenir \`a de multiples jeunes gens, class\'es parmi les meilleurs, les discours du Mar\'echal, de l'amiral et de M.\ Lamirand, o\`u sont expos\'eees les raisons de la collaboration; ainsi que des ouvrages fran\c{c}ais, our traduits en fran\c{c}ais, susceptibles de les mettre sur la voie d'une compr\'ehension de l'Allemagne et d'une collaboration des c\oe urs. 

S'il e\^ut exist\'e, l'on aurait peut-\^etre essay\'e d'utiliser la p\'eriode des vacances d'\'et\'e our tenter d'\'etablir un contact entre jeunesses fran\c{c}aise et allemande.

En tout cas, l'on peut s'attacher \`a organiser des contacts professionels de plus en plus larges.

L'on devra r\'eviser compl\`etement les manuels d'histoire, et, au moins autant, les manuels d'enseignement de la langue allemande, o\`u il y aura lieu de faire ressortir compr\'ehensivement les valeurs permanentes et les valeurs actuelles de l'\^ame et de l'esprit allemands.

Toutes les biblioth\`eques scolaires, universitaires, d'\'ecoles normales (tant qu'elles subsistent!) et municipales devront \^etre abondamment pourvues d'ouv\-rages classiques et r\'ecents r\'epondant \`a nos desiderata.

Il faudra imposer d\'esormais \`a tous les programmes d'examens et de concours de toutes cat\'egories la connaissance de la langue et de la civilisation allemandes.

N'oublions pas, par ailleurs, qu'il serait int\'eressant d'organiser des cours de langue et de civilisation \`a l'usage des Allemands.

Le haut commissariat, enfin, aurait \`a se pr\'eoccuper de se tenir en constant liaison avec les services allemands:  condition \'evidemment essentielle de la collaboration.  Ainsi s'effectuerait sans heurts la mise au point de questions qui restent pendantes et qui, de ce fait, demeurent irritantes.  Une volont\'e -- et une bonne volont\'e -- de collaboration aplanirait combien de difficult\'es dont on n'a pas vu jusqu'ici la solution bien nettement!

\begin{center}
* * *
\end{center}

Ainsi l'on pourrait escompter qu'un esprit nouveau s'\'etablirait progressivement en France.  Le haut commisariat, par son existence m\^eme, constituerait un encouragement et un reconfort pour tous ceux qui, dans l'Universit\'e et dans les milieux divers, s'attachent avec courage et foi \`a suivre le mot d'order du Mar\'echal, et se fondent sur les paroles de l'amiral, qui, ayant \`a opter pour notre pays entre la vie et la mort, a choisi la vie.

Le haut commisariat serait, enfin, un phare pour notre jeunesse actuelle, d\'esorient\'e et d\'esempar\'ee, qui sent que son avenir et son salut ne se trouvent pas, au del\`a des mers, mais sur le continent, en symbiose avec notre voisin imm\'ediat de l'Est; mais qui, par une d\'elicate pudeur, n'ose pas aller d'elle-m\^eme vers ce voisin syst\'ematiquement d\'ecri\'e jusqu'ici aupr\`es d'elle.

Jeunesse, notre espoir et notre raison d'\^etre, nous le crions, du fond de notre douleur morale et de notre chair meutrie (mais notre vie, \`a tous, doit \^etre l'histoire des victoires que nous remportons sur nous-m\^emes), nous le crions:

\textbf{C'est par la collaboration que tu pourras r\'ealiser tes l\'egitimes aspirations, que tu pourras devenir r\'eellemnet tout ce que tu es virtuellement;}

\textbf{Par la collaboration, dont la base est la compr\'ehension, et dont le but est la paix;}

\textbf{La collaboration, salut de la France et de l'Europe!}

\subsubsection*{English translation}

\noindent
French-German relations:  For the creation of a high commission

\smallskip

\noindent
by R.\ VOIZE, agr\'eg\'e of the university, professor at the lyc\'ee Louis-le-Grand

\smallskip
So Carcopino, minister of national education, is in the process of collapsing.  

He collapses without a jolt, shamefully, in the state of decomposition to which he was brought by his shortcomings and felonies:  a soft thud.  Shortcomings and felonies that denounced him by themselves to the Marshall, to the admiral, to the ambassador representing the French government in the occupied zone, and to the general secretaries of the presidency of the Council.

He was finished off thanks to the determined action of clearsighted resolute university academics, animated by the Marshall's spirit, faithful to the National Revolution, who formed the ``team for university reconstruction''.  They have torn away the mask of clipped ribbons that he was still holding to his face with his cramped fingers.  

Two-faced Janus, he professed his spirit of National Revolution and his good will towards the Collaboration:  and he covered himself on the side of the Popular Front, slyly counting on its return, and on rewards.

He gave them tokens of his fidelity.

He kept Luc in his fat sinecure as general director of technical education and also installed him as Lamirand's controller and crusher of zeal.  (We raised the alarm here:  see \textit{La Gerbe} for 3 July.)

Worse yet:  he was getting ready to reinstall in the important post of rector of the University of Paris the half-breed Roussy, who had been thrown out of the door of the Sorbonne this winter by our courageous national academics.  Roussy, Swiss by origin, who had brought from Switzerland his cheap condensed milk along with Bolshevik tracts.  Roussy the plutocratic demagogue (Nestl\'e provided him a million a year), tinpot scientist, who got himself photographed by \textit{l'Humanit\'e} raising his fist beside Vaillant-Couturier; undesirable straw scientist, who has peopled the university with mediocrities in his own image.  We will send him back to Switzerland, to his milk cows, or keep him in France in a camp that will not be a youth camp.

Roussy had a faithful servant:  Guyot.  Carcopino recycled him.  Guyot, former general secretary of the University of Paris, perfidious aborted thing, evil larva, keeper of all the dirty secrets of rectors who have been disorienting the youth of France morally for twenty years.  This Guyot, chased from the Sorbonne at the same time as Roussy, was fished back out by Carcopino, and he again devotes himself to his subversive activity\dots as head of the office of higher education.

Pulling back from the boundary of notoriety, Carcopino was going to take advantage of the university's torpor in the month of August to put into effect, at every level of education, massive appointments and promotions of professors and administrators driven by a systematic and irreducible hostility towards the government of National Revolution and the Marshall's collaboration.  This is when the team for university reconstruction, the team of the Marshall's men, stood up as a bloc:  the criminal of the University was stopped cold; and they struck him in his rackets.

His harmful activity can be summed up this way:  he deviously supported and encouraged all those in the University who stupidly or criminally stirred up the youth to rebel against German authority and proclaim themselves Gaullist; as for those who tried to orient this youth towards the path to the future indicated by the Marshall and the admiral, Carcopino disavowed, stifled, persecuted, and eliminated them.  In short, he exerted himself to park our youth in a sort of intellectual concentration camp closed to ideals and to reality, and he forbade them any escape, any zeal, any effort of understanding and collaboration.

Leave aside, in the rubble, this Bonaparte in a carton without flesh, this presumptuous clown who plays the emperors and traitors of tragic comedies, this colorful monkey with the grimacing smile.  Let him go find, in the catacombs outside the Temple, the dust of his dead gods.  We will not anger ourselves over this cadaver; let us simply keep it from embalming itself in the sacrophage of the \textit{Ecole normale sup\'erieure}; he is forever unworthy to guide the youth, this promoter of the partisanly anti-religious, demagogic, bolshevising and profiteering petty Mafia, whose arranger was the dervish dancer Durkheim and whose amply certified barker was the intellectual dolt Bougl\'e.

\begin{center}
* * *
\end{center}

Let us speak constructively:  so that your reign will come, Collaboration, salvation of France and of Europe.

Collaboration has understanding as its foundation and peace as its goal.

For a year, since the exodus, the Germans have patiently, tenaciously, and quietly deployed their efforts in every domain:  touching assistance the miserable hordes of the French populations (look along the paths of the retreat); work found in France and in Germany for French men and women who need to find sustenance for themselves and their families; generosity concerning the faults of the trembling vanquished\dots

In this article, let us limit our examination to the efforts at collaboration at the intellectual level. It is at this level that collaboration may present the least difficulty:  because, as Leibniz said, ``minds are the bodies that impede themselves the least.''  And then, we believe in the predominant value of the elite of the elites.  The elites are made to govern, and the masses to be governed.  And the elites, for whom reflex is overcome by reflection, are more particularly apt at understanding that to live is to overcome the past, a past mortgaged by all our shortcomings, our effeteness, our abdications.

The Germans are working at the rapprochement of the elites.  The German Institute, masterfully directed by Dr.\ Epting and his assistant Dr.\ Bremer, puts French of good will in contact with the German language and civilization:  an enterprise crowned with success, with incalculable repercussions.  In addition, eminent Germans, scientists and engineers, come to make contact with the French milieus corresponding to their specialties:  most recently the Reich's minister of agriculture.  

What are the French doing?  The French are acting.  The group Collaboration, with its various sections so well supplied and so active; the Institute of Germanic Studies at the Sorbonne; the different groupings and newspapers who invoke the Marshall and are driven by the spirit of the National Revolution, the French-German journals.

So the Germans act, and some French act. 

\begin{center}
* * *
\end{center}

But what is missing is a \textit{coordinating organism} for all these impressive efforts.

Only an official organism would be able to achieve this indispensable and urgent coordination, because it would appear perfectly qualified in everyone's eyes.

This coordinating organism could take the form of a \textit{high commission for French-German intellectual relations}.

It should be attached directly to the presidency of the Council, thus functioning as a super-ministerial organism, so that it can impose its authority on the various ministries, especially the Ministry of National Education, and prescribe measures that would apply equally to the unoccupied zone and the occupied zone.

If it had existed, it would have ordered the Education Ministry to take advantage of the awarding of prizes to provide many young people, ranked among the best, with the Marshall's, the admiral's, and Mr.\ Lamirand's speeches, where the reasons for the Collaboration are explained; as well as French books, or books translated into French, that could launch them on an understanding of Germany and a collaboration of hearts and minds.

If he had existed, it might have tried to use the summer vacation to try to create contact between the French and German youth.

In any case, we can set out to to organize more and more professional contacts.

We should completely revise our history textbooks, and at least as much, the textbooks for teaching the German language, where there will be space for bringing out in a comprehensive way the permanent and present values of the German soul and spirit.

All the libraries in schools, universities, normal schools (such as they are!) and municipalities should be generously provided with classic and recent works satisfying our desiderata. 

From now on, we should require knowledge of the German language and civilization in the programs of examinations and competitions of all categories.

What's more, do not forget that it would be interesting to organize courses in language and civilization for Germans to use.

Finally, the high commission should make sure that it stays in constant liaison with the German services:  an obviously essential condition for collaboration.  In this way the resolution of questions that are still pending, and therefore still irritating, can be carried out without clashes.  A will -- and a good will -- for collaboration would smooth out so many of the difficulties that we have not until now seen clearly how to solve.

\begin{center}
* * *
\end{center}

Thus we could count on a new spirit gradually establishing itself in France.  The high commission by its very existence, would constitute encouragement and comfort for everyone, in the University and various other milieus, who devote themselves with courage and faith to following the Marshall's marching orders and who anchor themselves in the words of the admiral, who, having to choose between life and death for our country, chose life.

Finally, the high commission would be a lighthouse for our present youth, disoriented and despairing, who feel that their future and their salvation is to be found not beyond the seas but on the continent, in symbiosis with our immediate neighbor to the East, but who out of a delicate sense of decency do dare not move by themselves towards a neighbor who, until now, is systematically denigrated to them.

Youth, our hope and our reason for being, we cry to you, from the bottom of our moral sorrow and our bruised flesh (but our live, for everyone, should be the history of our victories over ourselves), we cry to you:

\textbf{It is by collaboration that you can achieve your legitimate aspirations, that you can become in reality all that you are virtually;}

\textbf{By collaboration, whose foundation is understanding, and whose goal is peace;}

\textbf{Collaboration, salvation of France and of Europe!}

\subsection{Pistre's recommendation of Voize, 11 August 1941}\label{subsec:pistre}

\textit{Edmond Pistre-Caraguel to the Propaganda-Abteilung, 11 August 1941.}

\subsubsection*{French original}

J'ai l'honneur d'appeler votre attention sur l'article publi\'e par "La Gerbe" du 7 Ao\^ut 1941, page 5, sous la signature de R.Voize, Professeur au Lyc\'ee Louis-le-Grand, agr\'eg\'e pour la langue Allemande.

Apr\`es avoir lu cet article, j'ai pris contact avec Monsieur Voize en vue de le faire collaborer par ses conseils comp\'etents \`a ma mission de r\'eformer l'esprit des membres du Corps Enseignant.

Au cours de notre entretien, je me suis charg\'e de vous exprimer le d\'esir de Monsieur Voize d'\^etre re\c cu en audience par le Chef de la Propaganda Abteilung Frankreich susceptible d'appuyer son id\'ee d'un Haut Commissariat aux Rapports Intellectuels Franco Allemands.

L'int\'er\^et que pr\'esenterait, pour le succ\`es de ma mission, l'adoption de cette initiative est que je pourrais ainsi m'appuyer officiellement sur un organisme du Gouvernement Fran\c cais, autre que le Ministre de l'Education Nationale dont l'activit\'e r\'enovatrice a \'et\'e tr\`es discut\'ee jusqu'\`a pr\'esent.

M.\ Voize habite 28 Quai de B\'ethune Paris IV$^{\text{e}}$; il parle parfaitement la langue allemande; vous pourrez ainsi convoquer pour uner conversation directe, si vous le jugez utile.

\subsubsection*{English translation}

It is my honor to call your attention to the article by published by ``La Gerbe'' on 7 August 141, page 5, under the signature of R.\ Voize, professor at the lyc\'ee Louis le Grand, agr\'eg\'e in the German language.  

After reading the article, I got into contact with Mr.\ Voize, hoping to enlist his competent advice for my mission of reforming the mindset of the teaching corps.  

In the course of our conversation, I agreed to contact you about his desire for an audience with the head of the Propaganda Abteilung Frankreich so as to advance his idea of a High Commission for French-German Intellectual Relations.  

The adoption of this initiative would promote the success of my own mission, because it would permit me to rely officially on an organism of the French government other than the Ministry of National Education, whose renovating activity has been very questionable so far.

Mr.\ Voize lives at 28 Quai de B\'ethune Paris IV$^{\text{e}}$; he speaks perfect German; so you can call him in for a direct conversation if you judge this to be useful.

\subsection{Denunciation of Carcopino's appointments, 28 August 1941}\label{subsec:godefroi}

\textit{Excerpts from an article on page 4 of the 28 August 1941 issue of \textit{Le cri du peuple}.}

\subsubsection*{French original}

\noindent
Dictature

\smallskip
\noindent
Par Pierre Godefroi

\smallskip
Le mois de juin 1940 a permis aux Fran\c{c}ais de mesurer l'ampleur de leur d\'eb\^acle militaire.  Il leur apparut que le fond de l'abime \'etait atteint, que l'on ne pourrait tomber plus bas.  L'ann\'ee qui vient de s'\'ecruler leur offre le lamentable spectacle de la d\'eb\^acle morale dans l'Enseignement.  Par une contradiction v\'eritablement surprenante alors qu'en principe la France vit sous un r\'egime autoritaire qui professe l'anti-ma\c{c}onnisme, en r\'ealit\'e, la Ma\c{c}onnerie fait peser plus lourdement son joug sur les patriotes, tandis qu'elle \'el\`eve ses fid\`eles serviteurs.  Autre ab\^ime o\`u chaque jour qui passe nous p\'ecipite plus profond\'ement!  Comme on comprend bien l'ironie de ce proviseur d'un plus grands lyc\'ees de Paris, qui, il y a quelques jours \`a peine, fasiait allusion en riant aux nouvelles attestations demand\'ees aux professeurs concernant la Franc-Ma\c{c}onnerie!  Il savait pertinemment la vanit\'e de ces formules et de ces serments sur l'honneur, puisque les postes essentiels, les leviers de commande sont toujours occup\'es soit part des ma\c{c}ons, soit par des ma\c{c}onnisants, ce qui est plus grave.  Qu'importent donc des signatures plac\'ees au bas de formulaires administratifs, qui ne sont m\^eme pas \'ecrits en un fran\c{c}ais correct?  La secte reconnait toujours les siens et un parjure de plus les rend plus chers \`a son c\oe ur, et plus dignes de sa confiance.

Une longue suite de scandales jalonne le calvaire que la gestion des secr\'etaires d'Etat successifs \`a l'Education Nationale a fait gravir aux professeurs et instituteurs patriotes, fid\`eles au Mar\'echal et d\'evou\'es corps et \^ame \`a la R\'evolution Nationale:  nominations de MM.\ Charmoillaux, Maurain, Piobetta, Luc, Masbou, Chattelun, Santelli que ce journal a mises en lumi\`ere et sur lesquelles il ne se lassera pas de revenir.  Notons en contre-partie la sanction inqualifiable qui a frapp\'e Serge Jeanneret, coupable du crime de l\`ese-majest\'e vis-\`a-vis de la Franc-Ma\c{c}onnerie.  Jusqu'\`a pr\'esent, il est le seul patriote qui ait \'et\'e frapp\'e, mais tous ses amis furent l'objet de menaces plus ou moins voil\'ees.  Le jour n'est pas loin sans doute o\`u leur tour viendra.

Notre excellent confr\`ere ``Je suis Partout''  annon\c cait dans son num\'ero du lundi 4 ao\^ut que M.\ Carcopino mettrait les vacances \`a profit pour se livrer aux nominations les plus conformes \`a l'esprit du Front Populaire. Une premi\`ere nomination publi\'ee par la radiodiffusion nationale en apporte l'\'eclatante confirmation. Mme Edm\'ee Hatinguais, directrice du Lyc\'ee Racine, est nomm\'ee directrice de l'Ecole Normale Sup\'erieure de S\`evres en remplacement de Mme Cotton-Feytis.  Nul n'ignore dans l'enseignement f\'eminin combien cette personne ondoyante et diverse sait unir la souplesse \`a la fausset\'e. Amie de Juifs notoires, patronn\'ee par des Juifs, ma\c connisante de bonne classe, elle a profess\'e au cours de la derni\`ere ann\'ee scolaire par ses paroles, par ses entretiens avec les \'el\`eves, par ses attitudes, des sentiments gaullistes, ce qui ne l'a pas emp\^ech\'ee de faire en son temps un petit voyage \`a Vichy pour se couvrir, le cas \'ech\'eant, d'un manteau protecteur.  Cette politique vient pleinement de r\'eussir.  Mais patience le dossier de cette dame s'amplifie.  Il d\'ebordera un jour.  Nous l'ouvrions bient\^ot, devant nos lecteurs.

\subsubsection*{English translation}

\noindent
Dictatorship

\smallskip
\noindent
By Pierre Godefroi

\smallskip

The month of June 1940 allowed the French to measure that magnitude of their military debacle.  It looked to them like the bottom of the abyss had been reached, that we could not fall any lower.  The year just passed offers them the lamentable spectacle of the moral debacle in Education.  By a truly surprising contradiction, whereas in principle France lives under an authoritarian regime that professes anti-Freemasonism, in reality Freemasonry is making its yoke weigh more and more heavily on patriots, while it raises up its faithful servants.  Another abyss where every day that passes pushes us deeper.  How well we understand the irony of the deputy head of one of the largest lyc\'ees in Paris who, only a few days ago, alluded with laughter to the new demands made to professors concerning Freemasonry!  He understand precisely the emptiness of these formalities and oaths on one's honor, because the essential posts, the levers of command, are still occupied either by Masons or by the masonizing, which is worse.  So what does it matter that there are signatures at the bottom of administrative forms, which are not even written in correct French? The sect always recognizes its own, and yet one more perjury makes them that much more dear to its heart, more worthy of its confidence.

The Calvary imposed by the management of national education by successive ministers on patriotic professors and school teachers faithful to the Marshall and devoted body and soul to the National Revolution has been marked by a long sequence of scandals: appointments of Charmoillaux, Maurain, Piobetta, Luc, Masbou, Chattelun, Santelli, which this newspaper has brought to light and about which it will not get tired of bringing up.  Notice on the other side the unspeakable punishment that has served on Serge Jeanneret, guilty of the crime of \textit{l\`ese-majest\'e} towards Freemasonry.  He is the only patriot to be struck so far, but all of his friends have been the objects of more or less veiled threats.  The day is not far away when their turn will come.

Our excellent sister publication \textit{Je suis Partout} announced in its issue of 4 August that Mr.\ Carcopino would take advantage of the summer vacation to devote himself to making appointments totally in the spirit of the Popular Front.  The first appointment announced by the national radio confirms this spectacularly.  Mrs.\ Edm\'ee Hatinguais, director of the Lyc\'ee Racine, has been named to replace Mrs.\ Cotton-Feytis as director of the \textit{Ecole Normale Sup\'erieure} de S\`evres.%
\footnote{The \textit{Ecole Normale Sup\'erieure de S\`evres} was the female counterpart of Paris's \textit{Ecole Normale Sup\'erieure}.  The physicist Eug\'enie Cotton (1881-1967), director since 1936, had been forced to retire in June 1941.  Her husband, the physicist Aim\'e Cotton, was later arrested along with Borel.}
Everyone in women's education knows how this wriggling and protean personality succeeds in combining suppleness with dishonesty.  Friend of notorious Jews, patronized by the Jews, high-class Masonizer, she has exhibited her Gaullist sentiments this past academic year by her words and her conversations with students, even though this did not keep her from finding time for a little trip to Vichy to protect herself.

\subsection{Voize to Dahnke, 6 September 1941}\label{subsec:voizedahnke}

\textit{This is the body of a handwritten letter addressed to ``Monsieur'', and signed by Voize.}

\subsubsection*{French original}

A propos de M.\ Chattelun, je vous signale un article que l'on vient de me communiquer:  article parfaitement exact, et particuli\`erement courageux puisqu'il \'emane d'un subordonnn\'e direct de Chatelun, instituteur:  Paris, Jeanneret bien connu depuis longtemps pour ses tendances nationales.  Cet article a paru dans le ``Cri du Peuple'' de mercredi 17 Ao\^ut, en premi\`ere page:  ``Propagande [?] d'un \'energum\`ene.''

D'ici quelques jours, je vous enverrai mes consid\'erations sur le Haut Commissariat. Finalement, je ne briguerai pas ce poste, pr\'ef\'erant me cantonner dans une activit\'e purement priv\'ee : et je vais me consacrer \`a cr\'eer ``l'Institut Langues et Culture'' dont je vous ai parl\'e : je serai heureux, \`a ce sujet que vous me fassiez savoir quand je pourrai vous revoir.

\subsubsection*{English translation}

Concerning Mr.\ Chatelun, I am sending you an article that was just sent to me, an article that is perfectly accurate and particularly courageous because it comes from a direct subordinate of Chatelun's, the Paris schoolteacher Jeanneret, long well known for his national tendencies.  This article appeared in the \textit{Cri du Peuple} for Wednesday, 17 August, on the first page:  ``A Firebrand's Propaganda [?]''.

In a few days I will send you my ideas about the High Commission.  I have decided not to be a candidate for this post, preferring to limit myself to purely private activity, and I want to concentrate on creating the `Institute for Languages and Cultures' that I talked with you about.  In this connection I would be pleased if you could let me know when I can see you again.

\subsection{French police reports, October--November 1941}\label{subsec:frenchpolice}

\textit{These two reports, from the archives of Paris Police (DOSSIER BA 1798), show the Paris police first confirming the arrest of the four academicians and then reporting on the reaction in university circles.}

\subsubsection*{French originals}

\noindent
16 octobre 1941

\smallskip

A la suite d'une information signalant que cinq membres de l'Acad\'emie des sciences, MM.\ Langevin, Lapicque, Mauguin, Borel et Cotton avaient \'et\'e arr\^et\'es, il a \'et\'e proc\'ed\'e \`a des v\'erifications qui ont donn\'e les r\'esultats suivants

M.\ Paul Langevin, ancien professeur \`a la Facult\'e des Sciences, \textit{Coll\`ege de France}, en retraite depuis plus d'un an, n'habite plus 10, rue Vauquelin. Il est actuellement \`a Troyes et on n'a pas connaissance de son arrestation.

M.\ Lapicque, demeurant 17, rue Soufflot, est retrait\'e depuis plusieurs ann\'ees. Depuis juillet dernier, il se trouvait \`a Ploubalzanec (Cotes du Nord). Vendredi matin, 10 octobre, des membres de la police allemande se sont pr\'esent\'es \`a son domicile \`a Paris et en son absence ont proc\'ed\'e \`a une visite domiciliaire \`a la suite de laquelle les scell\'es ont \'et\'e appos\'es aux portes de son appartement. Une op\'eration a du \^etre effectu\'ee le m\^eme jour en Bretagne et M.\ Lapicque a \'et\'e incarc\'er\'e \`a Saint-Brieuc. Une premi\`ere visite avait d\'ej\`a eu lieu l'ann\'ee derni\`ere, par les autorit\'es d'occupation, au domicile de M.\ Lapicque.

M.\ Mauguin, professeur en activit\'e \`a la Sorbonne, 1 rue victor Cousin, a \'et\'e arr\^et\'e le 10 octobre courant dans cet \'etablissement. Emmen\'e par les autorit\'es occupantes, il a du \^etre conduit \`a son domicile \`a Thiais pour visite domiciliaire.

M.\ Borel Emile, n\'e le 7 janvier 1877 [sic] \`a Ste Affrique (Aveyron), membre de l'Institut, ancien professeur \`a la Sorbonne et ancien d\'eput\'e de l'Aveyron, est domicili\'e depuis six mois avec sa femme, n\'ee Marguerite Appell, 4 rue Froidevaux.  Il a \'et\'e appr\'ehend\'e le vendredi 10 octobre courant, vers 17h30, par les autorit\'es allemandes, \`a son domicile, \`a la suite d'une perquisition. On ignore les motifs de 
son arrestation.

\smallskip
\noindent
*****************
\smallskip

\noindent
7 novembre 1941

\smallskip

Prenant pr\'etexte des op\'erations qui auraient \'et\'e effectu\'ees r\'ecemment par la Police allemande au domicile de plusieurs savants fran\c cais, notamment MM.\ Langevin, Emile Borel, Louis Lapicq, Mauguin, Cotton, etc..., les dirigeants communistes s'efforcent actuellement de d\'evelopper leur propagande dans les milieux universitaires.

Le but principal de cette agitation est d'exploiter, pour les besoins de leur propagande, la sympathie de tous les milieux universitaires pour ces savants, en leur faisant conna\^\i tre ce qu'ils appellent: ``les bienfaits de la collaboration'' et en d\'eveloppant parmi tous les intellectuels un courant d'hostilit\'e aux Autorit\'es allemandes et au Gouvernement fran\c cais.

Dans ce but, les militants communistes ont re\c cu pour instructions de cr\'eer parmi tous les universitaires une inqui\'etude g\'en\'erale en soulignant le sort qui leur est r\'eserv\'e, s'ils ne s'opposent pas d\`es maintenant aux brimades allemandes et aux sanctions appliqu\'ees par le Gouvernement sous les pr\'etextes les plus futiles. 

Pour donner \`a cette agitation son maximum d'effets, il est recommand\'e aux propagandistes de s'adresser de pr\'ef\'erence aux intellectuels connus pour leurs sympathies au Front Populaire, \`a la Franc-Ma\c connerie, ou pour leur anglophobie [sic] ou leur germanophobie, puis gagner peu \`a peu ceux appartenant \`a des groupements de droite.

Cette propagande devra \^etre men\'ee de fa\c con \`a d\'emontrer aux intellectuels qu'en s'attaquant aux universitaires et aux savants l'Allemagne n'a d'autre but que de briser en France tout ce qui est science, pens\'ee honn\^ete et ind\'ependante, tout ce qui est national et qui emp\^eche de transiger avec l'ennemi, de s'entendre avec lui contre la Patrie, pour encourager tout ce qui est m\'epris de l'intelligence et exaltation de la force brutale.

\subsubsection*{English translations}

\noindent
16 October 1941

\smallskip

Following a report that five members of the \textit{Acad\'emie des Sciences}, Langevin, Lapicque, Mauguin, Borel, and Cotton, had been arrested, an investigation was conducted, producing the following results.

Mr.\ Paul Langevin, formerly professor at the Faculty of Sciences, College of France, in retirement for more than a year, no longer lives at 10, rue Vauquelin.  He is now at Troyes, and people do not know about his arrest.

Mr.\ Lapicque, living at 17, rue Soufflot, has been retired for many years.  Since last July, he has been at Ploubalzanec (Cotes du Nord).  On Friday morning, 10 October, members of the German police presented themselves at his home in Paris and searched it his absence, then placing seals on the doors of his apartment.  An operation was carried out on the same day in Brittany, and Mr.\ Lapicque was imprisoned at Saint-Brieuc.  The Occupation authorities had already visited his home last year.

Mr.\ Mauguin, currently a professor at the Sorbonne, 1 rue Victor Cousin, was arrested this 10 October in that establishment.  He must have been conducted by the occupying authorities to his home at Thiais for a search.

Mr.\ Borel Emile, born 7 January 1877 [sic] at Ste Affrique (Aveyron), member of the Institut, formerly professor at the Sorbonne, and formerly deputy for Aveyron, has been living for six months with his wife, born Marguerite Appell, at 4 rue Froidevaux.  He was apprehended at his home this last Friday, 10 October, around 17:30, by the German authorities, following a search.  The motives for his arrest are not known.

\smallskip
\noindent
*****************
\smallskip

\noindent
7 November 1941

\smallskip

Taking as a pretext the recent operations that may have been carried out by the German police at the homes of several French scientists, notably Langevin, Emile Borel, Louis Lapicque, Mauguin, Cotton, etc.\dots, communist leaders are now trying to extend their propaganda into university circles.

The main goal of this agitation is to exploit, for the purposes of their propaganda, the sympathy for these scientists in all university circles, making them aware of what they call ``the benefits of collaboration'' and developing among all intellectuals a current of hostility towards the German authorities and the French government.

To this end, the communist militants have been instructed to create general anxiety among all university academics by underlining the fate reserved for them if they do not oppose forthwith German intimidation and the punishments applied by the government under the most futile pretexts.

In order to give this agitation the greatest effect, the propagandists are advised to address first of all intellectuals known for their sympathies for the Popular Front and Freemasonry or for their anglophobia [sic] or their Germanophobia, and then to win over bit by bit those belonging to groups on the right.

This propaganda is to be conducted so as to demonstrate to intellectuals that the Germans have no other object in their attacks on the university academics and scientists than to destroy everything in France that is science, independent and honest thought, and everything that is national and impedes compromising with the enemy and coming to terms with them against the country, and to encourage everything that is contempt for intelligence and exaltation of brutal force.

\subsection{Gidel and von St\"ulpnagel, 15 December 1941}\label{subsec:gidel}

\textit{Excerpt from a report on the conversation on 15 December 1941 at the MBF headquarters at Hotel Majestic between Gilbert Gidel, rector of the Academy of Paris, and Otto von St\"upnagel, military commander in France.  Now in CARAN AJ/16/7117, in the file on the Basdevant affair.} 

\subsubsection*{French original}

Le g\'en\'eral d\'eclare au d\'ebut de l'entretien qu'il a donn\'e jadis son assentiment \`a la nomination de M.\ Gidel comme Recteur de l'Acad\'emie de Paris, en exprimant l'espoir que le nouveau Recteur saurait faire preuve, dans ses nouvelles fonctions, de l'\'energie n\'ecessaire pour \'eviter tous incidents qui pourraient conduire \`a des mesures de r\'epression \'egalement f\^acheuses pour la Sorbonne et les autorit\'es d'occupation.

Il sait que la jeunesse universitaire, par une fausse conception du patriotisme, peut se laisser entra\^\i ner \`a des actes de turbulence. Mais il faut que cette jeunesse sache \`a son tour que le g\'en\'eral pourrait \^etre amen\'e \'eventuellement \`a prendre de s\'ev\`eres mesures \`a l'\'egard de la Sorbonne. Il rappelle qu'il a \'et\'e dans l'obligation de fermer la Sorbonne apr\`es les incidents de novembre 1940 et assure qu'il regretterait profond\'ement d'avoir \`a prendre \`a nouveau pareille mesure. Car tous ses efforts ont tendu \`a rendre possible une atmosph\`ere de r\'eelle collaboration dans l'espoir d'un avenir heureux pour la France et pour l'Europe.

Le g\'en\'eral fait alors allusion au cas de M.\ Villey et de M.\ Basdevant. Je ne peux m'emp\^echer de croire, dit le G\'en\'eral, que parmi les professeurs de la Sorbonne, il y ait actuellement des \'el\'ements qui cherchent \`a fomenter le d\'esordre.

Le Recteur expose que ses efforts se sont appliqu\'es \`a faire r\'egner le calme dans une atmosph\`ere de travail. 'Et je crois y avoir r\'eussi, ajoute-t-il, puisque la journ\'ee du 11 novembre s'est d\'eroul\'ee gr\^ace aux mesures prises, sans le moindre incident.'

\subsubsection*{English translation}

The general declared at the beginning of the conversation that when he had agreed to Mr.\ Gidel's appointment as Rector of the Academy of Paris, he had expressed the hope that the new rector would manage to show in his new role the energy needed to avoid any incidents that could lead to repressive measures that would be equally unpleasant for the Sorbonne and the Occupation authorities.

He knows that the university youth, because of a false conception of patriotism, could allow themselves to be drug into disruptive actions.  But this youth should also know that the general could in turn be led to take severe measures with respect to the Sorbonne.  He recalled that he had been obliged to close the Sorbonne after the incidents of November 1940, and he vouched that he would profoundly regret having to take such a measure again.  Because all his efforts had been to make possible an atmosphere of real collaboration in the hope of a happy future for France and Europe.

The general then alluded to the Mr.\ Villey's and Mr.\ Basdevant's cases.  I cannot keep from thinking, the general said, that among the professors of the Sorbonne, there are now elements trying to foment disorder.

The rector explained that his efforts had been to establish calm and an atmosphere of work.  And I believe I have succeeded, he added, because the day 11 November unfolded, thanks to the measures taken, without the least incident.

\subsection{Letters from Borel to Lacroix, 1941--1944}\label{subsec:lacroix}

\noindent
\textit{These letters, from Emile Borel to Alfred Lacroix, permanent secretary of the \textit{Acad\'emie des Sciences}, are preserved in the academy's archives.}

\subsubsection*{French originals}

\noindent
Paris, 19 novembre 1941

\noindent
Mon cher confr\`ere,

Merci de votre sympathie, qui m'est pr\'eciseuse. Je me r\'ejouis de vous revoir lundi.

\noindent
Votre bien d\'evou\'e,

\noindent
Emile Borel

\smallskip
\noindent
*****************
\smallskip

\noindent
St Affrique, 16 octobre [1943?]

\noindent
Monsieur et cher confr\`ere,

Il ne m'est pas possible, \`a mon grand regret, de rentrer \`a Paris pour les \'elections; veuillez m'en
excuser aupr\`es de nos confr\`eres et en particulier aupr\`es du Bureau.
Vous seriez bien aimable si de nouvelles \'elections ou d'autres circonstances, comme des r\'eunions
de plusieurs commissions de prix, rendaient d\'esirable ma participation aux travaux de l'Acad\'emie,
de vouloir bien m'en informer assez longtemps d'avance pour qu'il me soit possible d'obtenir en
temps utile le laissez-passer me permettant de revenir \`a Paris.

Veuillez agr\'eer, Monsieur et cher confr\`ere, l'expression de mes sentiments d\'ef\'erents et
cordialement d\'evou\'es,

\noindent
Emile Borel

\smallskip
\noindent
*****************
\smallskip

\noindent
St Affrique, 28 d\'ecembre 1943

\noindent
Messieurs les Secr\'etaires Perp\'etuels,

J'ai bien re\c{c}u votre lettre relative aux places de correspondants de la section de g\'eom\'etrie. J'avais
d\'ej\`a correspondu \`a ce sujet avec plusieurs membres de la section et je pense que celle-ci serait
dispos\'ee \`a proc\'eder \`a une \'election, en r\'eservant sa d\'ecision au sujet de la seconde place.

J'\'ecris \`a M.\ Cartan de s'entendre avec vous pour la date de convocation de la section afin de faire
les pr\'esentations pour cette \'election.

Veuillez agr\'eer l'expression de ma haute consid\'eration et de mes sentiments d\'evou\'es.

\noindent
Emile Borel

\smallskip
\noindent
*****************
\smallskip

\noindent
Paris, 23 septembre 1944

\noindent
Monsieur le Secr\'etaire perp\'etuel et cher confr\`ere,

La bienveillance que vous m'avez t\'emoign\'ee me fait un devoir de vous tenir au courant de ce qui
suit.

D\`es la lib\'eration, j'ai cherch\'e \`a \^etre exactement renseign\'e sur les r\'eactions de l'Acad\'emie au
moment de mon incarc\'eration \`a Fresnes. M.\ Vincent a bien voulu me raconter en d\'etail ses
d\'emarches, en tant que Pr\'esident de l'Acad\'emie, aupr\`es de M.\  de Brinon. Celui-ci le persuada
que toute d\'emarche en faveur des acad\'emiciens incarc\'er\'es ne pourrait avoir que les plus graves
dangers pour eux-m\^emes et pour l'Acad\'emie. Par suite, lorsque, apr\`es le d\'ec\`es de M.\ Picard, un
courant se manifesta dans l'Acad\'emie en faveur de ma candidature, M.\ Vincent ainsi que
M.\ Esclangon, vice-pr\'esident, consid\'er\`erent comme de leur devoir de s'opposer \`a cette
candidature. M.\ Vincent a r\'ep\'et\'e cette conversation \`a M.\ Roussy et M.\ Esclangon me l'a confirm\'ee,
en ajoutant qu'\`a son avis, j'aurais \'et\'e s\^urement \'elu si je n'avais pas \'et\'e arr\^et\'e par les Allemands.
De l'avis de M.\ Vincent et de M.\ Esclangon, c'est donc cette arrestation seule qui a emp\^ech\'e mon
\'election puisque tout candidat qui m'\'etait oppos\'e devait avoir, outre des voix personnelles, toutes
celles de ceux qui partageaient l'opinion des Pr\'esidents sur les dangers de mon \'election. Malgr\'e
cela, au d\'ebut de janvier, j'avais de nombreuses promesses et mon succ\`es paraissait assur\'e. Le
coup de gr\^ace me fut alors donn\'e par M.\ Carcopino qui, comme vous le savez mieux que
personne, exigea le retrait de ma candidature, pour les m\^emes raisons sugg\'er\'ees \`a M.\ Vincent par
M.\  de Brinon.

Je suis donc en droit de penser que, sans le vouloir, l'Acad\'emie m'a inflig\'e une peine
suppl\'ementaire, s'ajoutant aux cinq semaines d'incarc\'eration. Il me semble que j'ai le droit de
demander qu'une r\'eparation me soit accord\'ee. La r\'eparation la plus compl\`ete serait la d\'emission
de M.\ Louis de Broglie et mon \'election. Certains amis de M.\ Louis de Broglie et son fr\`ere sont
d'accord pour la sugg\'erer, d\`es le retour \`a Paris du duc de Broglie. Ce serait le duc de Gramont qui
parlerait au duc de Broglie et M.\ Julia qui parlerait \`a M.\ Louis de Broglie. Il est vraisemblable que
vous serez alors consult\'e et j'ai toute confiance en votre esprit de justice.

Veuillez agr\'eer, Monsieur le Secr\'etaire Perp\'etuel et cher confr\`ere, l'expression de mes sentiments
cordialement d\'evou\'es.

\noindent
Emile Borel

\smallskip
\noindent
*****************
\smallskip

\noindent
Paris, le 6 octobre 1944

\noindent
Monsieur le secr\'etaire perp\'etuel et cher confr\`ere,

J'ai eu lundi avec Maurice de Broglie une conversation tr\`es cordiale. Il me para\^it \`a peu pr\`es
certain que Louis de Broglie donnera sa d\'emission de secr\'etaire perp\'etuel dans quelques
semaines; il reprendra donc dans la section de m\'ecanique la place vacante par le d\'ec\`es de
Jouguet.

Je compte causer lundi avec Louis de Broglie, mais j'ai tenu \`a vous mettre au courant sans tarder.

\noindent
Votre bien d\'evou\'e,

\noindent
Emile Borel

\subsubsection*{English translations}

\noindent
Paris, 19 November 1941

\noindent
My dear colleague,

Thank you for your sympathy, which is important to me.  I look forward to seeing you on Monday.

\noindent
Your very devoted,

\noindent
Emile Borel

\smallskip
\noindent
*****************
\smallskip

\noindent
St affrique, 16 October [1943?]%
\footnote{The year cannot be read clearly, but the stamp on the postcard features P\'etain.}

\noindent
Sir and dear colleague,

To my great regret, it is not possible for me to return to Paris for the elections; please excuse me to our colleagues and especially to the office.  If new elections or other circumstances, such as the convening of many prize committees, make my participation in the Academy's work desirable, it would be very kind of you to let me know soon enough in advance that I can obtain in good time the pass allowing me to return to Paris.

Please accept, sir and dear colleague, the expression of deferential and cordially devoted sentiments,

\noindent
Emile Borel

\smallskip
\noindent
*****************
\smallskip

\noindent
St Affrique, 28 December 1943

\noindent
Dear permanent secretaries,

I did indeed receive your letter concerning the posts of corresponding member in the geometry section.  I had already corresponded with several members of the section on this topic, and I think that the section would be disposed to proceed to an election while reserving a decision with respect to the second post.

I am writing to Mr.\ Cartan about agreeing on a date for convening the section for presentations for the election.

Please accept the expression of my high consideration and my devoted sentiments.

\noindent
Emile Borel

\smallskip
\noindent
*****************
\smallskip

\noindent
Paris, 23 September 1944

\noindent
Mr.\ permanent secretary and dear colleague,

The kindness you have shown me makes it my duty to keep you informed of what follows.

Since the Liberation, I have tried to learn exactly how the Academy reacted at the time of my imprisonment at Fresnes.  Mr.\ Vincent was glad to tell me in detail about his conversations, as president of the Academy, with Mr.\ de Brinon,%
\footnote{Fernand de Brinon (1885--1947) was the Vichy government's official representative in Paris.  After the war, he was tried, condemned to death, and executed for collaboration.}
who persuaded him that any action in favor of the imprisoned members of the academy could only create the greatest of dangers for them and for the academy.  Later, after the death of Mr.\ Picard, when sentiment emerged in the Academy in favor of my candidacy, Mr.\ Vincent and also the vice president, Mr.\ Esclangon, considered it their duty to oppose that candidacy.  Mr.\ Vincent repeated that conversation to Mr.\ Roussy, and Mr.\ Esclangon confirmed this to me, adding that in his opinion, I would surely have been elected if I had not been arrested by the Germans.  In the opinion of Mr.\ Vincent and Mr.\ Esclangon, my election was therefore blocked solely by that arrest, because any opposing candidate, aside from the votes that he attracted personally, would have all those of members who shared the presidents' opinion about the dangers of my election.  In spite of that, at the beginning of January I had many promises, and my success appeared assured.  The coup de gr\^ace came from Mr.\ Carcopino, who, as you know better than anyone, demanded the withdrawal of my candidacy, for the same reasons suggested by Mr.\ Vincent and Mr.\ de Brinon.

So I have the right to feel that the Academy, without wanting to do so, added a supplemental punishment to my five weeks of imprisonment.  It seems to me that I have the right to demand reparation.  The most complete reparation would be Mr.\ Louis de Broglie's resignation and my election.  Some of the friends of Mr.\ de Broglie and his brother have agreed to suggest it, as soon as the Duke de Broglie%
\footnote{At this point Louis's brother Maurice (1875--1960), also a physicist and member of the \textit{Acad\'emie des Sciences}, held the family title of Duke.  Louis became Duke upon Maurice's death.}
returns to Paris.  The Duke de Gramont should talk with the Duke de Broglie, and Mr.\ Julia should talk with Mr.\ Louis de Broglie.  It is likely that you will then be consulted, and I have complete confidence in your sense of justice.

Please accept, Mr.\ permanent secretary and dear colleague, the expression of my cordially devoted sentiments.

\noindent
Emile Borel

\smallskip
\noindent
*****************
\smallskip

\noindent
Paris, 6 October 1944

\noindent
Mr.\ permanent secretary and dear colleague, 

I had a very cordial conversation with Maurice de Broglie.  It seems nearly certain to me that Louis de Broglie will resign as permanent secretary in a few weeks; it will then take up again a place in the section on mechanics left vacant by the death of Jouguet.

I expect to chat with Louis de Broglie on Monday, but I wanted to keep you informed without delay.

\noindent
Your very devoted

\noindent
Emile Borel

\subsection{Eug\'enie Cotton's recollections, 1967}\label{subsec:cotton}

\textit{Excerpt from Eug\'enie Cotton's biography of her husband, Aim\'e Cotton \cite{cotton:1967}.}

\subsubsection*{French original}

J'ai toujours pens\'e que l'arrestation des quatre acad\'emiciens progressistes avait \'et\'e li\'ee aux manifestations d'\'etudiants qui avaient eu lieu \`a l'Arc de Triomphe de l'Etoile, en 1940, \`a l'occasion de l'armistice de 1918. Paul Langevin avait \'et\'e rendu responsable de ces manifestations et il avait \'et\'e arr\^et\'e et maintenu depuis en r\'esidence surveill\'ee. Pour \'eviter le renouvellement de semblables manifestations, les Allemands avaient pris les devants et arr\^et\'e les savants qui leur paraissaient susceptibles de pouvoir entra\^iner les \'etudiants \`a manifester le 11 novembre 1941. Arr\^et\'es les 10 et 11 octobre 1941 \`a la date o\`u la pr\'eparation de la manifestation aurait d\^u commencer, Borel, Cotton, Lapicque et Mauguin furent lib\'er\'es le 13 novembre sit\^ot pass\'ee la date anniversaire de l'armistice. 

Dans une lettre qu'il m'\'ecrivit apr\`es la mort de mon mari Charles Mauguin a rappel\'e ce qui expliquait \`a ses yeux leur commune arrestation: `` Il faut croire que nous avions quelques traits communs dans nos temp\'eraments et notre philosophie car nous \'etions tr\`es g\'en\'eralement d'accord dans nos mani\`eres de voir.

Naturellement bienveillant, il (Aim\'e Cotton) se donnait sans r\'eticence avec un d\'evouement \`a toute \'epreuve, \`a la cause qu'il estimait juste. Il prenait par contre avec courage et nettet\'e position contre ce qu'il r\'eprouvait.

C'est sans doute ce qui lui valut, en 1941, d'\^etre intern\'e par les Allemands \`a Fresnes, o\`u nous nous trouvions tous les deux, chacun \`a l'insu de l'autre. Je ne suis pas pr\`es d'oublier notre retour en commun dans le m\'etro, \`a notre sortie de prison, avec la t\^ete hirsute de brigands de grand chemin.''

\subsubsection*{English translation}

I have always thought that the arrest of the four progressive academicians was related to the student demonstrations that took place at the Arc de Triomph in 1940 on the occasion of the anniversary of the 1918 Armistice.  Paul Langevin had been held responsible for these demonstrations, and he had been arrested and put under house arrest.  To prevent the recurrence of such demonstrations, the Germans made the first move and arrested the scientists they thought might encourage the students to demonstrate on 11 November 1941.  Arrested on 10 and 11 October, when preparation for the demonstration might have begun, Borel, Cotton, Lapicque and Mauguin were freed on 13 November, soon after the anniversary of the Armistice.

In a letter written to me after my husband's death, Charles Mauguin recalled his opinion about the cause of their common arrest: ``You have to believe that we had some things in common in our temperaments and our philosophy, because there was a lot of agreement in our way of seeing things.  Naturally kind, he (Aim\'e Cotton) dedicated himself without hesitation, enduring every ordeal, to the cause he believed to be just.  On the other side, he declared himself against what he reproved with courage and clearness.  This is probably why he was detained by the Germans in 1941 at Fresnes, where the two of us were together without knowing it.  I will never forget our trip back together after leaving the prison, looking like highwaymen with our shaggy heads in the metro.''

\subsection{Camille Marbo's recollections, 1968}\label{subsec:marbo}

\textit{Excerpt from the autobiography of Emile Borel's wife, Camille Marbo (\cite{marbo:1968}, pp.\ 298--304).  This passage begins as Borel and Marbo leave Saint-Affrique for Paris in the autumn of 1940.}

\subsubsection*{French original}

En automne, Emile Borel relev\'e de ses fonctions de maire, nous rentrons \`a Paris.

Inutile d'insister sur un destin partag\'e par les trois quarts de la population parisienne:  queues aux portes des boutiques pour obtenir des rutabagas ou des choux brocolis, maigres rations contre tickets, navets rempla\c{c}ant les introuvables pommes de terres.  Saint-Affrique ne comprend pas les S.O.S.\ que nous lan\c{c}ons par allusions sur les cartes-questionnnaires impos\'ees.

Au probl\`eme du ravitaillement s'ajoute celui du chauffage.  Boulevard Haussman, les radiateurs sont glac\'es, ma m\`ere et mon mari frileux.  Je br\^ule une banale armoire \`a lingerie, une coiffeuse, une table et des chaise, empil\'ees au septi\`eme dans ma seconde chambre de bonne, inoccup\'ee.

Arrestation de Paul Langevin, incarc\'er\'e \`a la Sant\'e, auquel j'\'ecris et qui re\c{c}oit ma lettre, alors que, l'ann\'ee suivante, on n'autorisera auncune correspondance avec les intern\'es politiques.  Nous allons le voir apr\`es sa lib\'eration.  

D\'esireux de se rapprocher d'un centre intellectuel, Emile Borel d\'ecide de quitter le boulevard Haussmann pour le 4, rue Froidevaux, o\`u nous poss\'edons un petit appartement, dans l'immeuble achet\'e en commun avec des amis.  Ma m\`ere s'en va pr\`es de Pierre, rue Notre-Dame-des-Champs.

Au d\'ebut de l'\'et\'e 1941, nous pouvons nous rendre au sud de la Loire.  Pour peu de temps, Emile borel pensant alors que son devoir le rappelle \`a Paris, o\`u il a pris des contacts.  Au cours de ce voyage, je vois pour la derni\`ere fois Jean-Robert, en service \`a Toulon, et Jean Perrin \`a Lyon, o\`u il s'est retir\'e apr\`es avoir \'et\'e au Maroc.

Jean Perrin h\'esite \`a rejoindre son fils Francis, en mission scientifique \`a New York.  Nous passons un apr\`es-midi, seuls ensemble, au parc de la T\^ete-d'Or.  L\`a, en face d'une pelouse cern\'ee de tuilipes, il me fait des adieux que je pressens d\'efinitifs.  Ces heures compteront parmi les plus \'emouvants de ma vie.

En octobre, cette m\^eme ann\'ee 1941, mon mari est arr\^et\'e.  

Vers deux heures, une femme de m\'enage engag\'ee la veille, m'annonce, effar\'ee:

``Il y a, au salon, un Allemand en uniforme.''

C'est un officier, gourm\'e, qui \'echange avec mon mari des propos oiseux.  Des coups \'ebranlent notre porte.  Quatre soldats et un feldwebel apparaissent, derri\`ere lesquels p\'en\`etre ma m\`ere, tout interdite, venue pour sa visite quotidienne.  L'officier change de ton:

``Je dois faire perquisition.  Veuillez tous vous asseoir.''

Mon mari, ma m\`ere, la servante qu'un soldat ram\`ne, et moi, nous restons quatres heures, sans bouger, sous la surveillance du feldwebel, qui fume sans arr\^et, \`a califourchon sur une chaise, ordonnant brutalement de se taire, successivement, \`a la femme de m\'enage qui dit timidement:  ``Le gaz br\^ule \`a la cuisine''', puis \`a ma m\`ere, \`a bout de nerfs, qui constate:  ``Il pleut \`a verse\dots''

Nous entendons le fracas des bottes et le bruit de tiroirs jet\'es \`a terre. L'appartement est pass\'e au peigne fin. Non seulement l'armoire \`a linge vid\'ee, les draps, nappes, etc., d\'epli\'es et amoncel\'es en d\'esordre sur le plancher, comme tous les papiers, dossiers, manuscrits extraits des classeurs et des meubles, mais aussi le coffre \`a charbon renvers\'e dont le contenu s'\'eparpillera dans la cuisine.

A sept heures l'officier repara\^it:

``Monsieur le professeur, je vous prie de me suivre.''

-- Pourquoi? O\`u m'emmenez vous?

-- Cela ne me concerne pas. Que madame vous pr\'epare quelques effets pour la nuit.

-- Puis-je lui servir \`a manger?

-- Non. C'est notre affaire.''

Emile Borel m'embrasse et s'en va.  Suivi par ma m\`ere, boulevers\'ee, que j'envoie d\^iner chez Pierre, et par la femme de m\'enage, que je pense ne jamais revoir apr\`es ces funestes d\'ebuts.

Pendant plus de quinze jours, aucune nouvelle de mon mari, Pierre m'escorte pour quelques d\'emarches vaines. Paul, son fils, qui parle parfaitement allemand, court avec moi les t\'en\'ebreux bureaux de la Kommandantur, install\'es dans des h\^otels, et toutes les prisons, o\`u nous pr\'esentons un colis, refus\'e apr\`es consultation du registre d'\'ecrou:

``Pas de Borel ici.''

Enfin, gr\^ace \`a Paul, qui est retourn\'e \`a Fresnes \`a bicyclette, comme on y accepte le paquet qu'il prom\`ene, nous sommes fix\'es sur le lieu d'incarc\'eration et autoris\'es \`a apporter, une fois la semaine, du linge.  

Nous joignons au premier colis un jeu d'\'echecs.  L'officier de garde commence par le refuser. Paul explique que l'on peut s'exercer aux \'echecs sans partenaire.

``Nein!''

Un tr\`es jeune caporal, qui seconde l'officier, se risque \`a dire timidement que son p\`ere joue aux \'echecs tout seul. La bo\^ite est ouverte, les pi\`eces examin\'ees une par une, l\'echiquier de carton consid\'er\'e de pr\`es. Le jeu accept\'e. Emile Borel dira combien il lui a \'et\'e pr\'ecieux pendant les deux mois de cellule durant lesquels il n'eut aucun contact humain, sauf le matin, par le guichet au travers duquel on lui passait, pour vingt-quatre heures, une assiette de soupe, un bloc de graisse, un morceau de pain et une tranche de citron ``pour les vitamines r\'eglementaires'', et la visite hebdomadaire o\`u, sans un mot, on lui jetait le contenu de la valise que j'avais apport\'ee la veille.

Il n'eut pas droit au barbier, ne fit aucune promenade dans la cour. Dans l'obscurit\'e compl\`ete, de quatre ou cinq heures du soir jusqu'au matin, il se livrait \`a des r\'eflexions math\'ematiques et \`a ses souvenirs, couch\'e sur sa paillasse, ayant enfil\'e les uns par dessus les autres tous les habits que je lui envoyais, pour ne pas grelotter. Le jour, il marchait de long en large, \'evaluant les kilom\`etres parcourus, et s'interrompant pour une partie d'\'echecs.  Sur le couvercle en sapin de la bo\^ite d'\'echecs, \`a l'aide d'un bout de crayon \'echapp\'e \`a l fouille dans la poche de son gilet, il tra\c{c}ait, chaque matin, une petite barre, pour compter ses jours de captivit\'e.  En fait, il en inscrivit un de moins.  Je conserve la bo\^ite et son couvercle \`a Cornus.

D'autres savants ayant \'et\'e arr\^et\'es le m\^eme jour que lui , nous pensons que l'Acad\'emie des sciences, dont ils font tous partie, peut intervenir pour demander qu'ils soient humainement trait\'es, vu leur \^age, et rapidement traduits devant un tribunal.

Il nous est r\'epondu que c'est impossible. Il s'agirait d'une action politique, risquant de faire dissoudre l'acad\'emie. Le duc Maurice de Broglie, Elie Cartan, Paul Montel, le professeur Jolibois soutiennent sans succ\`es le principe d'une d\'emarche. Le ministre lui-m\^eme, que je vais voir, tout en m'assurant de sa profonde sympathie, m'oppose un refus tr\`es net: 

``H\'elas! j'ai les mains li\'ees\dots''

\textit{Ces messieurs} affirment que ``nos savants'' sont \`a l'infirmerie de Fresnes, bien couch\'es, pourvus de bonnes couvertures, convenablement nourris et que les gardiens leur t\'emoignent des \'egards. Or, les cellules de Fresnes pourvues d'installations sanitaires rudimentaires dispensent les gardiens de corv\'ees de vidange. Du reste, \`a la fa{\c c}on grossi\`ere dont ils traitent les m\`eres et les \'epouses apportant des colis, il est difficile de les imaginer simplement polis.

A pat la conduite affectueuse et fid\`ele de nos amis, deux choses me touch\`erent surtout durant cette \'epreuve.

D'abord la d\'ecision imm\'ediate de l'Acad\'emie royale de Su\`ede de nommer Emile Borel membere correspondant, en lui conf\'erant le droit de proposer des prix Nobel.  Ensuite la conduite de ma jeune femme de m\'enage qui, m'ayant pr\'evenue qu'elle ne pourrait entrer \`a mon service, r\'eguli\`erement, qu'une semaine plus tard, vint sonner \`a ma porte d\`es l'aube, le lendemain de l'arrestation, disant tranquillement:  ``J'ai envoy\'e promener l'autre dame.  Il ya du travail ici.''

Cette pr\'esence attentive d'une ``salari\'ee'' que je ne connaissais pas auparavant me fut d'un extr\^eme r\'econfort.  Je n'oublierai jamais cette jolie Espagnole, nomm\'ee Julienne.

Emile Borel me fut renvoy\'e \`a l'improviste. Maigri de sept kilos, h\^ave, le visage mang\'e de barbe. Il avait souffert du froid plus que de la faim. Et il avait d\'e affronter, sa valise \`a la main, les courants d'air du quai, en attendant le m\'etro qui le ram\`enerait. 

Le lendemain, avec $40^{\circ}$ de fi\`evre, une pneumonie double se d\'eclarait, dont il se remit tr\`es mal, n'\'etant pas aliment\'e comme il l'aurait fallu.

Pendant sa captivit\'e, les femmes d'autres prisonniers, parmi lesquelles Marie de Pange, Eug\'enie Cotton et moi nous formions un \^ilot d'angoisses au milieu d'une population r\'esign\'ee. Connaissant des alternatives de terreur lorsqu'on nous annon{\c c}ait un transfert en Allemagne, une fusillade \`a Vincennes. Ou lorsqu'on nous coupait, \`a la prison, la r\'eception de nos paquets, en laissant entendre que nos prisonniers \'etaient exp\'edi\'es on ne sait o\`u. Par contre, quel soulagement lorsque nous rapportions la valise de la livraison pr\'ec\'edente avec leur linge \`a laver, ce qui n'arriva que deux fois; la premi\`ere, mon mari avait inscrit, \`a l'aide du miraculeux bout de crayon oubli\'e \`a la fouille: O.K., tout va bien, sur la ceinture d'un cale{\c c}on.

Eug\'enie Cotton, dont le mari partageait le sort d'Emile Borel, et moi, nous nous agitions comme deux bourdons affol\'es. On nous exp\'edia chez un avocat inefficace. J'allai trouver l'aum\^onier de la rue Lhomond, l'abb\'e Stock qui, faute d'avoir obtenu le permis de communiquer, ne put visiter mon mari en lui portant une Bible, comme il me l'avait promis.

Lors d'une p\'eriode o\`u Fresnes refusait les paquets et o\`u l'on m'affirmait que nos prisonniers \'etaient transf\'er\'es ailleurs, d'apr\`es un nom donn\'e par Mme de Pange, je me risque avenue Foch, dans un coquet h\^otel, o\`u je ne sais pas que la Gestapo vient de s'installer, rempla{\c c}ant d'autres bureux allemands. L\`a, ayant demand\'e une certaine \textit{Fr\"aulein}, qui avait su renseigner Mme de Pange sur le lieu o\`u se trouvait son mari, me voici pouss\'ee de pi\`ece en pi\`ece, enferm\'ee \`a cl\'e pendant plus d'une heure, saisie par deux soldats qui me plantent dans un vestibule et, par une porte ouverte, parlent \`a un officier. J'ouvre une autre porte, esp\'erant trouver l'air libre et je tombe dans une longue salle \'etroite o\`u des soldats, assis derri\`ere des tr\'eteaux garnis de bougies cachettent des enveloppes jaunes avec de la cire noire. J'\'eprouve une frayeur \`a ce spectacle. Un feldwebel m'empoigne, me tra\^ine, me jette sur le trottoir de l'avenue Foch. Je n'ai jamais bien compris \`a quoi tout cela correspondait.

On m'avertit, par une voie myst\'erieuse, que si je vais, dans un local proche de la Sorbonne, signer l'engagement que ni moi ni lui nous ne nous opposerons \`a l'occupant, Emile Borel sera rel\^ach\'e dans les vingt-quatre heures. Je ne vais pas au rendez-vous. J'explique \`a Paul Montel, qui vient me voir, juste au moment o\`u me quitte l'interm\'ediaire officieux, d'abord que mon mari ne me pardonnerait pas de donner ma signature \`a cette renonciation, puis que je n'y ai aucun m\'erite, ne croyant pas \`a la parole des nazis.

``En somme, conclut Montel, vous ne voulez pas jouer les Tosca.''

\subsubsection*{English translation}

In the autumn, Emile Borel having been relieved of his functions as mayor, we returned to Paris.

There is no use dwelling on the destiny we shared with three-fourths of the Parisian population:  lines at the doors of stores to get rutabagas or winter cauliflower, tickets traded for meager rations, turnips replacing impossible-to-find potatoes.  Saint-Affrique did not understand the SOSs we tried to put between the lines in the questionnaire-like postcards we were allowed to send.  

On top of the problem of finding food, there was the problem of heating.  In our apartment at Boulevard Hausmann, the radiators were frozen, my mother and my husband were freezing.  I burned furniture that we had piled in the unoccupied second maid's room on the seventh floor:  an ordinary clothes chest, a dressing table, a table and chairs.

Then there was the arrest of Paul Langevin, thrown into the Sant\'e prison.  I wrote to him and he received my letter, whereas the following year no correspondence with political prisoners was allowed.  We went to see Langevin after he was released.

Wanting to get close to an intellectual center, Emile Borel decided to leave Boulevard Hausmann for 4 rue Froidevaux, where we had a small apartment in a building that we had bought together with friends.  My mother moved close to Pierre,%
\footnote{Pierre Appell (1887--1957), Marbo's brother, was in the French parliament from 1928 to 1936 and then went into business.}
rue Notre-Dame-des-Champs.

At the beginning of summer 1941, we were able to go south of the Loire.  But only for a short time, because Emile Borel thought it was his duty to return to Paris, where he had made contacts.  During that trip, I saw for the last time Jean-Robert,%
\footnote{Jean-Robert Appell, Pierre's son, was arrested while trying to reach England and deported to the Neuengamme concentration camp in Germany.  He died on 3 May 1945 in the tragedy at Lubeck, where ships meant to transport the camp's prisoners to their deaths were sunk by the British air force.}
who was in service at Toulon, and Jean Perrin in Lyon, where he had come after having been in Morocco.

Jean Perrin was hesitating about whether to join his son Francis in New York, where he was on a scientific mission.  We passed an afternoon alone together, at the T\^ete-d'Or park, across from a lawn surrounded by tulips, and he gave me farewells that seemed final to me.  Those hours are among the most moving of my life.%
\footnote{The Nobel prize winning physicist Jean Perrin was one of Emile Borel's closest allies.  He died in New York on April 1942.  The T\^ete-d'Or is a park in Lyon.  Borel and Marbo could have gone from Paris to south of the Loire without crossing into the occupied zone, but Toulon and Lyon were in the occupied zone.} 

In October of that same year, 1941, my husband was arrested.

Around two o'clock, the housekeeper I had hired the day before announced, with a fright:

``There is a German in uniform in the living room.''

It was a pretentious officer exchanging meaningless words with my husband.  Then our door was shaken by knocking.  Four soldiers and a sergeant appeared, and behind them my mother, disconcerted, who had come for her daily visit.  The officer's tone changed:

``I must conduct a search.  All of you please take a seat.''

The four of us, my husband, my mother, the housekeeper brought back by a soldier, and myself, waited without budging for four hours, under the surveillance of the sergeant, who smoked constantly, sitting in a chair with his legs crossed, brutally ordering to shut up first the housekeeper, who timidly said, ``The gas is on in the kitchen,'' and then my mother, who noted, at the end of her nerves, ``It's raining like\dots''.

We heard the noise of boots and of drawers being thrown on the floor.  They put the apartment through a fine-toothed comb.  Not only was the linen closet emptied, the sheets, covers, etc.\ unfolded and thrown in a mess on the floor, along with all the papers, files, and manuscripts from the desks and cabinets, but even the coal bucket was knocked over, its contents strewn across the kitchen.

The officer reappeared at seven o'clock:

``Mr.\ Professor, please follow me.''

-- Why?  Where are you taking me?

-- That is none of my business.  Have your wife pack you a bag for the night.

-- May I give him something to eat?

-- No.  That's our job.''

Emile Borel hugged me and left.  Followed by my distressed mother, whom I sent to Pierre's for dinner, and my housekeeper, whom I supposed I would never see again after this horrible beginning.

There was no news about my husband for more than two weeks, as Pierre accompanied me on several futile missions.  His son Paul, who spoke German perfectly, went with me into the numerous offices of the German military administration, located in hotels, and into all the prisons, where we offered a package, always refused after consultation of a register of those being held:

``No Borel here.''

Finally, thanks to Paul, who bicycled to Fresnes with a package they accepted, we learned the place of incarceration and were authorized to bring laundry once a week.  

We added a chess set to the first package.  At first the officer guarding the entry refused to take it.  Paul explained that you can play chess by yourself.

``Nein.''

A very young corporal, who was helping the officer, dared to say timidly that his father played chess by himself.  The box was opened, the pieces examined one by one, the cardboard chess board examined closely.  The chess set was accepted.  Emile Borel would say how precious it was to him during his two months in a cell with no human contact except through the grill each morning, when he received his daily ration---a bowl of soup, a square of fat, a piece of bread, and a slice of lemon ``for the statutory vitamins''---and the weekly visit where they silently threw in the contents of the bag I had brought the day before.  

He had no right to a barber.  He never took a walk in the courtyard.  In complete darkness, from four or five o'clock in the evening until morning, he occupied himself with mathematical reflections and his memories, lying on his straw mat, having put on all the clothing I had sent him,  one garment over the other, to avoid shivering.  During the day, he walked back and forth, calculating how many kilometers he had gone and pausing for a game of chess.  On the fir lid of the chess set, using a stub of a pencil that had been missed when the pocket of his jacket was searched, he made a little mark every morning, to count the days of his imprisonment.  In fact, he had made one mark too few.  I still have the chess set and its lid at Cornus.

Other scientists having been arrested the same day, we thought that the \textit{Acad\'emie des Sciences}, where they were members, could intervene to ask that they be treated humanely, in consideration of their age, and quickly brought before a court.  

We were told that this was impossible.  It would be a political act, risking the abolition of the Academy.  Duke Maurice de Broglie, Elie Cartan, Paul Montel, and Professor Jolibois supported taking action, without success.  I went to see the minister, who expressed his profound sympathy but made his refusal very clear:  ``Alas, my hands are tied\dots''

\textit{These gentlemen} declared that ``our scientists'' were well cared for at the Fresnes infirmary, with good bedding and proper food, and that the guards treated them with respect.  In fact, the cells at Fresnes did have rudimentary sanitary facilities, which relieved the guards from taking out the sewage.  But as for the rest, the rudeness with which they treated the mothers and wives who brought packages makes it difficult to imagine that they were even polite.

Aside from the affectionate and faithful conduct of our friends, two things especially touched me during this difficulty.

First, the immediate decision of the Swedish Royal Academy to name Emile Borel a corresponding member, while giving him the right to make nominations for the Nobel prize.  Secondly, the conduct of my young housekeeper.  Having earlier warned me that she could not work for me regularly for another week, she came to ring the bell at my door the very next morning, saying calmly, ``I dropped the work for the other woman.  There is work to do here.''

This attentive presence of an ``employee'' whom I had not known before was a great comfort to me.  I will never forget that pretty young Spaniard, named Julienne.

Emile Borel was returned to me unannounced.  Seven kilos thinner, haggard, his face obscured by his beard.  He had suffered more from the cold than from hunger.  And he had had to endure the wind, suitcase in hand, while waiting for the metro that brought him back home.

The next morning, he came down with double pneumonia, with a 40-degree fever.  He had a hard time recovering, not having been nourished as he should have been.

During his imprisonment, the wives of the other prisoners, including Madame de Pange and Eug\'enie Cotton, and I formed an island of anxiety in the midst of a resigned population.  We experienced a moment of terror whenever they announced a transfer of prisoners to Germany or an execution at Vincennes.  They would stop accepting our packages at the prison, leaving us to understand that our prisoners had been sent to some unknown destination.  What relief, on the other hand, when they gave us back the bag we had brought the time before, with the laundry to wash.  This happened twice.  The first time, my husband had written ``OK, everything is going well'' on the waist of a pair of underpants, with the help of the miraculous piece of pencil overlooked when he was searched.

Eug\'enie Cotton, whose husband shared Emile Borel's fate, and I buzzed around like two crazy bumblebees.  We were sent to see an ineffective lawyer.  I went to find the Lhomond Street chaplain, Father Stock, who was not able to keep his promise to visit my husband and take him a bible, because he had not acquired a driver's license.

During a period when Fresnes refused packages and they told me that our prisoners had been transferred somewhere else, following up on a name given me by Madame de Pange, I dared to go to Foch Avenue, to a stylish hotel where, unbeknownst to me, the Gestapo had just moved in, replacing some other German offices.  There, having asked for a certain \textit{Fr\"aulein}, who had been able to tell Madame de Pange right away where her husband was, I was sent from one room to another, locked in for an hour, then seized by two soldiers, who set me in a corridor and talked with an officer through an open door.  I opened another door, hoping to find the open air, and I found myself in a narrow room, where soldiers behind trestles with candles on them were sealing yellow envelopes with black wax.  It was a frightening spectacle.  A sergeant grabbed me, drug me away, and threw me onto the sidewalk of Foch Avenue.  I never quite understood what it was all about.

Once I got word, by a mysterious route, that if I went to a place near the Sorbonne and signed a guarantee that neither he nor I would oppose the occupier, Emile Borel would be released in twenty-four hours.  I did not go to the appointment.  I explained to Paul Montel, who came to see me just as the unofficial intermediary left, that in the first place my husband would never forgive me for signing such a renunciation, and then that I saw no point in it, as I did not trust the Nazis' word.  

``In short,'' Montel said, ``you do not want to play Tosca.''

\pagebreak
\tableofcontents

\end{document}